\documentclass[11pt]{amsart}

\usepackage{a4wide}
\usepackage{amssymb}
\usepackage{amsthm}
\usepackage{mathrsfs,mathtools}
\usepackage{enumerate}
\usepackage{mathscinet}
\usepackage{titletoc}
\usepackage[colorlinks=true, urlcolor=blue, linkcolor=blue, citecolor=blue]{hyperref}
\usepackage{kotex}
\usepackage{xcolor}
\usepackage{enumitem}
\usepackage[nameinlink]{cleveref}

\theoremstyle{plain}
\newtheorem{theorem}{Theorem}[section]
\newtheorem{lemma}[theorem]{Lemma}
\newtheorem{corollary}[theorem]{Corollary}

\theoremstyle{definition}
\newtheorem{definition}[theorem]{Definition}
\newtheorem{remark}[theorem]{Remark}

\newtheorem{notation}[theorem]{Notation}

\numberwithin{equation}{section}

\makeatletter
\@namedef{subjclassname@2020}{\textup{2020} Mathematics Subject Classification}
\makeatother

\title[Parabolic $p$-Laplace type equations]{Boundary H\"older gradient estimates for parabolic $p$-Laplace type equations}

\author{Se-Chan Lee}
\address{School of Mathematics, Korea Institute for Advanced Study, Seoul 02455, Republic of Korea}
\email{sechan@kias.re.kr}

\author{Yuanyuan Lian}
\address{Departamento de An\'{a}lisis Matem\'{a}tico,
Instituto de Matem\'{a}ticas IMAG, Universidad de Granada}
\email{lianyuanyuan.hthk@gmail.com; yuanyuanlian@correo.ugr.es}

\author{Hyungsung Yun}
\address{School of Mathematics, Korea Institute for Advanced Study, Seoul 02455, Republic of Korea}
\email{hyungsung@kias.re.kr}

\author{Kai Zhang}
\address{Departamento de Geometr\'{i}a y Topolog\'{i}a,
Instituto de Matem\'{a}ticas IMAG, Universidad de Granada}
\email{zhangkaizfz@gmail.com; zhangkai@ugr.es}


\thanks{Se-Chan Lee has been supported by the KIAS Individual Grant (No. MG099001) at Korea Institute for Advanced Study. Yuanyuan Lian has been supported by PID2020-
117868GB-I00 and PID2023-150727NB-I00 funded by MCIN/AEI. Hyungsung Yun has been supported by the KIAS Individual Grant (No. MG097801) at Korea Institute for Advanced Study. Kai Zhang has been supported by the Project PID2020-118137GB-I00 funded by MCIN/AEI}

\begin{document}

\begin{abstract}
In this paper, we study the boundary regularity for viscosity solutions of parabolic $p$-Laplace type equations. In particular, we obtain the boundary pointwise $C^{1,\alpha}$ regularity and global $C^{1,\alpha}$ regularity.
\end{abstract}

\maketitle

%
%

\tableofcontents

\section{Introduction}
In this paper, we study the boundary regularity for viscosity solutions of the following parabolic $p$-Laplace type singular/degenerate equation:
\begin{equation}\label{e1.1}
  \left\{
    \begin{aligned}
    Pu &=f && \mbox{in }\Omega\\
    u&=g&& \mbox{on }\partial_p \Omega,
  \end{aligned}
  \right.
\end{equation}
where $\Omega \subset \mathbb{R}^{n+1}$ is a bounded domain and the $p$-Laplace type operator $P$ is defined as
\begin{equation*}
    Pu\coloneqq u_t-|Du|^{\gamma}\left(\delta^{ij} + (p-2)\frac{u_iu_j}{|D u|^{2}}\right)u_{ij} = u_t -|Du|^{\gamma+2-p} \Delta_p u.
\end{equation*}
Throughout this paper, we always assume that
\begin{equation*}
    -1<\gamma<\infty \quad\text{and} \quad 1<p<\infty.
\end{equation*}

The first equation in \eqref{e1.1} includes two important examples. If $\gamma=p-2$, then it becomes the parabolic $p$-Laplace equation, which exhibits degeneracy when $p>2$ and singularity when $1<p<2$. The parabolic $p$-Laplace equation has been widely studied in the context of the regularity theory; we refer to the monograph \cite{MR1230384} by DiBenedetto. Just to name a few, the interior H\"older gradient estimate was developed by DiBenedetto and Friedman \cite{MR783531}, while the boundary H\"older gradient estimate under zero boundary conditions was established by Chen and DiBenedetto \cite{CD89}. Moreover, Lieberman proved the boundary H\"older gradient estimate for weak solutions with conormal boundary condition in \cite{Lie90} and with general Dirichlet boundary condition in \cite{MR1207530}. We would like to point out that the strategies in \cite{CD89, Lie90, MR1207530} strongly rely on the divergence structure of equations, particularly through integration by parts with suitably chosen test functions.

Another important example is the parabolic normalized $p$-Laplace equation, which corresponds to the case $\gamma=0$. Since the normalized $p$-Laplace equation does not admit an associated energy-like quantity, the notion of viscosity solutions (instead of weak solutions) becomes essential. Within the nondivergence framework, Jin and Silvestre \cite{JS17} established the interior H\"older gradient estimate for parabolic normalized $p$-Laplace equation with $f\equiv 0$. Attouchi and Parviainen \cite{AP18} extended the interior regularity result for general nonhomogeneous term $f$. The boundary regularity is rather straightforward, since one can apply the boundary H\"older gradient estimate for functions in the solution class; see \cite{LY24} for details. We also refer the reader to \cite{AS22, DPZZ20, FZ21, FZ23, FPS25, LLYZ25, PV20} for further results on interior and boundary regularity in various settings.

Recently, Imbert, Jin and Silvestre \cite{MR3918407} proved the interior H\"{o}lder gradient estimate for viscosity solutions of \eqref{e1.1} (with $f\equiv 0$) in a unified manner for the full ranges of $\gamma$ and $p$. Later, it was extended to equations with general nonhomogeneous term $f$ (see \cite[Theorem 1.1]{MR4122677} for $0<\gamma<+\infty$, \cite[Theorem 1.1]{MR3804259} for $\gamma=0$ and \cite[Theorem 1.1]{MR4128334} for $-1<\gamma<0$). Moreover, the interior H\"{o}lder gradient estimate has been extended to fully nonlinear equations by Lee, together with the first and third authors \cite{MR4676644}.

The goal of the present paper is to develop the boundary counterpart of \cite{MR3918407} and \cite{MR4122677,MR3804259,MR4128334}, i.e., the boundary H\"{o}lder gradient estimate for \eqref{e1.1}. To the best of our knowledge, no suitable technique has been developed for boundary $C^{1,\alpha}$ regularity for solutions of \eqref{e1.1}. This work proposes an approach tailored to address boundary $C^{1,\alpha}$ regularity for such structures. Our approach is purely non-variational in the sense that we do not make use of integral estimates. Compared to the techniques used in \cite{CD89,Lie90,MR1207530}, our technique is more flexible and can be applied to more complicated problems (e.g. fully nonlinear degenerate/singular parabolic equations). Moreover, we aim to impose sharper conditions on $g$ and $\partial \Omega$ that are consistent with the boundary regularity of solutions. Instead of flattening the curved boundary, we apply the perturbation argument combined with compactness method; see the recent results in \cite{MR4088470, LZ_Parabolic}. To the best of our knowledge, the boundary pointwise $C^{1,\alpha}$ regularity is new even for the standard elliptic $p$-Laplace equations (see \Cref{co1.0}).

Let us first introduce the definition of the pointwise $C^{1,\alpha}$ smoothness for a function.
\begin{definition}\label{d-f}
Let $\Omega\subset \mathbb{R}^{n+1}$ be a bounded set (may not be a domain) and $f:\Omega\rightarrow \mathbb{R}$ be a function. Given $0<\alpha\leq 1$ and $\gamma>-1$ satisfying $\alpha(1+\gamma)\leq 1$, we say that $f$ is $C^{1,\alpha}_{\gamma}$ at $(x_0,t_0)\in \Omega$ (denoted by $f\in C^{1,\alpha}_{\gamma}(x_0,t_0)$) if there exist positive constants $K$, $r_0$ and a linear polynomial $L(x)$ (independent of $t$) such that
\begin{equation}\label{m-holder}
  |f(x,t)-L(x)|\leq K(|x-x_0|^{1+\alpha}+|t-t_0|^{\frac{1+\alpha}{2-\alpha\gamma}}) \quad\text{for all }(x,t)\in \Omega\cap Q_{r_0}(x_0,t_0).
\end{equation}
Then define
\begin{equation*}
    Df(x_0,t_0)= DL, \quad \|f\|_{C^{1}(x_0,t_0)}=|L(x_0)|+|DL|,
\end{equation*}
and
\begin{equation*}
\begin{aligned}
    [f]_{C^{1,\alpha}_{\gamma}(x_0,t_0)}&= \min \left\{K \mid \eqref{m-holder} \mbox{ holds with }L\mbox{ and }K\right\}, \\
    \|f\|_{C^{1,\alpha}_{\gamma}(x_0,t_0)}&= \|f\|_{C^{1}(x_0,t_0)}+[f]_{C^{1,\alpha}_{\gamma}(x_0,t_0)}.
\end{aligned}
\end{equation*}
If $f\in C^{1,\alpha}_{\gamma}(x,t)$ for any $(x,t)\in \Omega$ with the same $r_0$ and
\begin{equation*}
  \|f\|_{C^{1,\alpha}_{\gamma}(\overline{\Omega})}\coloneqq \sup_{(x,t)\in \Omega} \|f\|_{C^{1, \alpha}_{\gamma}(x,t)}<+\infty,
\end{equation*}
we say that $f\in C^{1,\alpha}_{\gamma}(\overline{\Omega})$.
\end{definition}
\begin{remark}\label{re1.0}
If $\gamma=0$, we write $f\in C^{1,\alpha}$ instead of $f\in C^{1,\alpha}_{0}$ for simplicity. Note that in this case, \Cref{d-f} coincides with the standard definition of the pointwise $C^{1,\alpha}$ smoothness. In addition, it is clear from the definition that
\begin{equation*}
    C^{1,\alpha}_{\gamma_1}\subset C^{1,\alpha}_{\gamma_2} \quad\mbox{if } \gamma_1\geq \gamma_2 \qquad  \mbox{and} \qquad C^{1,\alpha_1}_{\gamma}\subset C^{1,\alpha_2}_{\gamma} \quad\mbox{if } \alpha_1\geq \alpha_2.
\end{equation*}
\end{remark}

\begin{remark}\label{re15.1}
The condition $\alpha(1+\gamma)\leq 1$ is to guarantee that $(1+\alpha)/(2-\alpha\gamma)\leq 1$.
\end{remark}

\begin{remark}\label{re15.4}
If $\Omega$ is the boundary of some domain, a linear polynomial $L$ in \eqref{m-holder} may not be unique (see \cite[Remark 1.3 and Remark 1.4]{MR4644419} for the explanation). However, this definition still works for our purpose to find a polynomial $L$ that approximates $f$ in $C_{\gamma}^{1,\alpha}$ sense.
\end{remark}

We next provide the definition of the $C^{1,\alpha}$ domain suggested in \cite[Definition 1.4]{LZ_Parabolic}.
\begin{definition}\label{d-domain}
Let $\Omega$ be a bounded domain, $\Gamma\subset \partial_p \Omega$ be relatively open and $(x_0,t_0)\in \Gamma$. Given $0<\alpha\leq 1$ and $\gamma>-1$ satisfying $\alpha(1+\gamma)\leq 1$, we say that $\Gamma$ is $C^{1,\alpha}_{\gamma}$ at $(x_0,t_0)$ (denoted by $\Gamma\in C^{1,\alpha}_{\gamma}(x_0,t_0)$) if there exist constants $K>0$, $0<r_0\leq 1$ and a new coordinate system $\{x_1,\cdots,x_n,t \}$ (by rotating and translating with respect to $x$ and only translating with $t$) such that $(x_0,t_0)=(0,0)$ in this coordinate system,
\begin{equation}\label{e-re}
    Q_{r_0} \cap \big\{(x',x_n,t) \mid x_n>K(|x-x_0|^{1+\alpha}+|t-t_0|^{\frac{1+\alpha}{2-\alpha\gamma}})\big\} \subset Q_{r_0}\cap \Omega
\end{equation}
and
\begin{equation}\label{e-re2}
    Q_{r_0} \cap \big\{(x',x_n,t) \mid x_n<-K(|x-x_0|^{1+\alpha}+|t-t_0|^{\frac{1+\alpha}{2-\alpha\gamma}})\big\} \subset Q_{r_0}\cap \Omega^c.
\end{equation}
Then define
\begin{equation*}
    [\Gamma]_{C^{1,\alpha}_{\gamma}(x_0,t_0)} =\min \left\{K \mid \eqref{e-re} \mbox{ and } \eqref{e-re2}
\mbox{ hold with } K\right\}.
\end{equation*}
If $\Gamma\in C^{1, \alpha}_{\gamma}(x,t)$ for any $(x,t)\in \Gamma$ with the same $r_0$ and
\begin{equation*}
  \|\Gamma\|_{C^{1,\alpha}_{\gamma}}\coloneqq \sup_{(x,t)\in \Gamma}~[\Gamma]_{C^{1, \alpha}_{\gamma}(x,t)}<+\infty,
\end{equation*}
we say that $\Gamma\in C^{1,\alpha}_{\gamma}$. If $\Gamma'\in C^{1,\alpha}_{\gamma}$ for any $\Gamma'\subset \joinrel \subset\Gamma$, we denote $\Gamma\in C^{1,\alpha}_{\gamma, \mathrm{loc}}$.

If \eqref{e-re} and \eqref{e-re2} are replaced by
\begin{equation}\label{e-re3}
    Q_{r_0} \cap \big\{(x',x_n,t) \mid x_n>K\big\} \subset Q_{r_0}\cap \Omega
\end{equation}
and
\begin{equation}\label{e-re4}
    Q_{r_0} \cap \big\{(x',x_n,t) \mid x_n<-K\big\} \subset Q_{r_0}\cap \Omega^c,
\end{equation}
we define
\begin{equation*}
   \underset{Q_{r_0}}{\mathrm{osc}} \, \partial_p\Omega =\min \left\{K \mid \eqref{e-re3} \mbox{ and } \eqref{e-re4}
\mbox{ hold with } K\right\}.
\end{equation*}
\end{definition}

\begin{remark}\label{re1.2}
For studying parabolic equations, a boundary $\partial_p \Omega$ is divided into two classes: lateral boundary and bottom boundary (see \cite[P. 29--30]{MR1135923} for the precise definitions). From \Cref{d-domain}, we know that if $\Gamma\in C^{1,\alpha}_{\gamma}(x_0,t_0)$, then $(x_0,t_0)$ must belong to the lateral boundary rather than the bottom boundary. Note that the regularity on the lateral boundary and on the bottom boundary are quite different (compare \cite[Section 2]{MR1139064} with \cite{MR1151267}).

Let $U\subset \mathbb{R}^n$ be a bounded domain, $0\in \partial U$ and $\Omega\coloneqq U\times (-1,0]$. If $\partial U\in C^{1,\alpha}(0)$ in the usual sense (see \cite[Definition 1.2]{MR4088470}), then $\partial_p \Omega \in C^{1,\alpha}_{\gamma}(0,0)$ for any $\gamma>-1$.
\end{remark}

The following are our main results. We always assume that $(0,0)\in\partial_p \Omega$ and study the pointwise regularity at $(0,0)$. In addition, if we use \Cref{d-f} or \Cref{d-domain} at $(0,0)$, then we always assume that $r_0=1$. In this paper, a constant $C>0$ is called universal if it depends only on $n$, $p$ and $\gamma$.
\begin{theorem}\label{th1.1}
Let $u\in C(\overline{\Omega\cap Q_1})$ be a viscosity solution of
\begin{equation}\label{e11.1}
\left\{\begin{aligned}
    Pu&=f&& \mbox{in } \Omega\cap Q_1\\
    u&=g&& \mbox{on } \partial_p \Omega\cap Q_1.
\end{aligned}\right.
\end{equation}
Suppose that
\begin{equation*}
	f\in C(\Omega\cap Q_1)\cap L^{\infty}(\Omega\cap Q_1)
\end{equation*}
and
\begin{equation*}
  \left\{
  \begin{aligned}
&g\in C^{1,\alpha}_{\gamma}(0,0), \quad \partial_p\Omega\cap Q_1\in C^{1,\alpha}_{\gamma}(0,0)\quad\mbox{if } \gamma> 0\\
&g\in C^{1,\alpha}(0,0), \quad \partial_p\Omega\cap Q_1\in C^{1,\alpha}(0,0)\quad\mbox{if } \gamma\leq 0\\
  \end{aligned}
  \right.
\end{equation*}
for some $\alpha \in(0,\bar{\alpha})$, where $0<\bar{\alpha}< \min\{1/2,1/2(1+\gamma)\}$ is universal (see \Cref{le3.7}).

Then
\begin{equation*}
  \left\{
  \begin{aligned}
u &\in C^{1,\alpha}(0,0)&&\mbox{if }\gamma> 0\\
u &\in C_{\gamma}^{1,\alpha}(0,0)&&\mbox{if }\gamma\leq 0,\\
  \end{aligned}
  \right.
\end{equation*}
i.e., there exists a linear polynomial $L$ (independent of $t$) such that for any $(x,t)\in \Omega\cap Q_1$,
\begin{equation}\label{e1.2}
    |u(x,t)-L(x)|\leq\left\{
  \begin{aligned}
&C (|x|^{1+\alpha}+|t|^{\frac{1+\alpha}{2}}) &&\mbox{if }\gamma> 0\\
  &C (|x|^{1+\alpha}+|t|^{\frac{1+\alpha}{2-\alpha\gamma}}) &&\mbox{if }\gamma\leq 0,
  \end{aligned}
  \right.
\end{equation}
and
\begin{equation*}
|DL|\leq  C,
\end{equation*}
where $C>0$ is a constant depending only on $n$, $p$, $\gamma$, $\alpha$, $\|u\|_{L^{\infty}(\Omega\cap Q_1)}$, $\|f\|_{L^{\infty}(\Omega\cap Q_1)}$, $\|g\|_{C^{1,\alpha}_{\gamma}(0,0)}$ (or $\|g\|_{C^{1,\alpha}(0,0)}$) and $\|\partial_p\Omega\cap Q_1\|_{C^{1,\alpha}_{\gamma}(0,0)}$ (or $\|\partial_p\Omega\cap Q_1\|_{C^{1,\alpha}(0,0)}$).
\end{theorem}
\begin{remark}\label{re15.2}
The assumption $f\in C(\Omega\cap Q_1)$ is due to the fact that we consider viscosity solutions in this paper.
\end{remark}

\begin{remark}\label{re10.1}
The universal constant $\bar{\alpha}$ originates from the model problem, i.e., the problem with a flat boundary and homogeneous boundary data (see \Cref{le3.7}). In general, $\bar{\alpha}\ll 1$; however, in some particular cases, it can be relatively large (close to $1$). The simplest case is: $\gamma=0$ and $p=2$. Then the operator $P$ reduces to the standard parabolic operator $u_t-\Delta u$, and so we can take $\bar{\alpha}=1$ in \Cref{th1.1}.

Furthermore, if $\gamma=0$ and $p$ is close to $2$, then $\bar{\alpha}$ can be chosen to be close to $1$. Precisely, for any $0<\bar{\alpha}<1$, there exists $\delta>0$ such that if $\gamma=0$ and $|p-2|\leq \delta$, then we have the boundary $C^{1,\bar\alpha}$ regularity for the model problem (see \cite[Theorem 1.4]{LY24}).
\end{remark}

\begin{remark}\label{re1.4}
Note that (cf. \Cref{re1.0}) $C^{1,\alpha}_{\gamma}$ is stronger (resp. weaker) than $C^{1,\alpha}$ if $\gamma\geq 0$ (resp. $\gamma\leq 0$). Hence, in \Cref{th1.1}, we obtain a weaker regularity (e.g., $C^{1,\alpha}$ regularity for $\gamma\geq 0$) upon a stronger assumption (e.g., $C^{1,\alpha}_{\gamma}$ assumption for $\gamma\geq 0$).

The reason is the following. The regularity depends on whether the equation is degenerate (i.e., $|DL|=0$) or nondegenerate (i.e., $|DL|\neq 0$) at $(0,0)$. The scaling transformations are different for these two cases. If the equation is degenerate, we use the two-parameter scaling and we can obtain the $C^{1,\alpha}_{\gamma}$ regularity. If the equation is nondegenerate, we use the usual parabolic scaling and we can obtain the $C^{1,\alpha}$ regularity. Since we do not know a priori whether the equation is degenerate or not, we have to make a stronger assumption such that it works for both cases (see \Cref{th3.2} for the detailed proof).
\end{remark}

\begin{remark}\label{re1.8}
Up to our knowledge, the first boundary $C^{1,\alpha}$ regularity for parabolic $p$-Laplace type equations was obtained by Chen and DiBenedetto \cite{CD89}. They proved the boundary regularity for parabolic systems under the assumption:
\begin{equation*}
\gamma=p-2, \quad \max\left\{1, \frac{2n}{n+2}\right\}<p<+\infty \quad \text{and} \quad g\equiv 0.
\end{equation*}
Since $g\equiv 0$, they could use an odd reflection to reduce the model problem (i.e., the boundary value problem with a flat boundary; see \eqref{e2.0}) to the interior problem (see \cite[Proposition 3.1]{CD89}). This technique fails for a general $g$.

At almost the same time, Lieberman \cite{MR1207530} obtained the boundary $C^{1,\alpha}$ regularity for a single parabolic equation with $\gamma=p-2$, a general $g$ and the full range of $p$ (i.e., $1<p<+\infty$). His proof is involved.

Note that both \cite{CD89} and \cite{MR1207530} considered equations in divergence form since they used techniques based on integral estimates. In addition, under the $C^{1,\alpha}$ assumptions (i.e., $g\in C^{1,\alpha}(0,0)$ etc.), they could only obtain the boundary $C^{1,\tilde\alpha}$ regularity for some $\tilde{\alpha}\leq \alpha$. Instead, we obtain the $C^{1,\alpha}$ regularity in \Cref{th1.1}. Moreover, their regularity results are not pointwise since they used the technique of flattening the boundary with some transformations.
\end{remark}

Since $\gamma>-1$, we have $(1+\alpha)/(2-\alpha\gamma)>1/2$. Hence, the second estimate in \eqref{e1.2} is indeed a pointwise $C^{1,\alpha_0}$ regularity in the usual parabolic distance with $\alpha_0=\alpha(2+\gamma)/(2-\alpha\gamma)$. Therefore, we have the following boundary pointwise $C^{1,\alpha}$ regularity in a conciser (but weaker) form.
\begin{corollary}
Let $u\in C(\overline{\Omega\cap Q_1})$ be a viscosity solution of
\begin{equation*}
\left\{\begin{aligned}
    Pu&=f&& \mbox{in } \Omega\cap Q_1\\
    u&=g&& \mbox{on } \partial_p \Omega\cap Q_1.
\end{aligned}\right.
\end{equation*}
Suppose that
\begin{equation*}
 \begin{aligned}
	f\in C(\Omega\cap Q_1)\cap L^{\infty}(\Omega\cap Q_1), \quad g\in C^{1,\alpha}(0,0) \quad \text{and} \quad \partial_p\Omega\cap Q_1\in C^{1,\alpha}(0,0).
  \end{aligned}
\end{equation*}
Then $u\in C^{1,\alpha_0}(0,0)$ for some $\alpha_0\in (0,\alpha]$. That is, there exists a linear polynomial $L$ such that
\begin{equation*}
    |u(x,t)-L(x)|\leq C (|x|^{{1+\alpha_0}} +|t|^{\frac{{1+\alpha_0}}{2}})\quad \text{for all }(x,t)\in \Omega\cap Q_1
\end{equation*}
and
\begin{equation*}
|DL|\leq  C,
\end{equation*}
where $C>0$ is a constant depending only on $n$, $p$, $\gamma$, $\alpha$, $\|u\|_{L^{\infty}(\Omega\cap Q_1)}$, $\|f\|_{L^{\infty}(\Omega\cap Q_1)}$, $\|g\|_{C^{1,\alpha}(0,0)}$ and $\|\partial_p\Omega\cap Q_1\|_{C^{1,\alpha}(0,0)}$.
\end{corollary}

Since any elliptic equation can be regarded as a special parabolic equation, we have the following corollary for the classical elliptic $p$-Laplace equations.
\begin{corollary}\label{co1.0}
Let $U\subset \mathbb{R}^n$ be a bounded domain and $u\in C(\overline{U })$ be a viscosity solution of
\begin{equation*}
\left\{\begin{aligned}
    \Delta_p u&=f&& \mbox{in } U \\
u&=g&& \mbox{on } \partial U,
\end{aligned}\right.
\end{equation*}
where $1<p<+\infty$. Suppose that $0\in \partial U$ and
\begin{equation*}
 \begin{aligned}
	f\in C(U )\cap L^{\infty}(U ), \quad  g\in C^{1,\alpha}(0), \quad
	\partial U\in C^{1,\alpha}(0)
  \end{aligned}
\end{equation*}
for some $\alpha\in (0,\bar{\alpha})$.

Then $u\in C^{1,\alpha}(0)$, i.e., there exists a linear polynomial $L$ such that
\begin{equation*}
    |u(x)-L(x)|\leq C |x|^{1+\alpha}\quad \text{for all }x\in U\cap B_1
\end{equation*}
and
\begin{equation*}
|DL|\leq  C,
\end{equation*}
where $C>0$ is a constant depending only on $n$, $p$, $\alpha$, $\|u\|_{L^{\infty}(U\cap B_1)}$, $\|f\|_{L^{\infty}(U\cap B_1)}$, $\|g\|_{C^{1,\alpha}(0)}$ and $\|\partial U\cap B_1\|_{C^{1,\alpha}(0)}$.
\end{corollary}

\begin{remark}\label{re10.5}
Extend $u$, $f$, $g$ to $\Omega\coloneqq U\times(-1,0]$ by the standard way (i.e., set $u(x,t)=u(x)$ etc.). Then $u$ is a viscosity solution of \eqref{e11.1} with $\gamma=p-2$. Note that
\begin{equation*}
g\in C^{1,\alpha}(0)\implies g\in C^{1,\alpha}_{\gamma}(0,0), \quad
\partial U\in C^{1,\alpha}(0)\implies \partial \Omega\in C^{1,\alpha}_{\gamma}(0,0)\quad \text{for all }\gamma\in (-1,+\infty).
\end{equation*}
Hence, by \Cref{th1.1}, we have $u\in C^{1,\alpha}(0,0)$ or $u\in C^{1,\alpha}_{\gamma}(0,0)$ with $\gamma=p-2$. Then by transferring to the domain $U$, we have $u\in C^{1,\alpha}(0)$.
\end{remark}

\begin{remark}\label{re10.4}
As far as we know, there is no boundary pointwise $C^{1,\alpha}$ regularity for the elliptic $p$-Laplace equation and hence \Cref{co1.0} is new. For the $p$-Laplace equations, the notion of viscosity solution is equivalent to the notion of weak solution (see \cite{MR2915869,MR1871417,MR3918385}), hence \Cref{co1.0} is also valid for weak solutions.
\end{remark}

\begin{remark}\label{re15.3}
As pointed out in \Cref{re10.1}, if $p$ is close to $2$, $\bar{\alpha}$ is close to $1$. Then we can obtain a higher boundary regularity.
\end{remark}

By combining the interior regularity with the boundary pointwise regularity, we have the following global $C^{1,\alpha}$ regularity.
\begin{theorem}\label{th1.2}
Let $u\in C(\overline{\Omega\cap Q_1})$ be a viscosity solution of
\begin{equation*}
\left\{\begin{aligned}
Pu&=f&& \mbox{in } \Omega\cap Q_1\\
u&=g&& \mbox{on } \partial_p \Omega\cap Q_1.
\end{aligned}\right.
\end{equation*}
Suppose that
\begin{equation*}
 \begin{aligned}
	f\in C(\Omega\cap Q_1)\cap L^{\infty}(\Omega\cap Q_1)
  \end{aligned}
\end{equation*}
and
\begin{equation*}
  \left\{
  \begin{aligned}
  g&\in C^{1,\alpha}_{\gamma}(\overline{\partial_p\Omega\cap Q_1}), \quad \partial_p\Omega\cap Q_1\in C^{1,\alpha}_{\gamma}\quad\mbox{if }\gamma> 0\\
  g&\in C^{1,\alpha}(\overline{\partial_p\Omega\cap Q_1}), \quad \partial_p\Omega\cap Q_1\in C^{1,\alpha}\quad
\mbox{if }\gamma\leq 0\\
  \end{aligned}
  \right.
\end{equation*}
for some $0<\alpha< \bar{\alpha}$.

Then
\begin{equation*}
  \left\{
  \begin{aligned}
u&\in C^{1,\alpha}(\overline{\Omega\cap Q_{1/2}})\quad\mbox{if }\gamma> 0\\
u&\in C_{\gamma}^{1,\alpha}(\overline{\Omega\cap Q_{1/2}})\quad\mbox{if }\gamma\leq 0\\
  \end{aligned}
  \right.
\end{equation*}
and
\begin{equation*}
  \left\{
  \begin{aligned}
  \|u\|_{C^{1,\alpha}(\overline{\Omega\cap Q_{1/2}})}&\leq C \quad\mbox{if }\gamma> 0\\
\|u\|_{C^{1,\alpha}_{\gamma}(\overline{\Omega\cap Q_{1/2}})}&\leq C \quad\mbox{if }\gamma\leq 0,
  \end{aligned}
  \right.
\end{equation*}
where $C>0$ is a constant depending only on $n$, $p$, $\gamma$, $\alpha$, $\|u\|_{L^{\infty}(\Omega\cap Q_1)}$, $\|f\|_{L^{\infty}(\Omega\cap Q_1)}$, $\|g\|_{C^{1,\alpha}_{\gamma}(\overline{\partial_p\Omega\cap Q_1})}$ (or $\|g\|_{C^{1,\alpha}(\overline{\partial_p\Omega\cap Q_1})}$) and $\|\partial_p\Omega\cap Q_1\|_{C^{1,\alpha}_{\gamma}}$ (or $\|\partial_p\Omega\cap Q_1\|_{C^{1,\alpha}}$).
\end{theorem}

As a direct corollary, we have
\begin{corollary}\label{co1.1}
Let $\Omega=U\times (-1,0]$ and $u\in C(\overline{\Omega})$ be a viscosity solution of
\begin{equation*}
\left\{\begin{aligned}
Pu&=f&& \mbox{in } \Omega\\
u&=g&& \mbox{on } \partial U \times (-1,0].
\end{aligned}\right.
\end{equation*}
Suppose that
\begin{equation*}
 \begin{aligned}
	f\in C(\Omega)\cap L^{\infty}(\Omega), \quad
	\partial U\in C^{1,\alpha}
  \end{aligned}
\end{equation*}
and
\begin{equation*}
  \left\{
  \begin{aligned}
g&\in C^{1,\alpha}_{\gamma}(\partial U\times (-1,0])\quad\mbox{if }\gamma> 0\\
g&\in C^{1,\alpha}(\partial U\times (-1,0])\quad\mbox{if }\gamma\leq 0\\
  \end{aligned}
  \right.
\end{equation*}
for some $0<\alpha< \bar{\alpha}$.

Then
\begin{equation*}
  \left\{
  \begin{aligned}
u&\in C^{1,\alpha}(\overline{U}\times [-1/2,0]) &&\mbox{if }\gamma> 0\\
u&\in C_{\gamma}^{1,\alpha}(\overline{U}\times [-1/2,0])&&\mbox{if }\gamma\leq 0\\
  \end{aligned}
  \right.
\end{equation*}
and
\begin{equation*}
  \left\{
  \begin{aligned}
\|u\|_{C^{1,\alpha}(\overline{U}\times [-1/2,0])} &\leq C &&\mbox{if }\gamma> 0\\
\|u\|_{C^{1,\alpha}_{\gamma}(\overline{U}\times [-1/2,0])}&\leq C &&\mbox{if }\gamma\leq 0,
  \end{aligned}
  \right.
\end{equation*}
where $C>0$ is a constant depending only on $n$, $p$, $\gamma$, $\alpha$, $\|u\|_{L^{\infty}(\Omega)}$, $\|f\|_{L^{\infty}(\Omega)}$, $\|g\|_{C^{1,\alpha}_{\gamma}(\partial U\times (-3/4,0])}$ (or $\|g\|_{C^{1,\alpha}(\partial U\times (-3/4,0])}$) and $\|\partial U\|_{C^{1,\alpha}}$.
\end{corollary}


Let us provide some remarks on the proof of our main theorems. First of all, since we are concerned with nonhomogeneous Dirichlet boundary data $g$ and we need to apply approximation, normalization and scaling techniques, we deal with operators $P^{\varepsilon}_{a, \nu}$ in a more general form than $P$; see \eqref{e2.7} for the precise definition. In particular, the term $|Du|^{\gamma}$, which represents the degeneracy or singularity in the original operator $P$, is transformed to $|\nu Du+a|^{\gamma}$ for some $a \in \mathbb{R}^n$ and $\nu \in [0,1]$.

Moreover, when we discuss the regularity for solutions of such generalized operators $P^{\varepsilon}_{a, \nu}$, the proofs and the associated scalings (for solutions and equations) belong to one of the two essential schemes: \emph{nondegenerate} or \emph{degenerate}. If the equation is nondegenerate (i.e., $|\nu Du+a| \approx 1$), then we use rather classical method for uniformly parabolic equations. Roughly speaking, this situation happens when either
\begin{enumerate}[label=(\roman*)]
    \item $|a| \gg \nu$, or
    \item $|a| \ll \nu$ and $|Du|$ is not small.
\end{enumerate}
In this nondegenerate scheme, we consider the classical scaling for uniformly parabolic equations given by (for some $r>0$)
\begin{equation*}
y=\frac{x}{r}, \quad s=\frac{t}{r^2} \quad\text{and}\quad v(y,s)=\frac{u(x,t)}{r^\kappa}
\end{equation*}
for some $\kappa\in \mathbb{R}$ together with the usual parabolic cylinders $Q_r$ and $Q_r^+$. We note that the small perturbation theorems play an important role in some stages of this nondegenerate scheme, and so we present appropriate versions of such theorems for both interior and boundary contexts in \Cref{S9}.

On the other hand, if the equation is degenerate (i.e., $|\nu Du+a| \ll 1$), then we prove the decay of oscillation of $|Du|$ by following the strategy of Imbert, Jin and Silvestre in \cite{MR3918407}. Roughly speaking, this situation happens when $|a| \ll \nu$ and $|Du|$ is small. In this degenerate scheme, we use the so-called two-parameter family of scaling (see \cite{MR783531,MR3918407,MR4676644}) given by (for some $r>0$ and $\rho>0$):
\begin{equation*}
y=\frac{x}{r}, \quad s=\frac{t}{\rho^{-\gamma} r^2} \quad\text{and}\quad v(y,s)=\frac{u(x,t)}{r\rho},
\end{equation*}
and the following special cylinders:
\begin{equation*}
Q^{\rho}_r\coloneqq B_r\times (-\rho^{-\gamma}r^2,0] \quad \text{and}\quad
Q^{\rho+}_r\coloneqq B_r^+\times (-\rho^{-\gamma}r^2,0].
\end{equation*}
The advantage of this \emph{intrinsic} scaling is that an equation like \eqref{e1.1} in $Q^{\rho}_r$ can be transformed to an equation in the same form in $Q_1$.  Moreover, we point out that the so-called \emph{Cutting Lemma}, which was used in the elliptic setting \cite{MR2995669} to remove degeneracy, is no longer applicable to parabolic equations due to the presence of time derivatives (i.e., $u_t$).

Let us finally illustrate sequential steps to arrive at the boundary pointwise $C^{1, \alpha}$ regularity for general boundary data $g$ on general boundary $\partial \Omega$ (\Cref{th1.1}).  We first obtain the boundary $C^{1,\alpha}$ regularity for the model problem (see \eqref{e2.1}), i.e., the problem with a flat boundary and homogeneous boundary data. Then we use the perturbation technique to derive the full regularity. To be more precise:
\begin{enumerate}[label=(\roman*)]
    \item (Interior $C^{0, 1}$ and $C^{1, \alpha}$ regularity; \Cref{S1}) We first utilize the Ishii--Lions method to establish the interior H\"older and then Lipschitz estimate in space for the general operator $P^{\varepsilon}_{a, \nu}$ (\Cref{lem-lipschitz-space}). Once we have the interior Lipschitz regularity, the interior $C^{1, \alpha}$ regularity (\Cref{th2.1}) follows from the existing regularity for $\nu=1$ (see \cite{MR4122677, MR3804259, MR4128334}) and the small perturbation theory provided in \Cref{L.K2-1}.

    \item (Boundary and global $C^{0, 1}$ regularity for the model problem; \Cref{S2}) We construct suitable barrier functions and apply the comparison principle to control the solution $u$ near the flat boundary. The construction of barriers depends on whether $\nu =1$ (\Cref{le2.1}) or $|a|$ is big (\Cref{le2.1-2}). Then by combining the interior regularity with the boundary regularity, we obtain the global $C^{0, 1}$ regularity (\Cref{le2.2}).

    \item (Boundary $C^{1, \alpha}$ regularity for the model problem when $|a| \gg \nu$; \Cref{S3}) Since this case can be regarded as the uniformly parabolic context, we develop the rather classical tools such as strong maximum principle, Harnack inequality and Hopf lemma. Then the boundary $C^{1, \alpha}$ regularity (\Cref{le2.3}) follows in a standard way. In fact, with the aid of the small perturbation regularity, we have higher regularity (\Cref{le3.4}) in this case.

    \item (Boundary estimates for the model problem when $|a| \ll \nu$; \Cref{S4}) In this case, there are two further essential situations depending on whether $|Du|$ is small or not. If $Du$ is close to the vector $e_n=(0, \cdots, 0, 1)$ in a measure sense, then $u$ must be close to the linear function $x_n$ in a smaller cylinder (\Cref{le3.6}). Then the boundary $C^{1,\alpha}$ regularity follows from the small perturbation theory again.

        On the other hand, if $|Du|$ is small in a measure sense, then this smallness information itself guarantees a decay of $|Du|$ in a smaller cylinder (\Cref{le3.3}). To be precise, we first combine the measure-type smallness assumption on $\partial_nu$ and the homogenous Dirichlet boundary condition to obtain the decay of $\partial_n u$ on the flat boundary (\Cref{le3.2}). This step is inspired by \cite[Lemma 1.2]{MR1207530}. Then by applying an argument similar to the proof of \cite[Lemma 4.1]{MR3918407} to obtain the decay of $|Du|$ in a smaller cylinder.

    \item (Boundary $C^{1, \alpha}$ regularity for the model problem; \Cref{S6}) Based on the consequences derived in \Cref{S3} and \Cref{S4}, we use an iteration presented in \cite[Theorem 4.8]{MR3918407} to obtain the boundary $C^{1,\alpha}$ regularity, regardless of the values of $|a|$ and $\nu$. We point out that the result of \Cref{S3} is necessary since we need to deal with an additional case coming from the presence of $a$, which does not happen for the interior regularity as in \cite{MR3918407}.

    \item (Boundary and global $C^{1, \alpha}$ regularity for general problems; \Cref{S7}) We use the perturbation argument (based on the compactness method) to prove the boundary pointwise $C^{1,\alpha}$ regularity for general problems (\Cref{th3.2}, \Cref{th3.3} and \Cref{th1.1}). We adopt an iteration formula inspired by \cite[Lemma 4.3]{MR4122677} and \cite[Corollary 3.3]{MR4128334}. If the equation is degenerate at the point we are concerned, the iteration can continue to infinity, which implies the boundary $C^{1,\alpha}$ regularity. Otherwise, by the small perturbation regularity, we have he boundary $C^{1,\alpha}$ regularity as well. Finally, by combining the interior $C^{1, \alpha}$ regularity with the boundary $C^{1, \alpha}$ regularity, we obtain the global regularity (\Cref{th1.2}).
\end{enumerate}

\begin{notation}\label{no1.1}
We summarize some basic notation as follows.
\begin{enumerate}
	\item Standard basis of $\mathbb{R}^n$: $\mathfrak{B} =\{e_i\}^{n}_{i=1}$, where   $e_i=(0,\cdots, 0,\underset{i^{\textnormal{th}}}{1},0, \cdots, 0)\in \mathbb{R}^n$.
	\item Points: $x'=(x_1, \cdots ,x_{n-1}) \in \mathbb{R}^{n-1}$, $x=(x',x_n) \in \mathbb{R}^{n}$ and $(x,t)=(x',x_n,t)\in \mathbb{R}^{n+1}$.
	\item Norms:
	The Euclidean norm is defined as $|x|=\left(\sum_{i=1}^{n} x_i^2\right)^{1/2}$ for $x\in \mathbb{R}^n$.
	The parabolic norm is defined as $|(x,t)|= (|x|^2+|t|)^{1/2}$ for $(x,t)\in \mathbb{R}^{n+1}$.
	\item $\mathbb{R}^n_+=\{x\in \mathbb{R}^n \mid x_n>0\}$ and $\mathbb{R}^{n+1}_+=\{(x,t)\in \mathbb{R}^{n+1} \mid x_n>0\}$.
	\item $B_r(x_0)=\{x\in \mathbb{R}^{n} : |x-x_0|<r\}$ and $B_r=B_r(0)$.
	\item $B_r^+(x_0)=B_r(x_0)\cap \mathbb{R}^n_+$ and $B_r^+=B^+_r(0)$.
	\item $Q_r(x_0,t_0)=B_r(x_0)\times (t_0-r^{2},t_0]$ and $Q_r=Q_r(0,0)$.
	\item $Q_r^+(x_0,t_0)=Q_r(x_0,t_0)\cap \mathbb{R}^{n+1}_+$ and $Q_r^+=Q^+_r(0,0)$.
    \item $Q^\rho_r=B_r\times (-\rho^{-\gamma}r^2,0]$, $Q^{\rho+}_r=B_r^+\times (-\rho^{-\gamma}r^2,0]$. Similarly, we can define $Q^\rho_r(x_0,t_0)$ etc.
	\item $T_r(x_0) =B_r(x_0) \cap \{x \in \mathbb{R}^{n} \mid x_n=0\}$ and $T_r=T_r(0)$.
	\item $S_r(x_0,t_0)\ =T_r(x_0)\times (t_0-r^{2},t_0]$ and $S_r=S_r(0,0)$.
    \item $\Omega_r=\Omega\cap Q_r$, $(\partial_p\Omega)_r=\partial_p\Omega\cap Q_r$, $\Omega_r^{\rho}=\Omega\cap Q^\rho_r$ and $(\partial_p\Omega)_r^{\rho}=\partial_p\Omega\cap Q_r^{\rho}$.
\item $\Omega^c$: the complement of $\Omega$; $\overline \Omega $: the closure of $\Omega$, where $\Omega\subset \mathbb{R}^{n+1}$.
    \item $\|M\|$: the spectral radius of an $n\times n$ symmetric matrix $M$.
    \item $u^+ \coloneqq \max\{u,0\}$, the positive part of $u$; $u^- \coloneqq \max\{-u,0\}$, the negative part of $u$ for a function $u: \Omega \to  \mathbb{R}$.
    \item $u_i=\partial u/\partial x_i$, $u_{ij}=\partial^2 u/\partial x_i \partial x_j$, $Du=(u_1,\cdots,u_n)$: the gradient of $u$.
\end{enumerate}
\end{notation}

%
%
\section{Interior \texorpdfstring{$C^{0, 1}$}{C0,1} and \texorpdfstring{$C^{1,\alpha}$}{C1,a} regularity}\label{S1}
Since we need to apply approximation, normalization, scaling etc. in later sections, it is useful to consider an operator $P^{\varepsilon}_{a,\nu}$ in the following more general form:
\begin{equation}\label{e2.7}
    P^{\varepsilon}_{a,\nu}u\coloneqq  u_t-(|\nu D u+a |^2+\varepsilon^2)^{\gamma/2} \left(\delta_{ij}+(p-2)\frac{(\nu u_i+a_i)(\nu u_j+a_j)}{|\nu D u+a|^2+\varepsilon^2}\right)u_{ij},
\end{equation}
where $\varepsilon$, $a$, $\nu$ are three parameters satisfying
\begin{equation*}
    0\leq \varepsilon\leq 1, \quad
    a\in \mathbb{R}^n \quad\text{and}\quad
    0\leq \nu\leq 1.
\end{equation*}
We also write for simplicity:
\begin{equation*}
    P^{\varepsilon}_{a}\coloneqq P^{\varepsilon}_{a,1}
    \quad\text{and}\quad
    P_a\coloneqq P^{0}_{a}.
\end{equation*}
In particular, the original operator in \eqref{e1.1} can be written as $P=P_0=P^{0}_0=P^0_{0,1}$.

The constant $\varepsilon$ is the approximation parameter and $\varepsilon=0$ is allowed in this section. In fact, in \Cref{S6}, we first prove uniform $C^{1,\alpha}$ estimates (independent of $\varepsilon$) for $P^{\varepsilon}_a$ (see \Cref{le5.1}) and then obtain the $C^{1,\alpha}$ regularity for $P_a$ by letting $\varepsilon\to 0$ (see \Cref{th3.1}). The vector $a$ and constant $\nu$ appear when we deal with scaling.

In this section, we prove the interior $C^{1,\alpha}$ regularity for a viscosity solution $u$ of
\begin{equation} \label{eq-approx}
    P^{\varepsilon}_{a,\nu}u=f\quad\mbox{in } Q_1
\end{equation}
 and always assume that
\begin{equation}\label{e2.17}
    |a| \leq 1, \quad
    \|u\|_{L^{\infty}(Q_1)}\leq 1 \quad\text{and}\quad
    \|f\|_{L^{\infty}(Q_1)}\leq 1.
\end{equation}
The assumptions on $a,u$ and $f$ are not restrictive, since we can use some normalization scheme to transfer a general case into \eqref{e2.17} (see \Cref{th2.2}).

We use the notion of viscosity solution introduced by Ohnuma and Sato \cite[Definition 2.4]{MR1443043} and Demengel \cite[Definition 1]{MR2804550}, which are equivalent (see \cite[Appendix]{MR2804550}). Recall that a function $\varphi $ touching $u$ from above at $(x,t)$ means: there exists an open neighborhood $Q$ of $(x,t)$ such that
\begin{equation*}
u \leq \varphi  \quad \mbox{in } Q \quad  \mbox{and} \quad u(x,t) = \varphi(x,t).
\end{equation*}
We define $\varphi$ touching $u$ from below similarly.
\begin{definition} [Viscosity solution]\label{de2.1}
Let $u$ and $f$ be continuous functions in $\Omega$. We say that $u$ is a \textit{viscosity subsolution} (resp. \textit{supersolution}) of
\begin{equation*}
    P^{\varepsilon}_{a,\nu}u=f\quad\mbox{in } \Omega
\end{equation*}
if for any $(x_0, t_0) \in \Omega$, both of the following two statements hold:
\begin{enumerate} [label=(\roman*)]
    \item For any smooth function $\varphi$ touching $u$ from above (resp. below) at $(x_0,t_0)$ with $\nu D\varphi(x_0,t_0)+a\neq 0$, we have
\begin{equation*}
    P^{\varepsilon}_{a,\nu}\varphi(x_0,t_0)\leq (\mathrm{resp. }\geq)\, f(x_0,t_0).
\end{equation*}
    \item For any $\varphi\in C^1(t_0-\delta, t_0+\delta)$ (for some $\delta> 0$) such that
\begin{equation*}
\nu u(x_0,t_0)+a\cdot x_0-\varphi(t_0)\geq (\mathrm{resp. }\leq)\, \nu u(x_0,t)+a\cdot x_0-\varphi(t)
\end{equation*}
for all $t\in (t_0-\delta, t_0+\delta)$ and
\begin{equation*}
\sup_{t\in (t_0-\delta, t_0+\delta)} (\mathrm{resp. }\inf)\,(\nu u(x,t)+a\cdot x-\varphi(t))
=\nu u(x_0,t_0)+a\cdot x_0-\varphi(t_0)
\end{equation*}
for all $x\in B_{\delta}(x_0)$, we have
\begin{equation*}
\varphi'(t_0)\leq (\mathrm{resp. }\geq)\, \nu f(x_0,t_0).
\end{equation*}
\end{enumerate}
\end{definition}
\begin{remark}\label{re1.1}
It follows from the definition that if $\nu\neq 0$, then $u$ is a viscosity solution of \eqref{eq-approx} if and only if $v(x,t):=\nu u(x,t)+a\cdot x$ is a viscosity solution of
\begin{equation*}
P^{\varepsilon}_{0,1}v=\nu f\quad\mbox{in}~~Q_1
\end{equation*}
in the sense of \cite[Definition 2.4]{MR1443043} and \cite[Definition 1]{MR2804550}. Hence, if $\nu\neq 0$, the \Cref{de2.1} is essentially the same as \cite[Definition 2.4]{MR1443043} and \cite[Definition 1]{MR2804550}.

If $\nu=0$ and $a\neq 0$, the case (ii) in \Cref{de2.1} can not happen for any $\varphi$ and the case (i) must happen. Then the \Cref{de2.1} is the same as the classical definition of viscosity solution for uniformly parabolic equation (see \cite[Definition 3.4]{MR1135923}).

Note that it cannot happen that $\nu=0$ and $|a|=0$ simultaneously throughout this paper.
\end{remark}

\begin{remark}\label{re1.3}
From the definition, we can prove directly the following comparison principle, which will be used to construct barriers (see \Cref{le2.1}, \Cref{le2.1-2} etc.). Let $u$ be a viscosity subsolution of \eqref{eq-approx} and $v$ be a smooth supersolution of \eqref{eq-approx}. If we assume that $v \geq u$ on $\partial_p\Omega$ and $\nu Dv+a\neq 0$ in $\Omega$, then
\begin{equation*}
v\geq u\quad\mbox{in}~~\Omega.
\end{equation*}
Indeed, if there exists $(x_0,t_0)\in \Omega$ such that $v(x_0,t_0)<u(x_0,t_0)$, then for some constant $c>0$, $v+c$ will touch $u$ from above at some $(x_1,t_1)\in \Omega$. By the definition of viscosity solution,
\begin{equation*}
P^{\varepsilon}_{a,\nu}v(x_1,t_1)\leq f(x_1,t_1),
\end{equation*}
which is a contradiction.
\end{remark}

The strategy is the following. We first use Ishii--Lions technique to obtain some uniform estimates with respect to $x$ (see \Cref{lem-lipschitz-space}) and $t$ (see \Cref{lem-holder-time}). These estimates provide the desired compactness and allow us to develop small perturbation regularity for \eqref{eq-approx} (we postpone this regularity to \Cref{S9}). Then by combining this small perturbation regularity with the existing interior $C^{1,\alpha}$ regularity when $\nu=1$, we can obtain the interior $C^{1,\alpha}$ regularity for \eqref{eq-approx} (see \Cref{th2.1}).

First, recall the following uniform Lipschitz regularity for the special case of $\nu=1$ and $a=0$, which was developed in \cite[Lemma 6.2]{MR4122677}.
\begin{lemma}[Uniform Lipschitz estimate when $\nu=1$ and $a=0$]\label{lem-IJS}
Let $u$ be a viscosity solution of \eqref{eq-approx} with $\nu=1$ and $a=0$. Then there exists a universal constant $C>0$ such that for every $(x, t), (y, t) \in Q_{1/2}$,
\begin{equation*}
|u(x, t)-u(y, t)| \leq C|x-y|.
\end{equation*}
\end{lemma}
\begin{corollary}[Uniform Lipschitz estimate: temporary version]\label{lem-temp}
Let $u$ be a viscosity solution of \eqref{eq-approx} with $\nu \in (0, 1]$. Then there exists a universal constant $C>0$ such that for every $(x, t), (y, t) \in Q_{1/2}$,
\begin{equation*}
|u(x, t)-u(y, t)| \leq C\left(\nu^{-\frac{\gamma}{1+\gamma}}+\nu^{-1}\big(|a|+|a|^{\frac{1}{1+\gamma}}\big)\right)|x-y|.
\end{equation*}
\end{corollary}
\begin{proof}
	
It immediately follows from \Cref{lem-IJS} combined with standard scaling and translation argument.
\end{proof}

We are going to improve this temporary version of Lipschitz estimate by developing uniform estimates with respect to $\nu$. For this purpose, we first prove the uniform H\"older estimates in a different setting when $\nu$ is relatively smaller than $|a|$. We use Ishii--Lions method (see \cite[Section VII]{MR1031377} and \cite[Theorem 5]{MR1341739}) as in \cite[Lemma 4]{MR2995669}, \cite[Section 2]{MR3918407}, \cite[Section 6]{MR4122677} etc.
\begin{lemma}[Uniform H\"older estimate in space]\label{lem-holder-space}
Suppose that $1/2\leq |a| \leq 1$. Then there exist universal constants $\kappa \in (0,1)$, $\nu_0 \in (0,1)$ and $C>0$ such that if $u$ is a viscosity solution of \eqref{eq-approx} with $0<\nu \leq \nu_0$, then
\begin{equation*}
|u(x, t)-u(y, t)| \leq C|x-y|^{\kappa} \quad \text{for any $(x, t), (y, t) \in Q_{1/2}$}.
\end{equation*}
\end{lemma}

\begin{proof}
Without loss of generality, we prove the H\"{o}lder continuity at $(0,0)$. We only need to prove
\begin{equation}\label{claim}
		M\coloneqq\max_{\substack{x, y \in \overline{B_{1/2}}  \\ t \in [-1/4, 0] }} \left[u(x, t)-u(y,t)-K_1\phi(|x-y|)-\frac{K_2}{2}|x|^2-\frac{K_2}{2}|y|^2-\frac{K_2}{2}t^2 \right] \leq 0,
\end{equation}
where
\begin{equation*}
		\phi(r)\coloneqq r^{\kappa} \quad \text{for $\kappa \in (0,1)$ to be determined soon.}
\end{equation*}
Indeed, if \eqref{claim} holds, by setting $(y,t)=(0,0)$ in \eqref{claim}, we have
\begin{equation*}
u(x,0)-u(0,0)\leq K_1|x|^{\kappa}+\frac{K_2}{2}|x|^2\leq C|x|^{\kappa} \quad \text{for all } x\in B_{1/2}.
\end{equation*}
Similarly,  by setting $(x,t)=(0,0)$ in \eqref{claim}, we have
\begin{equation*}
u(0,0)-u(y,0)\leq C|y|^{\kappa}\quad \text{for all } y\in B_{1/2}.
\end{equation*}
Hence, $u$ is $C^{\kappa}$ (with respect to $x$) at $(0,0)$.

We prove \eqref{claim} by contradiction: suppose that the positive maximum $M$ is attained at $t \in [-1/4, 0]$ and $x, y \in \overline{B_{1/2}}$. It immediately follows that $x \neq y$ and
\begin{equation*}
K_1\phi(|x-y|)+\frac{K_2}{2}|x|^2+\frac{K_2}{2}|y|^2+\frac{K_2}{2}t^2 \leq 2\|u\|_{L^{\infty}(Q_1)}\leq 2.
\end{equation*}
In particular,
\begin{equation*}
|x|^2+|y|^2+|t|^2 \leq \frac{4}{K_2}
\end{equation*}
and
\begin{equation}\label{eq-theta}
\phi(\theta)\leq \frac{|u(x, t)-u(y, t)|}{K_1} \leq \frac{2}{K_1}, \quad \text{where $\theta=|b|$ and $b=x-y$}.
\end{equation}
We fix $K_2>0$ large enough to ensure $t \in (-1/4, 0]$ and $x, y \in B_{1/2}$.
	
We now apply the parabolic version of Jensen--Ishii's lemma \cite[Theorem 8.3]{MR1118699}, to obtain that, for every $\delta>0$ sufficiently small, there exist $(\sigma_x, q_x, X) \in \overline{\mathcal{P}}^{2, +}u(x, t)$ and $(\sigma_y, q_y,  Y) \in \overline{\mathcal{P}}^{2, -}u(y, t)$ such that
\begin{itemize}
			\item
			$
			\begin{pmatrix}
				X & 0 \\
				0 & -Y
			\end{pmatrix}
			\leq K_1
			\begin{pmatrix}
				Z & -Z \\
				-Z & Z
			\end{pmatrix}
			+(2K_2+\delta)
			\begin{pmatrix}
				I & 0 \\
				0 & I
			\end{pmatrix}
			;$
			\item $\sigma_x-\sigma_y=K_2t$,
\end{itemize}
where
\begin{align*}
			q&\coloneqq K_1 \phi'(\theta)\hat{b}, \quad q_x \coloneqq q+K_2x, \quad q_y \coloneqq q-K_2y,\\
			Z&\coloneqq\phi''(\theta) \hat{b} \otimes \hat{b}+\frac{\phi'(\theta)}{\theta} (I-\hat{b}\otimes \hat{b}) \quad \text{and} \quad \hat{b}\coloneqq\frac{b}{|b|}=\frac{x-y}{|x-y|}.
\end{align*}
For the convenience of notation, we let
\begin{equation*}
			A[z]=I+(p-2)\frac{z_iz_j}{\varepsilon^2+|z|^2} \quad \text{and} \quad B[z]=(\varepsilon^2+|z|^2)^{\frac{\gamma}{2}}A[z] \quad \text{for $z \in \mathbb{R}^n$}.
\end{equation*}
Then it follows from the observations above and the definition of viscosity solution that
\begin{equation}\label{e2.15}
\begin{aligned}
    K_2t &\leq \mathrm{Tr}(B[\nu(q+K_2x)+a]X)-\mathrm{Tr}(B[\nu(q-K_2y)+a]Y)+f(x, t)-f(y, t)\\
    &\leq\underbrace{\mathrm{Tr} \left( (B[\nu(q+K_2x)+a]-B[\nu(q-K_2y)+a])  X \right)}_{\eqqcolon T_1}\\
    &\quad+\underbrace{\mathrm{Tr}(B[\nu(q-K_2y)+a](X-Y))}_{\eqqcolon T_2}+2 
\end{aligned}
\end{equation}

Before we estimate two terms $T_1$ and $T_2$, we first provide appropriate $L^{\infty}$-bounds for $q$, $q_x$, $q_y$, $X$ and $Y$. By choosing $K_1$ large enough, $\theta$ will be small, $|\phi'(\theta)|$ and $|q|$ will be large, and so we have
\begin{equation*}
\frac{|q|}{2} \leq |q+K_2x| \leq 2|q| \quad \text{and} \quad \frac{|q|}{2} \leq |q-K_2y| \leq 2|q|.
\end{equation*}
Moreover, by \eqref{eq-theta} and \Cref{lem-temp}, we have
\begin{equation}\label{eq-phi}
\nu|q|=\nu\kappa K_1 |x-y|^{\kappa-1} \leq \nu\kappa \frac{|u(x, t)-u(y, t)|}{|x-y|}
\leq \kappa C \left(\nu +\nu^{\frac{1}{1+\gamma}}+(|a|+|a|^{\frac{1}{1+\gamma}})\right)
\leq \frac{|a|}{4},
\end{equation}
if we choose $\kappa \in (0,1)$ sufficiently small, which is universal. Therefore, we conclude that
\begin{equation*}
\frac{|a|}{2} \leq |\nu(q+K_2x)+a| \leq 2|a| \quad \text{and} \quad \frac{|a|}{2} \leq |\nu(q-K_2y)+a| \leq 2|a|.
\end{equation*}
	
On the other hand, the matrix inequality for $X$ and $Y$ together with the fact that $\phi''(\theta)<0$ yields that
\begin{align*}
X, -Y \leq K_1 \frac{\phi'(\theta)}{\theta}(I-\hat{b}\otimes \hat{b})+K_2I.
\end{align*}
By combining these estimates and recalling that
\begin{equation*}
\min\{p-1, 1\} I \leq A[z] \leq 	\max\{p-1, 1\} I,
\end{equation*}
we have
\begin{equation*}
\begin{aligned}
    &\mathrm{Tr}(A[\nu(q+K_2x)+a]X) \\
    &\quad \geq (\varepsilon^2+|\nu (q+K_2x)+a|^2)^{-\frac{\gamma}{2}}K_2t\\
    &\qquad+ \left(\frac{\varepsilon^2+|\nu (q-K_2y)+a|^2}{\varepsilon^2+|\nu (q+K_2x)+a|^2} \right)^{\frac{\gamma}{2}}\mathrm{Tr}(A[\nu(q-K_2y)+a]Y)\\
				&\quad \geq -C\left(|a|^{-\gamma}+K_1 \frac{\phi'(\theta)}{\theta}+1\right) \geq -C\left(K_1 \frac{\phi'(\theta)}{\theta}+1\right).
\end{aligned}
\end{equation*}
Therefore, we conclude that
\begin{align*}
\|X\|, \|Y\| \leq C\left(K_1 \frac{\phi'(\theta)}{\theta}+1 \right).
\end{align*}

We are now ready to estimate $T_1$ and $T_2$. For $T_1$, an application of mean value theorem gives
\begin{equation}\label{eq-T_1}
T_1 \leq C\nu|a|^{\gamma-1}|x+y| \|X\| \leq C\nu\left(K_1 \frac{\phi'(\theta)}{\theta}+1 \right).
\end{equation}
For $T_2$, we again utilize the previous matrix inequality. First, by evaluating a vector of the form $(\xi, \xi)$ for any $\xi \in \mathbb{R}^n$, we have
\begin{align*}
(X-Y)\xi \cdot \xi \leq 3K_2|\xi|^2,
\end{align*}
which implies that any eigenvalues of $X-Y$ are less than $3K_2$. Next, by considering a special vector $(\hat{b}, -\hat{b})$, we arrive at
\begin{equation}\label{e2.16}
(X-Y)\hat{b} \cdot \hat{b} \leq 4K_1 \phi''(\theta)+3K_2.
\end{equation}
In other words, at least one eigenvalue of $X-Y$ is less than $4K_1 \phi''(\theta)+3K_2$. Therefore, due to the uniform ellipticity of $A$, we have
\begin{equation} \label{est:T2}
T_2 \leq C\left(K_1\phi''(\theta)+1 \right)|a|^{\gamma}\leq  C\left(K_1\phi''(\theta)+1 \right).
\end{equation}
By combining the two estimates for $T_1$ and $T_2$ with \eqref{e2.15}, it holds that
\begin{equation*}
-K_1\phi''(\theta) \leq C \left( \nu K_1\frac{\phi'(\theta)}{\theta}+1\right).
\end{equation*}
By recalling that
\begin{equation*}
\phi''(\theta)=-\kappa(1-\kappa)\theta^{\kappa-2}=-(1-\kappa)\frac{\phi'(\theta)}{\theta},
\end{equation*}
we can choose $\nu_0$ small enough such that
\begin{equation*}
-K_1 \phi''(\theta) \leq C.
\end{equation*}
It leads to the contradiction if we choose $K_1$ large enough.
\end{proof}
\begin{remark}[The choice of $\phi$]\label{re2.3}
Since $|a|$ is relatively larger than $\nu$, one can expect that the degeneracy or singularity term $|a+\nu q_x|$ is comparable to $|a|$. In order to justify this fact, we controlled the term $\nu|q|$ as in \eqref{eq-phi}, where the particular choice of H\"older function $\phi(r)=r^{\kappa}$ played a crucial role. In fact, this is the essential reason why we first develop the H\"older estimate of $u$, instead of the Lipschitz estimate of $u$ directly.

Furthermore, by Jensen--Ishii’s lemma and the definition of viscosity solution, we have the inequality \eqref{e2.15}. The key for obtaining a contradiction at the final stage is that $X-Y$ becomes very negative along the direction $x-y$ (see \eqref{e2.16}). This property leads to that $T_2\ll -1$. Then based on \eqref{e2.15}, we obtain a contradiction since $K_2t$ and $T_1$ can be controlled by $T_2$. At this step, the choice of $\phi$ is again exploited to guarantee that $\phi'' \ll -1$.
\end{remark}

We next improve the $C^{\kappa}$ regularity to $C^{0,1}$ regularity.
\begin{lemma}[Uniform Lipschitz estimate in space]\label{lem-lipschitz-space}
Let $\nu_0$ be chosen in \Cref{lem-holder-space} and suppose that $1/2  \leq |a| \leq 1$. Then there exist universal constants $\nu_1 \in (0,\nu_0)$ and $C>0$ such that if $u$ is a viscosity solution of \eqref{eq-approx} with $0 < \nu \leq \nu_1$, then
\begin{equation*}
|u(x, t)-u(y, t)| \leq C|x-y| \quad \text{for all $(x, t), (y, t) \in Q_{1/2}$}.
\end{equation*}
\end{lemma}

\begin{proof}
The proof is similar to the one in \Cref{lem-holder-space}, but the estimate is improved by exploiting the H\"older regularity of a viscosity solution. Indeed, we again consider the regularity at $(0,0)$ and use the same notation as in the proof of \Cref{lem-holder-space} except for $\phi$.

It is enough to show \eqref{claim}, where
\begin{equation*}
\phi(r)\coloneqq
\begin{dcases}
r-\frac{1}{2-\gamma_0}r^{2-\gamma_0} \quad &\text{for $r \in [0, 1]$}\\
1-\frac{1}{2-\gamma_0} \quad &\text{for $r > 1$},
\end{dcases}
\end{equation*}
where the constant $\gamma_0 \in (1/2, 1)$ will be determined later. As before, we suppose that the positive maximum $M$ is attained at $t \in [-1/4, 0]$ and $x, y \in \overline{B_{1/2}}$ by choosing $K_2$ large enough. Let $\kappa$ be chosen in \Cref{lem-holder-space}. Note that \Cref{lem-holder-space} provides
\begin{equation*}
K_1 \phi(|x-y|)+\frac{K_2}{2}\left(|x|^2+|y|^2+|t|^2\right) \leq u(x, t)-u(y,t) \leq C|x-y|^{\kappa},
\end{equation*}
which implies that
\begin{equation*}
|t|+|x|+|y| \leq C \theta^{\kappa/2}.
\end{equation*}
We next estimate $|\nu(q+K_2x)+a|$ and $|\nu(q-K_2y)+a|$. By choosing $K_1$ large enough, $\theta$ will be small, $|q|$ will be large, and so we have
\begin{equation*}
\frac{|q|}{2} \leq |q+K_2x| \leq 2|q| \quad \text{and} \quad \frac{|q|}{2} \leq |q-K_2y| \leq 2|q|.
\end{equation*}
Moreover, if we choose $\nu_1=\nu_1(n, p, \gamma, K_1) \in (0,1)$ sufficiently small, then we have
\begin{equation*}
\begin{aligned}
\nu |q| \leq \nu K_1 \leq \frac{|a|}{4},
\end{aligned}
\end{equation*}
which implies 	
\begin{equation*}
\frac{|a|}{2} \leq |\nu(q+K_2x)+a| \leq 2|a| \quad \text{and} \quad \frac{|a|}{2} \leq |\nu(q-K_2y)+a| \leq 2|a|.
\end{equation*}
The constant $K_1$ will be chosen later and its value does not depend on this choice of $\nu_1$.
	
Let us now suggest a better estimate for $|x+y|$ and so for $T_1$ as follows (compare it to the previous estimate \eqref{eq-T_1}):
\begin{align*}
T_1 \leq C\nu|a|^{\gamma-1}|x+y| \|X\| \leq  C\left(K_1 \frac{\phi'(\theta)}{\theta}+1 \right) \theta^{\kappa/2}.
\end{align*}
Therefore, by combining this with the estimate \eqref{est:T2} for $T_2$ and taking $\gamma_0 = 1-\kappa/4$, we have
\begin{equation*}
K_2 t \leq C(1+K_1 \theta^{\kappa/2-1}) + C\left(1-K_1 \theta^{\kappa/4-1} \right),
\end{equation*}
which implies that
\begin{equation*}
K_1 \theta ^{\kappa/4-1} \leq C (1+K_1 \theta^{\kappa/2-1}).
\end{equation*}

By choosing $K_1$ large enough, we again notice that $\theta$ becomes small enough, and hence we arrive at the contradiction when we choose $K_1$ further large enough. Note that this choice of $K_1$ does not depend on the previous choice of $\nu_1$ and hence we should first choose $K_1$ large and then $\nu_1$ small.
\end{proof}

We then make use of the uniform Lipschitz estimate in the space variable to prove the uniform H\"older estimate in the time variable. We point out that the proof strongly depends on the assumption that $\nu$ is relatively smaller than $|a|$; see \cite[Lemma 3.1]{MR3918407} for a similar result in the case of $\nu=1$ and $a=0$.
\begin{lemma}[Uniform H\"older estimate in time]\label{lem-holder-time}
Let $\nu_1$ be chosen in \Cref{lem-lipschitz-space} and suppose that $1/2 \leq |a| \leq 1$. Then there exist universal constants $\kappa' \in (0, 1)$, $\nu_2 \in (0,\nu_1)$ and $C>0$ such that if $u$ is a viscosity solution of \eqref{eq-approx} with $0<\nu \leq \nu_2$, then
	\begin{equation*}
		|u(x, t)-u(x, s)| \leq C|t-s|^{\kappa'/2} \quad \text{for all $(x, t), (x, s) \in Q_{1/2}$}.
	\end{equation*}
\end{lemma}
\begin{proof}
	Since $\nu \leq \nu_1$, we apply \Cref{lem-lipschitz-space} to have
	\begin{equation*}
		\|Du\|_{L^{\infty}(Q_{3/4})} \leq C.
	\end{equation*}
	In particular, if we choose $\nu_2$ small enough, then we can guarantee that
	\begin{equation*}
		\frac{|a|}{2} \leq |\nu Du+a| \leq 2|a|.
	\end{equation*}
	Therefore, we can understand $u$ as a viscosity solution of
	\begin{equation*}
		u_t \geq \mathcal{M}_{\lambda, \Lambda} ^-(D^2u)+f \quad \text{and} \quad u_t \leq \mathcal{M}_{\lambda, \Lambda} ^+(D^2u)+f,
	\end{equation*}
	where
	\begin{equation*}
		\lambda=
        \begin{cases}
             \min\{1, p-1\} \cdot(|a|/2)^{\gamma} & \text{if $\gamma \geq 0$} \\
             \min\{1, p-1\} \cdot(1+4|a|^2)^{\gamma/2} & \text{if $\gamma < 0$}
        \end{cases}
	\end{equation*}
    and
    \begin{equation*}
        \Lambda=
        \begin{cases}
             \max\{1, p-1\} \cdot(1+4|a|^2)^{\gamma/2} & \text{if $\gamma \geq 0$} \\
             \max\{1, p-1\} \cdot(|a|/2)^{\gamma} & \text{if $\gamma < 0$}.
        \end{cases}
    \end{equation*}

    Finally, the desired H\"older regularity of $u$ in time follows from the Krylov--Safonov theory (see \cite[Theorem 4.19]{MR1135923}).
\end{proof}

The estimates \Cref{lem-lipschitz-space} and \Cref{lem-holder-time} provide necessary compactness to prove the small perturbation regularity (see \Cref{S9}). In fact, with aid of the small perturbation regularity, we can obtain the following interior $C^{1,\alpha}$ regularity.
\begin{theorem}\label{th2.1}
Let $u$ be a viscosity solution of \eqref{eq-approx}. Suppose that $\nu=1$ or $1/2\leq |a|\leq 1$. Then there exist universal constants $C>0$ and $0<\bar{\alpha} <\min\{1/2,1/2(1+\gamma)\}$ such that $u\in C^{1,\bar\alpha}_{\gamma}(\overline{Q_{1/2}})$ and
\begin{equation}\label{e2.11}
   \|u\|_{C^{1,\bar\alpha}_{\gamma}(\overline{Q_{1/2}})}\leq C.
\end{equation}
\end{theorem}
\begin{proof}
For $\nu=1$, by the transformation $v=u+a\cdot x$, \eqref{e2.11} can be derived immediately by the interior $C^{1,\bar\alpha}$ regularity for the operator $P^{\varepsilon}_{0,1}$; see \cite[Theorem 1.1]{MR4122677} for $0<\gamma<+\infty$, \cite[Theorem 1.1]{MR3804259} for $\gamma=0$ and \cite[Theorem 1.1]{MR4128334} for $-1<\gamma<0$ (see also \cite[Theorem 1.1]{MR3918407} for the case $f\equiv 0$). For $1/2\leq |a|\leq 1$, the proof is divided into two cases: $\nu$ is small or big.
\begin{enumerate} [label=(\roman*)]
    \item If $\nu$ is small, then the interior $C^{1,\bar\alpha}_{\gamma}$ regularity follows from the perturbation theory. More precisely, if $\nu\leq \eta$, where $0<\eta<1$ (universal) is from \Cref{L.K2-1} with $\beta=1/2$ there, then we have the interior $C^{1,1/2}$ regularity by \Cref{L.K2-1}. Since $\bar{\alpha}<\min\{1/2,1/2(1+\gamma)\}$,
we have the interior $ C^{1,\bar\alpha}_{\gamma}$ regularity.
    \item If $\nu$ is big, the interior $ C^{1,\bar\alpha}_{\gamma}$ regularity follows from \cite{MR4122677, MR3804259, MR4128334} again. Indeed, if $\nu\geq \eta$, let
\begin{equation*}
v(x,t)=\nu u(x,t)+a\cdot x.
\end{equation*}
Then $v$ is a solution of
\begin{equation*}
    P^{\varepsilon}_{0,1}v=f\quad \mbox{in } Q_1.
\end{equation*}
By the interior $C^{1,\bar\alpha}_{\gamma}$ regularity again (see \cite{MR4122677, MR3804259, MR4128334}), $v\in C^{1,{\bar\alpha}}_{\gamma}(\overline{Q_{1/2}})$ and so $u\in C^{1, {\bar\alpha}}_{\gamma}(\overline{Q_{1/2}})$ with the uniform estimate \eqref{e2.11}.
\end{enumerate}
\end{proof}

\begin{remark}
The proof above shows that for the interior regularity, considering an equation with the term $|Du+a|^{\gamma}$ (instead of $|Du|^{\gamma}$) is not more difficult since we can use the transformation $v=u+a\cdot x$ to convert the former to the later. However, for the boundary regularity, the former is indeed more difficult. For example, consider the following problem:
\begin{equation*}
\left\{\begin{aligned}
    P^{\varepsilon}_{a}u&=0&& \mbox{in } Q_1^+ \\
	u&=0&& \mbox{on } S_1.
\end{aligned}\right.
\end{equation*}
If we use the previous transformation, then
\begin{equation*}
v=\sum_{i=1}^{n-1}a_ix_i\quad\mbox{on}~~S_1,
\end{equation*}
i.e., the original homogenous boundary condition becomes a nonhomogenous one, which is more complicated.
\end{remark}

By using a normalization technique, we have the following interior $C^{1,\alpha}$ regularity without the assumption \eqref{e2.17}.
\begin{theorem}\label{th2.2}
Let $u$ be a viscosity solution of \eqref{eq-approx} (without the assumption \eqref{e2.17}). Suppose that
\begin{equation*}
    \nu=1 \quad \text{or} \quad |a|\geq 1/2.
\end{equation*}
Then $u\in C^{1,\bar\alpha}_{\gamma}(\overline{Q_{1/2}})$ and
\begin{equation*}
   \|u\|_{C^{1,\bar\alpha}_{\gamma}(\overline{Q_{1/2}})}\leq C,
\end{equation*}
where $\bar{\alpha}$ is as in \Cref{th2.1} and $C>0$ depends only on $n$, $p$, $\gamma$, $|a|$, $\|u\|_{L^{\infty}(Q_1)}$ and $\|f\|_{L^{\infty}(Q_1)}$.
\end{theorem}
\begin{proof}
We just need to make some normalization such that \eqref{e2.17} holds and the assumptions of \Cref{th2.1} are satisfied. Consider the following transformation ($\rho_1 \in [1, \rho_2]$ to be specified later):
\begin{equation}\label{e2.20}
\begin{aligned}
    \rho_2&=\| u\|_{L^{\infty}(Q_1)}^2+\|f\|_{L^{\infty}(Q_1)}^2+|a|+1, \quad r=\rho_2^{-1/2}, \\
    \tilde{x}&=\frac{x}{r}, \quad \tilde{t}=\frac{t}{\rho_1^{-\gamma}r^2} \quad\text{and} \quad
    \tilde{u}(\tilde x,\tilde t)=\frac{u(x,t)}{\rho_2 r}.
\end{aligned}
\end{equation}
Then $\tilde{u}$ is a viscosity solution of
\begin{equation*}
    P^{\tilde \varepsilon}_{\tilde a,\tilde\nu}\tilde u=\tilde f\quad\mbox{in}~~Q_1,
\end{equation*}
where
\begin{equation*}
    \tilde{\varepsilon}= \rho_1^{-1} \varepsilon, \quad \tilde{a}=\rho_1^{-1}a, \quad \tilde{\nu}=\rho_1^{-1}\rho_2\nu \quad \text{and} \quad  \tilde{f}(\tilde{x},\tilde{t})=r\rho_1^{-\gamma}\rho_2^{-1}f(x,t).
\end{equation*}

We first observe that
\begin{equation*}
    \|\tilde u\|_{L^{\infty}(Q_1)}\leq 1 \quad \text{and} \quad \|\tilde f\|_{L^{\infty}(Q_1)}\leq 1.
\end{equation*}
If $\nu=1$, then we choose $\rho_1=\rho_2$ so that
\begin{equation*}
\tilde{\nu}=1 \quad \text{and} \quad |\tilde{a}|\leq 1.
\end{equation*}
Next, we consider the case $|a|\geq 1/2$.
\begin{enumerate} [label=(\roman*)]
    \item If $\rho_2\nu\geq 2|a|$, then we set $\rho_1=\rho_2 \nu$ ($\leq \rho_2$) so that
\begin{equation*}
    \tilde{\nu}=1 \quad \text{and} \quad |\tilde{a}|\leq 1/2.
\end{equation*}
    \item If $\rho_2\nu< 2|a|$, then we choose $\rho_1=2|a|$ so that
\begin{equation*}
    \tilde{\nu} < 1 \quad \text{and} \quad  |\tilde{a}|=1/2.
\end{equation*}
\end{enumerate}

In any cases, \Cref{th2.1} yields that $\tilde{u}\in C^{1,\bar\alpha}_{\gamma}(\overline{Q_{1/2}})$. By transforming back to $u$ and using standard covering arguments, we arrive at the desired estimate.
\end{proof}

We also have the following version of $C^{1,\alpha}$ regularity in a scaled cylinder $Q_r$, which will be used in next section.
\begin{corollary}\label{co2.1}
Let $u$ be a viscosity solution of
\begin{equation*}
    P_{a, \nu}^{\varepsilon}u=f\quad\mbox{in } Q_{r}
\end{equation*}
for some $r>0$. Suppose that
\begin{equation*}
    \nu=1 \quad \text{or} \quad |a|\geq 1/2.
\end{equation*}
Then $u\in C^{1,\bar\alpha}_{\gamma}(\overline{Q_{r/2}})$ and for any $(x,t),(x,s),(y,s) \in Q_{r/2}$, the following estimates hold:
\begin{equation*}
\begin{aligned}
    |u(x,t)-u(y,s)|&\leq C(|x-y|+|t-s|^{1/2}),\\
    |Du(x,t)-Du(y,s)|&\leq C\left(r^{-\bar\alpha}|x-y|^{\bar{\alpha}}
+r^{-\frac{2\bar{\alpha}}{2-\bar{\alpha}\gamma}}|t-s|^{\frac{\bar{\alpha}}{2-\bar{\alpha}\gamma}}\right),\\
    |u(x,t)-u(x,s)|&\leq Cr^{-\frac{\bar{\alpha}(2+\gamma)}{2-\bar{\alpha}\gamma}}|t-s|^{\frac{1+\bar{\alpha}}{2-\bar{\alpha}\gamma}},
\end{aligned}
\end{equation*}
where $\bar{\alpha}$ is as in \Cref{th2.1} and $C>0$ is a constant depending only on $n$, $p$, $\gamma$, $|a|$, $r^{-1}\|u\|_{L^{\infty}(Q_r)}$ and $r\|f\|_{L^{\infty}(Q_r)}$.
\end{corollary}
\begin{proof}
Under the transformation
\begin{equation*}
\begin{aligned}
    \tilde{x}=\frac{x}{r}, \quad \tilde{t}=\frac{t}{r^2} \quad \text{and} \quad
    \tilde{u}(\tilde x,\tilde t)=\frac{u(x,t)}{r},
  \end{aligned}
\end{equation*}
$\tilde{u}$ is a viscosity solution of
\begin{equation*}
    P^{\varepsilon}_{a,\nu}\tilde u=\tilde f\quad\mbox{in } Q_1,
\end{equation*}
where $\tilde{f}(\tilde{x},\tilde{t})=rf(x,t)$. Then \Cref{th2.2} shows that $\tilde{u}\in C^{1,\bar\alpha}_{\gamma}(\overline{Q_{1/2}})$ and
\begin{equation*}
   \|\tilde u\|_{C^{1,\bar\alpha}_{\gamma}(\overline{Q_{1/2}})}\leq C,
\end{equation*}
where $C>0$ is a constant depending only on $n$, $p$, $\gamma$, $|a|$, $\|\tilde u\|_{L^{\infty}(Q_1)}$ and $\|\tilde f\|_{L^{\infty}(Q_1)}$. By transforming back to $u$, we have $u\in C^{1,\bar\alpha}_{\gamma}(\overline{Q_{r/2}})$ with the desired estimate.
\end{proof}

%
%
\section{Boundary and global \texorpdfstring{$C^{0,1}$}{C0,1} regularity for the model problem}\label{S2}
As explained in the introduction, our strategy is to prove the boundary $C^{1,\alpha}$ regularity for the following model problem first:
\begin{equation}\label{e2.1}
\left\{\begin{aligned}
    P^{\varepsilon}_{a,\nu}u&=0&& \mbox{in } Q_1^+ \\
	u&=0&& \mbox{on } S_1
\end{aligned}\right.
\end{equation}
and then use the perturbation technique to derive the full regularity \Cref{th1.1}.

In this section, we prove the boundary and global $C^{0,1}$ regularity for \eqref{e2.1} and always assume that
\begin{equation}\label{e12.1}
|a| \leq 1 \quad \text{and} \quad \|u\|_{L^{\infty}(Q_1^+)}\leq 1.
\end{equation}
As in the previous section, we allow $\varepsilon=0$ in this section.

We prove the boundary $C^{0,1}$ regularity for two cases $\nu=1$ and $|a|\geq 1/2$, respectively.
\begin{lemma}\label{le2.1}
Let $u$ be a viscosity solution of \eqref{e2.1} with $\nu=1$. Then
\begin{equation*}
    |u(x,t)|\leq Cx_n \quad \text{for all } (x,t) \in Q_{1/2}^+,
\end{equation*}
where $C>0$ is universal.
\end{lemma}

\begin{proof}
Let
\begin{equation}\label{e2.9}
 	v(x,t)=C\left(1-|x+e_n|^{-\beta}\right)-t.
\end{equation}
By taking $\beta$ large enough first and then $C$ large enough (here we use $\gamma>-1$), $v$ satisfies
\begin{equation*}
\left\{\begin{aligned}
P^{\varepsilon}_{a}v &> 0&& \mbox{in } Q^+_{1}\\
    v&\geq 0&& \mbox{on } S_{1}\\
    v&\geq 1&& \mbox{on } \partial_p Q^+_{1}\setminus S_{1}.
\end{aligned}\right.
\end{equation*}
Note that $|Dv+a|\neq 0$ in $Q_1^+$. Then by the definition of viscosity solution (see also \Cref{re1.3}), we have
\begin{equation*}
    u\leq v\quad\mbox{in } Q_1^+.
\end{equation*}
By a direct calculation, we have (noting $v(0,0)=0$)
\begin{equation*}
    -Cx_n\leq -v\leq u\leq v\leq Cx_n \quad \mbox{on }
    \{(0',x_n,0) \mid 0\leq x_n \leq 1/2 \}.
\end{equation*}
By considering $v(x'-x'_0,x_n,t-t_0)$  for $(x'_0,t_0)\in S_{1/2}$ and similar arguments, we finish the proof.
\end{proof}
\begin{lemma}\label{le2.1-2}
Let $u$ be a viscosity solution of \eqref{e2.1} with $1/2\leq |a|\leq 1$
Then
\begin{equation*}
    |u(x,t)|\leq Cx_n \quad \text{for all } (x,t) \in Q_{1/2}^+,
\end{equation*}
where $C>0$ is universal.
\end{lemma}

\begin{proof}
Let
\begin{equation}\label{e2.8}
    \psi(x,t)=e^{-|x|^2/t}
\end{equation}
and
\begin{equation}\label{e10.7}
 	v(x,t)=2e^{\beta}\left(\psi^{\beta}(e_n,1)-\psi^{\beta}(x+e_n,t+1)\right)
 =2e^{\beta}\left(e^{-\beta}-\psi^{\beta}(x+e_n,t+1)\right).
\end{equation}
Choose $\beta$ large enough (universal) such that if
\begin{equation*}
1/4\leq |\tilde{a}|\leq 5/4 \quad \text{and}\quad    0\leq \varepsilon\leq 1
\end{equation*}
then
\begin{equation*}
 P^{\varepsilon}_{\tilde{a},0}v > 0 \quad \mbox{in } Q^+_{1}.
\end{equation*}

If $\nu \leq 1/(16\beta)$, we have $\nu |Dv|\leq 1/4$. Then $1/4 \leq |\nu Dv+a|\leq 5/4$ and hence
\begin{equation*}
\left\{\begin{aligned}
   P^{\varepsilon}_{a,\nu}v &> 0&& \mbox{in } Q^+_{1}\\
    v&\geq 0&& \mbox{on } S_{1}\\
    v&\geq 1&& \mbox{on } \partial_p Q^+_{1}\setminus S_{1}.
\end{aligned}\right.
\end{equation*}
As in the proof of \Cref{le2.1}, since $|\nu Dv+a|\neq 0$ in $Q_1^+$, by the definition of viscosity solution, we have
\begin{equation*}
u\leq v\quad\mbox{in}~~Q_1^+.
\end{equation*}
As in the last lemma, we have (noting $v(0,0)=0$)
\begin{equation}\label{e2.10}
    -Cx_n\leq -v\leq u\leq v\leq Cx_n \quad \mbox{on }
    \{(0',x_n,0) \mid 0\leq x_n \leq 1/2 \}.
\end{equation}

If $\nu\geq 1/(16\beta)$, we can take the same barrier as in \Cref{le2.1} (see \eqref{e2.9}), in which the constant $C$ depends also on $\beta$ now. Then we obtain \eqref{e2.10} again.

By considering $v(x'-x'_0,x_n,t-t_0)$  for $(x'_0,t_0)\in S_{1/2}$ and similar arguments, we finish the proof.
\end{proof}

\begin{remark}\label{re2.2}
The auxiliary function $\psi$ is a simplified version of the fundamental solution of the heat equation. For uniformly parabolic equations, constructing barriers based on the modification of the fundamental solution is a common method (e.g., \cite[P. 154]{MR1139064}). The basic properties of $\psi$ are the following.
\begin{itemize}
  \item $\psi$ is increasing in $t$ and decreasing in $|x|$;
  \item $\psi(\cdot,0)\equiv 0$ except at the origin $(0,0)$;
  \item For any uniformly parabolic equation, $\psi^{\beta}$ is a subsolution if we take $\beta$ large enough.
\end{itemize}
Based on above properties, we can construct the desired barrier easily.

If $|a|$ is big (e.g., $|a|\geq 1/2$), we treat the equation in the sprit of uniformly parabolic equations in this paper. Hence, in this case, we always construct barriers based on $\psi$.
\end{remark}

By combining the interior $C^{0,1}$ regularity with the boundary $C^{0,1}$ regularity, we can obtain the global $C^{0,1}$ regularity.
\begin{lemma}\label{le2.2}
Let $u$ be a viscosity solution of \eqref{e2.1}. Suppose that $\nu=1$ or $1/2\leq |a|\leq 1$. Then $u\in C^{0,1}(\overline{Q^+_{1/2}})$ and
\begin{equation}\label{e.Lipschitz-2}
	\|u\|_{C^{0,1}(\overline{Q^+_{1/2}})}\leq C,
\end{equation}
where $C>0$ is universal.
\end{lemma}

\begin{proof}
 Given $(\tilde{x},\tilde{t})\in Q^+_{1/4}$, denote $r=\tilde{x}_n$. Then $Q_r(\tilde{x},\tilde{t})\subset Q_1^+$. For any $(x,t)\in Q^+_{1/4}$ with $t\leq \tilde{t}$, if $(x,t)\in Q_{r/2}(\tilde{x},\tilde{t})$, then the interior Lipschitz estimate in $Q_r(\tilde{x},\tilde{t})$ (see \Cref{co2.1}) yields that
\begin{equation} \label{int_lip}
	|u(x,t)-u(\tilde{x},\tilde{t})|\leq C(|x-\tilde{x}|+|t-\tilde{t}|^{1/2}),
\end{equation}
where $C>0$ depends only on $n$, $p$, $\gamma$ and $r^{-1}\|u\|_{L^{\infty}(Q_r(\tilde{x},\tilde{t}))}$. By \Cref{le2.1} and \Cref{le2.1-2}, we have
\begin{equation*}
	\|u\|_{L^{\infty}(Q_r(\tilde{x},\tilde{t}))}\leq Cr,
\end{equation*}
where $C>0$ is universal. Hence, the constant $C$ in \eqref{int_lip} is indeed universal.

On the other hand, if $(x,t)\notin Q_{r/2}(\tilde{x},\tilde{t})$, by \Cref{le2.1} and \Cref{le2.1-2}, we have
\begin{equation}\label{e2.4}
\begin{aligned}
	|u(x,t)-u(\tilde{x},\tilde{t})|\leq& |u(x,t)|+|u(\tilde{x},\tilde{t})|\leq C(x_n+\tilde{x}_n)\\
\leq& C\left(|x_n-\tilde{x}_n|+\tilde{x}_n\right)\leq C(|x-\tilde{x}|+|t-\tilde{t}|^{1/2}),
\end{aligned}
\end{equation}
where $C>0$ is universal.

By combining \eqref{int_lip} and \eqref{e2.4}, we have $u\in C^{0,1}(\overline{Q^+_{1/4}})$. Then it is standard that $u\in C^{0,1}(\overline{Q^+_{1/2}})$ with the uniform estimate \eqref{e.Lipschitz-2}.
\end{proof}

With the aid of the two-parameter scaling, we have the following corollary.
\begin{corollary}\label{co3.1}
Let $u$ be a viscosity solution of \eqref{e2.1} (without the assumption \eqref{e12.1}). Suppose that $\nu=1$ or $|a|\geq 1/2$. Then $u\in C^{0,1}(\overline{Q^+_{1/2}})$ and
\begin{equation*}
	\|u\|_{C^{0,1}(\overline{Q^+_{1/2}})}\leq C,
\end{equation*}
where $C>0$ depends only on $n$, $p$, $\gamma$, $|a|$ and $\|u\|_{L^{\infty}(Q_1)}$.
\end{corollary}

\begin{proof}
As in the proof of \Cref{th2.2}, we just need to make some normalization such that \eqref{e12.1} holds and the assumptions of \Cref{le2.2} are satisfied. Consider the following transformation (same as \eqref{e2.20}, $\rho_1 \in [1, \rho_2]$ to be specified later):
\begin{equation*}
\begin{aligned}
    \rho_2&=\| u\|_{L^{\infty}(Q_1)}^2+|a|+1, \quad r=\rho_2^{-1/2}, \\
    \tilde{x}&=\frac{x}{r}, \quad \tilde{t}=\frac{t}{\rho_1^{-\gamma}r^2} \quad\text{and} \quad
    \tilde{u}(\tilde x,\tilde t)=\frac{u(x,t)}{\rho_2 r}.
\end{aligned}
\end{equation*}
Then $\tilde{u}$ is a viscosity solution of
\begin{equation*}
    P^{\tilde \varepsilon}_{\tilde a,\tilde\nu}\tilde u=0 \quad\mbox{in}~~Q_1,
\end{equation*}
where
\begin{equation*}
    \tilde{\varepsilon}= \rho_1^{-1} \varepsilon, \quad \tilde{a}=\rho_1^{-1}a \quad  \text{and} \quad \tilde{\nu}=\rho_1^{-1}\rho_2\nu.
\end{equation*}

Clearly, $\|\tilde u\|_{L^{\infty}(Q_1)}\leq 1$. Then as in \Cref{th2.2}, by choosing a proper $\rho_1$, the assumptions of \Cref{le2.2} are satisfied. Hence, $\tilde{u}\in C^{0,1}(\overline{Q^+_{1/2}})$. By transforming back to $u$ and using standard covering arguments, we arrive at the desired estimate.
\end{proof}

%
%
\section{Boundary \texorpdfstring{$C^{1, \alpha}$}{C1,a} regularity for the model problem when
\texorpdfstring{$|a| \gg \nu$}{a>>ν}}\label{S3}
In this section, we prove the boundary $C^{1,\alpha}$ regularity for the viscosity solution of \eqref{e2.1} when $|a|$ is big and we always assume that
\begin{equation*}
1/2\leq |a|\leq 2.
\end{equation*}

We use the classical technique for uniformly parabolic equations. Usually, one prove the boundary $C^{1,\alpha}$ regularity based upon the interior Harnack inequality and the Hopf lemma (e.g., \cite[Theorem 2.1]{MR1139064} and \cite[Theorem 2.8]{LZ_Parabolic}, see also \cite[Lemma 3.1]{MR3246039} and \cite[Lemma 2.12]{MR4682939} for elliptic equations). However, it seems not possible to prove the Harnack inequality for \eqref{e2.1} by the classical method (e.g., \cite[Chapter 4]{MR1351007} and \cite[Section 4]{MR1135923}). Instead, we first prove the strong maximum principle, which is relatively easy (a barrier is enough). Then we obtain the interior Harnack inequality with the aid of the compactness (interior $C^{0,1}$ regularity).

First, we prove the strong maximum principle:
\begin{lemma}[\textbf{Strong maximum principle}]\label{le4.1}
Let $u$ be a nonnegative viscosity supersolution of
\begin{equation*}
    P^{\varepsilon}_{a,\nu} u= 0\quad\mbox{in}~~Q_1.
\end{equation*}
Suppose that $u(0,0)=0$. Then $u\equiv 0$ in $Q_1$.
\end{lemma}
\begin{proof}
We prove the lemma by contradiction. Suppose that the lemma is false. Let
\begin{equation*}
    \Omega\coloneqq \left\{(x,t)\in Q_1: u(x,t)=0\right\},
\end{equation*}
which is a non-empty closed set in $Q_1$. Then there exist $(x_0,t_0)\in \Omega^c$ and $r_1>r_0>0$ such that
\begin{equation*}
Q_{r_0}(x_0,t_0)\subset \Omega^c, \quad  \widetilde Q\coloneqq B_{r_1}(x_0)\times (t_0-r_0^2,t_0]\subset Q_1, \quad (\partial B_{r_1}(x_0)\times \{t_0\})\cap \Omega\neq \varnothing.
\end{equation*}

Let $ Q =\widetilde Q \setminus Q_{r_0}(x_0,t_0)$ and consider the auxiliary function
\begin{equation*}
v(x,t)=c\left(\psi^{\beta}(x-x_0,t-(t_0-r_0^2))-\psi^{\beta}(r_1e_n,r^2_0)\right),
\end{equation*}
where $\psi$ is as in \eqref{e2.8}. Take $\beta$ large enough and then $c$ small enough so that
\begin{equation*}
\left\{\begin{aligned}
   P^{\varepsilon}_{a,\nu}v &<0&& \mbox{in } Q\\
    v&\leq u&& \mbox{on }  \partial_p Q_{r_0}(x_0,t_0)\\
    v&\leq 0&& \mbox{on } \partial_p Q\setminus\partial_p Q_{r_0}(x_0,t_0).
\end{aligned}\right.
\end{equation*}
By the definition of viscosity solution, we have
\begin{equation*}
v\leq u\quad\mbox{in}~~Q.
\end{equation*}
Note also that
\begin{equation*}
    v\leq 0\leq u\quad\mbox{in } \widetilde Q^c\cap \left\{(x,t)\in Q_1: t_0-r_0^2\leq t\leq t_0\right\}.
\end{equation*}
Since $(\partial B_{r_1}(x_0)\times \{t_0\})\cap \Omega\neq \varnothing$, we can choose $(x_1,t_0)\in (\partial B_{r_1}(x_0)\times \{t_0\}) \cap \Omega$. Then by combining the last two inequalities with $v(x_1,t_0)=u(x_1,t_0)=0$, we conclude that $v$ touches $u$ by below at $(x_1,t_0)$. This leads to a contradiction since
$P^{\varepsilon}_{a,\nu} v<0$.
\end{proof}
\begin{remark}\label{re3.2}
The idea of the proof is originated from Hopf \cite{Hopf_1927} and has been used by Nirenberg \cite{MR0055544} to obtain the strong maximum principle for parabolic equations.
\end{remark}
\begin{remark}\label{re3.4}
In the proof above, we use $|a|\geq 1/2$ in an essential way to construct the barrier $v$. Indeed, with this condition, if $Dv$ is small (guaranteed by choosing $c$ small enough), the equation becomes uniformly parabolic. Then we can construct the barrier as usual.
\end{remark}

Based on the strong maximum principle and the $C^{0,1}$ regularity of solutions, we can obtain the following Harnack inequality.
\begin{lemma}[\textbf{Harnack inequality}]\label{le4.2}
Let $u$ be a nonnegative viscosity solution of
\begin{equation*}
    P^{\varepsilon}_{a,\nu} u= 0\quad\mbox{in}~~Q_1.
\end{equation*}
Suppose that $u(0,-1/2)\ge1$ and $\|u\|_{L^{\infty}(Q_1)}\leq 4$. Then
\begin{equation*}
u\geq c>0\quad\mbox{in}~~Q_{1/2},
\end{equation*}
where $c$ is universal.
\end{lemma}
\begin{proof}
Suppose not. Then there exist sequences of $u_m,\varepsilon_m,a_m,\nu_m$ and $(x_m,t_m)\in Q_{1/2}$ such that $u_m$ is a nonnegative viscosity solution of
\begin{equation*}
    P^{\varepsilon_m}_{a_m,\nu_m} u_m= 0\quad\mbox{in } Q_1
\end{equation*}
with
\begin{equation*}
    u_m(0,-1/2)\ge1, \quad \|u_m\|_{L^{\infty}(Q_1)}\leq 4,\quad
1/2\leq |a_m|\leq 2, \quad 0\leq \varepsilon_m\leq 1, \quad 0\leq \nu_m\leq 1
\end{equation*}
and
\begin{equation*}
    u_m(x_m,t_m)\to 0 \quad \text{as } m \to \infty.
\end{equation*}
By the interior $C^{1,\alpha}$ regularity (see \Cref{th2.2}),
\begin{equation*}
    \|u_m\|_{C^{0,1}(\overline{Q_{3/4}})}\leq C \quad \text{for all } m\ge 1,
\end{equation*}
where $C$ is universal. Then there exist subsequences (denoted by $u_m$ etc. again) and $\bar{u}$, $\bar{\varepsilon}$, $\bar{a}$, $\bar{\nu}$, $\bar{x}$, $\bar{t}$ such that
\begin{equation*}
u_m\to \bar{u}\quad\mbox{in}~~L^{\infty}(Q_{3/4}), \quad \varepsilon_m\to \bar\varepsilon
, \quad a_m\to a, \quad \nu_m\to \bar{\nu}, \quad (x_m,t_m)\to (\bar{x},\bar{t})\in \overline{Q_{1/2}}.
\end{equation*}
By the stability of viscosity solutions (see \cite[Theorem 6.1]{MR1443043} or \cite[Proposition 3]{MR2804550}), $\bar{u}$ is a nonnegative viscosity solution of
\begin{equation*}
P^{\bar{\varepsilon}}_{\bar{a},\bar{\nu}}\bar{u}=0\quad\mbox{in}~~Q_{3/4}.
\end{equation*}
Since $\bar{u}(\bar{x},\bar{t})=0$ and $\bar{t} \in [-1/4,0]$, by the strong maximum principle \Cref{le4.1},
\begin{equation*}
  \bar{u}\equiv 0\quad\mbox{in }Q_{3/4}\cap \left\{(x,t): t<-1/4\right\},
\end{equation*}
which contradicts $\bar{u}(0,-1/2)\ge1$.
\end{proof}
\begin{remark}\label{re3.6}
As pointed out by Moser \cite[P. 577]{MR0159138}, the Harnack inequality is quantitative (and hence stronger) version of the strong maximum principle. The \Cref{le4.2} show that with the aid of the compactness, we can derive a quantitative property (the Harnack inequality) from a qualitative property (the strong maximum principle). This demonstrates the power of the compactness and carries the implication why the compactness method is so powerful in the regularity theory.
\end{remark}

Next, we present the Hopf lemma, which is again proved by constructing a barrier.
\begin{lemma}[\textbf{Hopf lemma}]\label{le4.3}
Let $u$ be a nonnegative viscosity solution of \eqref{e2.1}. Suppose that $u(e_n/2,-3/4) \geq 1$ and $\|u\|_{L^{\infty}(Q_1)}\leq 4$. Then
\begin{equation}\label{e.Hopf}
u\geq c x_n\quad\mbox{in}~~Q_{1/2}^+,
\end{equation}
where $c \in (0,1)$ is universal.
\end{lemma}
\begin{proof}
The proof is standard since we have the Harnack inequality. In fact, by applying the Harnack inequality \Cref{le4.2} and noting $u(e_n/2,-3/4) \geq 1$, we have
\begin{equation*}
   u\geq c_0\quad\mbox{in } Q \coloneqq B_{1/4}(e_n/2) \times (-1/4,0],
\end{equation*}
where $c_0>0$ is universal. Let
\begin{equation*}
  v(x,t)=c\left(\psi^{\beta}(x-e_n/2,t-1/4)-\psi^{\beta}(-e_n/2,-1/4)\right),
\end{equation*}
where $\psi$ is as in \eqref{e2.8}. By taking  $\beta$ large enough and $c$ small enough, $v$ satisfies
\begin{equation*}
\left\{\begin{aligned}
    P^{\varepsilon}_{a,\nu}&< 0&&  \mbox{in } Q_{1/2}(e_n/2, 0)\setminus Q\\
    v&\leq c_0&& \mbox{on } \partial_p Q\cap Q_{1/2}(e_n/2,0)\\
    v&\leq 0&& \mbox{on } \partial_p Q_{1/2}(e_n/2,0)\setminus \overline{Q}.
\end{aligned}\right.
\end{equation*}
As before, from the definition of viscosity solution and noting that $v(0,0)=0$ and $v\geq 0$ for $t=0$, we have
\begin{equation*}
  u\geq v\geq cx_n \quad \mbox{on } \{(0',x_n,0) \mid 0< x_n< 1/4\},
\end{equation*}
where $c$ is universal.

By considering $v(x'-x'_0,x_n,t-t_0)$  for $(x'_0,t_0)\in S_{1/2}$ and similar arguments, we obtain
\begin{equation*}
    u\geq cx_n \quad \mbox{in }\left\{(x', x_n, t)\mid (x',0,t)\in S_{1/2} \text{ and } 0<x_n<1/4\right\}.
\end{equation*}
Finally, by the Harnack inequality again,
\begin{equation*}
    u(x,t)\geq cu(e_n/2,-3/4)\geq cx_n \quad \text{for all }(x,t)\in Q^+_{1/2} \cap\{x_n\ge 1/4\}.
\end{equation*}
Therefore, \eqref{e.Hopf} follows.
\end{proof}

Now, we can prove the boundary $C^{1,\alpha}$ regularity.
\begin{lemma}\label{le2.3}
Let $u$ be a viscosity solution of \eqref{e2.1}. Suppose that
\begin{equation*}
a_n=0, \quad 1/2\leq |a|\leq 1, \quad \|Du\|_{L^{\infty}(Q_1^+)}\leq 1.
\end{equation*}
Then $u\in C^{1,\alpha}(0, 0)$, i.e., there exists a constant $A\in [-1,1]$ such that
\begin{equation*}
    |u(x,t)-Ax_n|\leq Cx_n(|x|^{\alpha}+|t|^{\alpha/2})\quad \text{for all } (x,t)\in Q^+_{1/2},
\end{equation*}
where $0<\alpha<1$ and $C$ are universal.
\end{lemma}

\begin{proof}
The proof is standard. It is enough to prove that there exist a nonincreasing sequence $A_m$ and a nondecreasing sequence $B_m$ ($m\geq 0$) such that for all $m\geq 1$,
\begin{equation}\label{e.fbC1a.k}
\begin{aligned}
  &B_mx_n\leq u\leq A_mx_n~~~~\mbox{   in}~~~~Q^+_{2^{-m}},\\
  &0\leq A_m-B_m\leq (1-\mu)(A_{m-1}-B_{m-1}),\\
  &A_0=-B_0=1,
\end{aligned}
\end{equation}
where $0<\mu<1/2$ is universal.

We prove the above by induction. Since $|Du|\leq 1$, \eqref{e.fbC1a.k} holds obviously for $m=1$. Assume that \eqref{e.fbC1a.k} holds for $m$ and we need to prove it for $m+1$. Since the proof is finished when $A_m=B_m$, we consider the case $A_m>B_m$.

Let $r=2^{-m}$. Since \eqref{e.fbC1a.k} holds for $m$, there are two possible cases:
\begin{flalign*}
&\mbox{\textbf{Case 1}:}~~~~u(re_n/2,-3r^2/4)\geq \frac{A_m+B_m}{2}\cdot \frac{r}{2},&&\\
&\mbox{\textbf{Case 2}:}~~~~u(re_n/2,-3r^2/4)< \frac{A_m+B_m}{2}\cdot \frac{r}{2}.&&
\end{flalign*}
Without loss of generality, we suppose that \textbf{Case 1} holds. Let
\begin{equation*}
    \tilde{x}=\frac{x}{r}, \quad
    \tilde t=\frac{t}{r^2} \quad\text{and}\quad
    \tilde u(\tilde x,\tilde t)=\frac{4(u(x,t)-B_mx_n)}{(A_m-B_m)r}.
\end{equation*}
Then $\tilde u$ satisfies
\begin{equation}\label{e3.6}
  \left\{\begin{aligned}
    &P^{\varepsilon}_{\tilde a,\tilde{\nu}}\tilde u=0&& \mbox{in } Q_{1}^+\\
    &0\leq \tilde u\leq 4&& \mbox{in }Q_{1}^+\\
    &\tilde u(e_n/2,-3/4)\geq 1,
  \end{aligned}\right.
\end{equation}
where
\begin{equation*}
    \tilde{a}=a+\nu B_me_n \quad\text{and} \quad
    \tilde{\nu}=\frac{(A_m-B_m)\nu}{4}.
\end{equation*}
Since $a_n=0$ and $|A_m|,|B_m|\leq 1$, we have
\begin{equation*}
\frac{1}{2}\leq |a|\leq |\tilde{a}|\leq |a|+1\leq 2\quad\text{and} \quad \tilde{\nu}\leq 1.
\end{equation*}

By applying \Cref{le4.3} to \eqref{e3.6}, we have
\begin{equation*}\label{e.fbC1a-w}
    \tilde u(\tilde x,\tilde t)\geq c\tilde{x}_n \quad \mbox{for all }(\tilde x,\tilde t) \in Q_{1/2}^+,
\end{equation*}
where $0<c<1/2$ is universal. By rescaling back to $u$,
\begin{equation*}\label{e.fbC1a-w2}
    u(x,t)\geq (B_m+\mu(A_m-B_m))x_n \quad \mbox{in } Q_{r/2}^+,
\end{equation*}
where $\mu=c/4$. Let $A_{m+1}=A_m$ and $B_{m+1}=B_m+\mu(A_m-B_m)$. Then
\begin{equation*}\label{e.reg-ak}
A_{m+1}-B_{m+1}=(1-\mu)(A_m-B_m).
\end{equation*}
Hence, \eqref{e.fbC1a.k} holds for $m+1$ and the proof is completed by induction.
\end{proof}

With the assistance of the small perturbation regularity \Cref{L.K2-2}, we have the following higher regularity.
\begin{lemma}\label{le3.4}
Let $u$ be a viscosity solution of \eqref{e2.1}. Suppose that
\begin{equation*}
a_n=0, \quad 1/2\leq |a|\leq 1, \quad \|Du\|_{L^{\infty}(Q_1^+)}\leq 1.
\end{equation*}
Then $u\in C^{1,\beta}(0, 0)$ for any $\beta \in (0,1)$, i.e., there exists a constant $A\in[-1,1]$ such that
\begin{equation*}
    |u(x,t)-Ax_n|\leq C(|x|^{1+\beta}+|t|^{\frac{1+\beta}{2}})\quad \text{for all } (x,t)\in Q^+_{1/2}
\end{equation*}
where $C$ depends only on $n$, $p$, $\gamma$ and $\beta$.
\end{lemma}

\begin{proof}
By \Cref{le2.3}, $u\in C^{1,\alpha}(0,0)$ for some $0<\alpha<1$. That is, there exists $A$ such that
\begin{equation*}
    |u(x,t)-Ax_n|\leq Cx_n(|x|^{\alpha}+|t|^{\alpha/2})\quad \text{for all } (x,t)\in Q^+_{1/2}.
\end{equation*}
Consider the following transformation for $0<r<1/2$:
\begin{equation*}
    \tilde{x}=\frac{x}{r}, \quad \tilde{t}=\frac{t}{r^2} \quad\text{and}\quad
    \tilde{u}(\tilde x,\tilde t)=\frac{u(x,t)-Ax_n}{r}.
\end{equation*}
Then $\tilde{u}$ is a solution of
\begin{equation*}
\left\{\begin{aligned}
    P^{\varepsilon}_{\tilde{a},\nu} \tilde u&= 0&&  \mbox{in } Q_{1}^+\\
\tilde u&=0&& \mbox{on } S_1,
\end{aligned}\right.
\end{equation*}
where $\tilde{a}=a+{\nu}Ae_n$. Clearly,
\begin{equation*}
 1/2\leq |a|\leq  |\tilde{a}|\leq |a|+|A|\leq 2.
\end{equation*}

For any $\beta \in (0,1)$, let $0<\eta<1$ be the constant from \Cref{L.K2-2} with
\begin{equation*}
f\equiv 0, \quad g\equiv 0 \quad\text{and}\quad \Omega=Q_{1}^+.
\end{equation*}
Take $r$ small enough such that
\begin{equation*}
\|\tilde{u}\|_{L^{\infty}(Q_1^+)}\leq Cr^{\alpha}\leq \eta.
\end{equation*}
Then from \Cref{L.K2-2}, $\tilde{u}\in C^{1,\beta}(0, 0)$. By rescaling back to $u$, we have $u\in C^{1,\beta}(0, 0)$.
\end{proof}

%
%
\section{Boundary estimates for the model problem when \texorpdfstring{$|a| \ll \nu$}{|a|<<ν}}\label{S4}
In this section, we prove some boundary estimates for smooth solutions of
\begin{equation}\label{e2.0}
	\left\{\begin{aligned}
		P^{\varepsilon}_{a}u& =0&& \mbox{in } Q_1^+ \\
		u&=0&& \mbox{on } S_1,
	\end{aligned}\right.
\end{equation}
when $|a|$ is small. Note that $P^{\varepsilon}_{a}=P^{\varepsilon}_{a,1}$ in \eqref{e2.0}. Throughout this section, we always assume
\begin{equation*}
    0<\varepsilon\leq 1, \quad
    0\leq |a|\leq 1/2 \quad\text{and}\quad  \|Du\|_{L^{\infty}(Q_1^+)}\leq 1
\end{equation*}
unless stated otherwise. To obtain the gradient estimates in this section, we need to find the equation satisfied by the gradient of $u$. For this purpose, we regularize the original equation with the approximation parameter $\varepsilon$, which allows us to deal with smooth solutions. In addition, since $|Du+a|=0$ may occur, we require that $\varepsilon$ is strictly positive.


Since $a$ is small (may be $0$), we cannot use the classical method for uniformly parabolic equations as in the last section. In fact, for this degenerate/singular case, there are only interior Harnack inequalities in some weak forms (see \cite{MR2865434}, \cite[Theorem 2.1, Chapter VI and Theorem 1.1, Chapter VII]{MR1230384}), which are not adequate for the boundary $C^{1,\alpha}$ regularity. Furthermore, the strong maximum and the Hopf lemma fail in general (see the counterexamples in \cite[Section 4]{MR3217145}).

In this section, we follow closely the strategy of Imbert, Jin and Silvestre \cite{MR3918407} (see also \cite{MR4676644} for fully nonlinear equations). The idea is that we prove the $C^{1,\alpha}$ estimate according to two cases: non-degenerate and degenerate.
\begin{itemize}
    \item If $Du$ is close to a unit vector (non-degenerate), we can use uniformly parabolic equation theory (or small perturbation regularity theory) to obtain the $C^{1,\alpha}$ estimate.
    \item Otherwise, $|Du|\leq l$ for some constant $l<1$ in a set with a positive measure (degenerate). Then by the weak Harnack inequality (for uniformly parabolic equations), $|Du|\leq 1-\delta$ for some $0<\delta<1$ in a smaller scale, that is, $|Du|$ has a decay. By iteration, we have the $C^{1,\alpha}$ estimate.
\end{itemize}

We remark here that this strategy (considering the equation in degenerate/non-degenerate cases separately) can be tracked to Uhlenbeck (see \cite[Section 5, Proposition 5.1]{MR474389}) and has been widely used for $p$-Laplace type equations, e.g., \cite[Section 2 and Section 3]{MR672713}, \cite[Proposition 4.1 and Proposition 4.2]{MR709038}, \cite[Proof of Proposition 3]{MR727034}, \cite[Section 4]{MR1264526} and \cite[Section 4]{MR1206157} etc.

The first lemma concerns the non-degenerate case.
\begin{lemma} \label{le3.6}
Let $u$ be a smooth solution of \eqref{e2.0}. For any $0<\eta<1$, there exist $\varepsilon_0,\varepsilon_1>0$ (small enough) depending only on $n$, $p$, $\gamma$ and $\eta$ such that if
\begin{equation*}
	| \{(x,t)\in Q^+_1 : |Du(x,t)-e_n| > \varepsilon_0 \}| \leq \varepsilon_1,
\end{equation*}
then
\begin{equation}\label{e4.21}
	|u(x,t)-x_n| \leq \eta \quad \text{for all } (x,t) \in Q^+_{1/2}.
\end{equation}

Similarly, if
\begin{equation*}
	| \{(x,t)\in Q^+_1 : |Du(x,t)+e_n| > \varepsilon_0 \}| \leq \varepsilon_1,
\end{equation*}
we have
\begin{equation}\label{e4.22}
	|u(x,t)+x_n| \leq \eta \quad \text{for all } (x,t) \in Q^+_{1/2}.
\end{equation}
\end{lemma}
\begin{proof}
We essentially follow the proof of \cite[Lemma 4.6]{MR3918407}, which becomes simpler due to the boundary condition $u=0$ on $S_1$. Let
\begin{equation*}
    f(t)\coloneqq | \{ x \in B_1^+ : |Du(x,t)-e_n| > \varepsilon_0 \} | \quad \text{for $t\in (-1,0]$}.
\end{equation*}
By the assumptions and Fubini's theorem,
\begin{equation*}
\int_{-1}^0 f(t) \, dt=| \{(x,t)\in Q^+_1 : |Du(x,t)-e_n| > \varepsilon_0 \}| \leq \varepsilon_1.
\end{equation*}
Let $E\coloneqq \{t\in (-1,0]: f(t)\ge\sqrt{\varepsilon_1}\}$. Then for $t\in (-1,0] \setminus E$, by the Morrey's inequality and noting $u=0$ on $S_1$ and $|Du|\leq 1$, we have
\begin{equation}\label{e4.19}
\|u(\cdot,t)-x_n\|_{L^{\infty}(B_{1/2}^+)}\leq C\|Du(\cdot,t)-e_n\|_{L^{2n}(B_{1/2}^+)}\leq C(\varepsilon_0+\varepsilon_1^{1/(4n)}),
\end{equation}
where $C$ depends only on $n$.

On the other hand,
\begin{equation*}
|E|\le \frac{1}{\sqrt{\varepsilon_1}}\int _{E} f(t) \, dt\le \frac{1}{\sqrt{\varepsilon_1}}\int _{-1}^0 f(t) \, dt\le \sqrt{\varepsilon_1}.
\end{equation*}
Hence, for any $s\in E$, there exists $t\in (-1, 0] \setminus E$ such that $|t-s|\leq \sqrt{\varepsilon_1}$. With the aid of the $C^{0,1}$ regularity \Cref{le2.2},
\begin{equation}\label{e4.20}
  \begin{aligned}
\|u(\cdot,s)-x_n\|_{L^{\infty}(B_{1/2}^+)}\leq& \|u(\cdot,s)-u(\cdot,t)\|_{L^{\infty}(B_{1/2}^+)} +\|u(\cdot,t)-x_n\|_{L^{\infty}(B_{1/2}^+)}\\
\leq&  C(|s-t|^{1/2}+\varepsilon_0+\varepsilon_1^{1/(4n)})\leq
C(\varepsilon_1^{1/4}+\varepsilon_0+\varepsilon_1^{1/(4n)}),\\
  \end{aligned}
\end{equation}
where $C$ is universal. By \eqref{e4.19}, \eqref{e4.20} and choosing $\varepsilon_0,\varepsilon_1$ small enough, we arrive at the conclusion.
\end{proof}

If \eqref{e4.21} or \eqref{e4.22} holds, the $C^{1,\alpha}$ regularity follows from the small perturbation regularity (see \Cref{le3.7} for details). Next, we move our attention to the degenerate case, which is more difficult. The main difficulty in applying the strategy of \cite{MR3918407}, which addresses the interior regularity, to the boundary regularity setting lies in the lack of information regarding the values of $u_n$ on $S_1$. This difficulty can be overcome by demonstrating that $u_n$ has a strict decay on $S_{1/2}$, if it is small on a set with a positive measure. This idea was inspired by Lieberman \cite[Lemma 1.2]{MR1207530}.
\begin{lemma}\label{le3.2}
Let $u$ be a smooth solution of \eqref{e2.0}. Suppose that
\begin{equation*}
	|\{(x,t)\in Q^+_1: u_n(x,t)\leq l\}|> \mu |Q^+_1|
\end{equation*}
for some $3/4<l<1$ and $\mu>0$. Then
\begin{equation*}
u_n\leq l_0\quad\mbox{on }S_{1/2},
\end{equation*}
where $l_0 \in [l,1)$ is a constant depending only on $n$, $p$, $\gamma$, $l$ and $\mu$.

Similarly, if
\begin{equation*}
	|\{(x,t)\in Q^+_1: -u_n(x,t)\leq l\}|> \mu |Q^+_1|,
\end{equation*}
we have
\begin{equation*}
	-u_n\leq l_0\quad\mbox{on }S_{1/2}.
\end{equation*}
\end{lemma}

%

\begin{proof}
	First, since $\|Du\|_{L^{\infty}(Q_1^+)}\leq 1$ and $u=0$ on $S_1$,
\begin{equation}\label{e3.4}
	u(x,t)=\int_{0}^{x_n}u_n(x',s,t) \, ds\leq x_n \quad\mbox{in } Q_1^+.
\end{equation}
Next, by \cite[Lemma 4.1]{MR3918407}, there exist $\tau_1=\tau_1(n,\mu) \in (0,1/4)$, $\tau=\tau(n, p,  \gamma, l, \mu) \in (0,\tau_1]$ and $\delta_1=\delta_1 (n, p, \gamma, l, \mu) \in (0,1)$  such that
\begin{equation} \label{e3.3}
	u_n\leq 1-\delta_1\quad\mbox{in } \widetilde{Q} \coloneqq \{(x,t) \in Q_1^+ : |x'|<\tau, \, |x_n-1/2|<\tau, \, -\sigma <t \leq 0 \}
\end{equation}
and $\widetilde{Q} \subset Q_{\tau_1}$, where $\sigma\coloneqq (1-\delta_1)^{-\gamma}\tau^2$. Let
\begin{equation*}
	\mathcal{Q}\coloneqq \{(x,t) \in Q_1^+ : |x'|<\tau, \, |x_n-1/2|<\tau/2, \,  -\sigma <t \leq0 \}.
\end{equation*}
Then for any $(x,t)\in \overline{\mathcal{Q}}$, by \eqref{e3.3},
\begin{equation}\label{e3.5}
\begin{aligned}
	u(x,t)&=\int_{0}^{1/2-\tau}u_n(x',s,t)\,ds + \int_{1/2-\tau}^{x_n}u_n(x',s,t)\,ds\\
	&\leq \frac{1}{2} -\tau +(1-\delta_1)\left(x_n-\frac{1}{2} + \tau\right)\\
	&\leq x_n - \frac{1}{2} \tau \delta_1.\\
  \end{aligned}
\end{equation}
For $\psi$ defined in \eqref{e2.8}, if we set
\begin{equation*}
	v(x,t)=x_n-\varphi(x,t) \quad \text{and} \quad
	\varphi(x,t)=c\left(\psi^{\beta}(x-e_n/2,t+\sigma)-\psi^\beta(-e_n/2,\sigma)\right),
\end{equation*}
then it is easy to check that
\begin{equation*}
	\varphi(0,0)=0 \quad\text{and}\quad
    \varphi\leq 0 \quad\mbox{on } \partial_p Q,
\end{equation*}
where $Q\coloneqq B_{1/2}(e_n/2)\times (-\sigma,0]$. Hence, by combining with \eqref{e3.4},
\begin{equation*}
	u\leq v\quad\mbox{on } \partial_p Q.
\end{equation*}
As before, with the aid of \eqref{e3.5}, by taking $\beta$ large enough and $c$ small enough, we have
\begin{equation*}
\left\{\begin{aligned}
    P^{\varepsilon}_{a}v &> 0&& \mbox{in } Q \setminus \overline{\mathcal{Q}}\\
    v&\geq u&& \mbox{on }\partial_p \mathcal{Q}\\
    v&\geq u&& \mbox{on }\partial_p Q
\end{aligned}\right.
\end{equation*}
and so
\begin{equation*}
	u\leq v\quad\mbox{in } Q\setminus \overline{\mathcal{Q}}.
\end{equation*}
Since $u(0,0)=v(0,0)=0$,
\begin{equation*}
	u_n(0,0)\leq v_n(0,0)=1-\delta_2.
\end{equation*}
By considering any other point $(x_0,t_0)\in S_{1/2}$ and similar arguments, we have
\begin{equation*}
	u_n\leq 1-\delta_3\quad\mbox{on } S_{1/2}.
\end{equation*}
Therefore, we obtain the conclusion by choosing $l_0=\max\{1-\delta_3,l\}$.
\end{proof}
\begin{remark}\label{re3.1}
The proof is mainly inspired by \cite[Lemma 1.2]{MR1207530}. That is, we first show the decay in the interior and then construct a barrier to obtain the decay on the boundary.
\end{remark}

Now, we can use the technique from \cite[Lemma 4.1]{MR3918407} to prove a decay for $Du$.
\begin{lemma}\label{le3.3}
	Let $u$ be a smooth solution of \eqref{e2.0}. Suppose that
\begin{equation*}
	|\{(x,t)\in Q^+_1:u_n(x,t)\leq l\}|> \mu |Q^+_1|
\quad\text{and} \quad
	|\{(x,t)\in Q^+_1:-u_n(x,t)\leq l\}|> \mu |Q^+_1|
\end{equation*}
for some $3/4<l<1$ and $\mu>0$. Then
\begin{equation}\label{e3.2}
|Du|\leq 1-\delta\quad\mbox{in }Q^{(1-\delta)+}_{\tau},
\end{equation}
where $\tau,\delta\in (0,1/4)$ are constants depending only on $n$, $p$, $\gamma$, $l$ and $\mu$.
\end{lemma}

\begin{proof}
The proof is almost the same as that of \cite[Lemma 4.1]{MR3918407}. The only difference is the definition of $w$ (see Line 6, Page 853 in \cite{MR3918407}). Given a unit vector $b=(b_1,\dots,b_n)$, without loss of generality, we may assume $b_n\geq 0$. Define
\begin{equation*}
l_1=\frac{1+l_0}{2}, \quad \rho=\frac{1-l_0}{4}
\end{equation*}
and
\begin{equation}\label{e4.23}
w=\left(Du\cdot b-l_1+\rho |Du|^2\right)^+,
\end{equation}
where $l_0$ is from \Cref{le3.2}. Then we have
\begin{equation*}
|Du|>3/4 \quad\mbox{in }\Omega_+\coloneqq \left\{(x,t)\in Q_{1}^+: w(x,t)>0\right\}
\end{equation*}
and so
\begin{equation}\label{e4.24}
|Du+a|\geq 1/4  \quad\mbox{in }\Omega_+.
\end{equation}
Hence, we differentiate the equation \eqref{e2.0} and proceed as in the proof of \cite[Lemma 4.1]{MR3918407} to show that $w$ is a subsolution of some uniformly parabolic equation (see Line 13, Page 853 in \cite{MR3918407}). That is,
\begin{equation}\label{e4.5}
w_t\le a_{ij} w_{ij} + c_1|D w|^2\quad\mbox{in }\Omega_+,
\end{equation}
where $c_1>0$ is a constant depending only on $n$, $p$, $\gamma$ and $l$. From \Cref{le3.2}, $u_n\leq l_0$ on $S_{1/2}$. By combining with $0\leq b_n\leq1$ and $u_i=0$ on
$S_{1/2}$ ($1\leq i\leq n-1$), we have
\begin{equation*}
 w=0\quad\mbox{on}~~S_{1/2}.
\end{equation*}
Then we can take the zero extension of $w$ to $Q_{1/2}^-\coloneqq Q_{1/2}\cap \{x_n<0\}$ such that $w$ is a viscosity subsolution of \eqref{e4.5} in $Q_{1/2}$. The rest of the proof is the same as that of \cite[Lemma 4.1]{MR3918407} and we omit it. In conclusion, we have
\begin{equation*}
Du\cdot b< 1-\delta\quad\mbox{in }Q^{(1-\delta)+}_{\tau},
\end{equation*}
where $\tau,\delta\in (0,1/4)$ depend only on $n$, $p$, $\gamma$, $l$ and $\mu$. Since $b$ is arbitrary, we obtain \eqref{e3.2}.
\end{proof}
\begin{remark}\label{re4.2}
The technique that differentiates the equation and considers an auxiliary function $w$ like \eqref{e4.23}, which is a subsolution of some linear uniformly parabolic equation after some calculation, is due to Lady\v zenskaya and Ural\cprime tseva \cite{MR173855} (see also \cite[Chapter 6.1]{MR244627}, \cite[Chapter VI.1]{MR241822} and \cite[Chapter 13.3]{MR1814364}). The innovation of \cite{MR3918407} is that it can deal with non-uniformly parabolic equations. To be more precise, to prove the decay of $|Du|$, we consider two cases. If $|Du|$ is small, this is just what we want. If $|Du|$ is big, the equation becomes uniformly parabolic, then we can use the technique of Lady\v zenskaya and Ural\cprime tseva.
\end{remark}

\begin{remark}\label{re4.3}
For \eqref{e2.0}, if \eqref{e4.24} holds, the equation becomes uniformly parabolic. Then we can use the technique from \cite{MR3918407}. Hence, the smallness of $|a|$ is used in an essential way.
\end{remark}

As a corollary, we have the following scaling version of \Cref{le3.3}. Since its proof is exactly the same as that of \cite[Corollary 4.2]{MR3918407}, we omit it.
\begin{lemma}\label{le3.5}
Let $l$, $\mu$, $\delta$ and $\tau$ be as in \Cref{le3.3} and let $u$ be a smooth solution of \eqref{e2.0}. For any $k\geq 0$ satisfying
\begin{equation}\label{e4.6-0}
    k\leq \min\left\{\frac{\ln (2\varepsilon)}{\ln (1-\delta)}, \frac{\ln (2|a|)}{\ln (1-\delta)}\right\},
\end{equation}
\begin{equation}\label{e4.6-1}
\left|\{(x,t)\in Q^{(1-\delta)^i+}_{\tau^i}:u_n(x,t)\leq l(1-\delta)^i\}\right|> \mu \left|Q^{(1-\delta)^i+}_{\tau^i}\right| \quad \text{for all }  i=0,1,\cdots, k
\end{equation}
and
\begin{equation}\label{e4.6}
\left|\{(x,t)\in Q^{(1-\delta)^i+}_{\tau^i}:-u_n(x,t)\leq l(1-\delta)^i\}\right|> \mu \left|Q^{(1-\delta)^i+}_{\tau^i}\right| \quad \text{for all }  i=0,1,\cdots, k,
\end{equation}
we have
\begin{equation*}
|Du|\leq (1-\delta)^{i+1}\quad\mbox{in }Q^{(1-\delta)^{i+1}+}_{\tau^{i+1}} \quad \text{for all } i =0,1,\cdots, k.
\end{equation*}
\end{lemma}

%
%
\section{Boundary \texorpdfstring{$C^{1,\alpha}$}{C1,a} regularity for the model problem}\label{S6}
Based on the estimates derived in last two sections and the small perturbation regularity (see \Cref{L.K2-2}), we can finally drive the boundary $C^{1,\alpha}$ estimate for smooth solutions of the model problem \eqref{e2.0}. Then by an approximation, we obtain the boundary $C^{1,\alpha}$ regularity for viscosity solutions.
\begin{lemma}\label{le3.7}
Let $u$ be a smooth solution of \eqref{e2.0} with
\begin{equation*}
  0<\varepsilon\leq 1 , \quad  a_n=0, \quad  |a|\leq 1 \quad\text{and}\quad \|Du\|_{L^{\infty}(Q_1^+)}\leq 1.
\end{equation*}
Then $u\in C^{1,\bar\alpha}_{\gamma}(0,0)$, i.e., there exists $A\in [-1,1]$ such that
\begin{equation}\label{e5.5}
	|u(x,t)-Ax_n|\leq  C (|x|^{1+\bar\alpha}+|t|^{\frac{1+\bar\alpha}{2-\bar\alpha\gamma}}) \quad \text{for all } (x,t)\in Q_{1/2}^+,
\end{equation}
where $0<\bar\alpha< \min\{1/2,1/2(1+\gamma)\}$ and $C>0$ are universal.
\end{lemma}
\begin{remark}\label{re2.1}
The constant $\bar{\alpha}$ may be different from the one in \Cref{th2.1} for the interior regularity. We can choose the smaller one such that both \Cref{th2.1} and \Cref{le3.7} hold for the same $\bar{\alpha}$. From now on, $\bar{\alpha}$ is fixed throughout this paper.
\end{remark}

\begin{proof}
We prove the lemma following the outline of \cite[Theorem 4.8]{MR3918407}, but there exists one additional case due to the presence of $a$. We first determine various constants. Let $0<\eta<1$ be the constant from \Cref{L.K2-2} subjected to
\begin{equation*}
    \beta=1/2 \quad\text{and}\quad \Omega=Q_{1/2}^+.
\end{equation*}
Hence, $\eta$ is universal. For this $\eta$, we fix $\varepsilon_0,\varepsilon_1>0$ such that \Cref{le3.6} holds with them. Set
\begin{equation*}
l=1-\varepsilon_0^2/2\quad\text{and} \quad \mu=\varepsilon_1/|Q_1^+|.
\end{equation*}
As in the proof of \cite[Theorem 4.8]{MR3918407}, if
\begin{equation*}
|\{(x,t)\in Q^+_1:u_n(x,t)\leq l\}|\leq \mu |Q^+_1|~~ (\mbox{or}~~|\{(x,t)\in Q^+_1:-u_n(x,t)\leq l\}|\leq \mu |Q^+_1|),
\end{equation*}
then we have
\begin{equation}\label{e4.8}
| \{(x,t)\in Q^+_1 : |Du-e_n| > \varepsilon_0 \}| \leq \varepsilon_1~~ (\mbox{or}~~| \{(x,t)\in Q^+_1 : |Du+e_n| > \varepsilon_0 \}| \leq \varepsilon_1).
\end{equation}
Let $\tau$, $\delta$ be the constants in \Cref{le3.3} depending on $l$ and $\mu$. From the choice of $l$, $\mu$, we know that $\tau$, $\delta$ are universal, i.e., they depend only on $n$, $p$ and $\gamma$. Moreover, we choose $\tau$ small enough such that
\begin{equation*}
\tau<(1-\delta)^{1+\gamma}.
\end{equation*}
Finally, we take $0<\bar\alpha<1/(1+\gamma)$ determined by
\begin{equation*}
\tau^{\bar\alpha}=1-\delta.
\end{equation*}

Let $k_0\geq 1$ be the smallest integer such that \eqref{e4.6-0}-\eqref{e4.6} hold for all $k\leq k_0-1$ but one of them does not hold for $k_0$. By \Cref{le3.5} with $k=k_0-1$,
\begin{equation}\label{e4.18}
|Du|\leq \tau^{i\bar\alpha}\quad\mbox{in }Q^{\tau^{i\bar\alpha}+}_{\tau^{i}} \quad \text{for all } i =0,1,\cdots, k_0,
\end{equation}
which implies
\begin{equation}\label{e4.9}
    |Du(x,t)|\leq C(|x|^{\bar{\alpha}}+|t|^{\frac{\bar{\alpha}}{2-\bar{\alpha}\gamma}})\quad\mbox{in }
Q_1^+\setminus Q^{\tau^{(k_0+1)\bar\alpha}+}_{\tau^{k_0+1}}.
\end{equation}

In the following, we prove the lemma according to three cases:

\noindent\textbf{Case 1:} \eqref{e4.6-0} fails for $k_0$ with
\begin{equation}\label{e10.13}
k_0-1\leq \frac{\ln (2\varepsilon)}{\ln (1-\delta)}=\min\left\{\frac{\ln (2\varepsilon)}{\ln (1-\delta)}, \frac{\ln (2|a|)}{\ln (1-\delta)}\right\}<k_0.
\end{equation}
Introduce the following transformation:
\begin{equation}\label{e4.10}
r=\tau^{k_0}, \quad \tilde{x}=\frac{x}{r}, \quad \tilde{t}=\frac{t}{r^{2-\bar\alpha\gamma}}
    \quad\text{and}\quad
    \tilde{u}(\tilde{x},\tilde{t})=\frac{u(x,t)}{r^{1+\bar\alpha}}.
\end{equation}
Then $\tilde{u}$ is a smooth solution of
\begin{equation}\label{e4.15}
	\left\{\begin{aligned}
		P^{\tilde\varepsilon}_{\tilde a}\tilde{u}& =0&& \mbox{in } Q_1^+ \\
		\tilde{u}&=0&& \mbox{on } S_1,
	\end{aligned}\right.
\end{equation}
where
\begin{equation*}
\tilde{\varepsilon}=\frac{\varepsilon}{r^{\bar{\alpha}}}\quad\text{and}\quad \tilde{a}=\frac{a}{r^{\bar{\alpha}}}.
\end{equation*}
By \eqref{e4.18}, we have $|D\tilde{u}|\leq 1$ in $Q_1^+$. In addition, by the property of $k_0$ (see \eqref{e10.13}), we have
\begin{equation*}
  1/2< \tilde{\varepsilon} < 1 \quad\text{and}\quad  |\tilde{a}| < 1.
\end{equation*}
Thus, $\tilde{u}$ satisfies a quasi-linear uniformly parabolic equation with smooth coefficients. By the interior regularity (see \cite[Theorem 4.4, P. 560]{MR241822}), there exists $\tilde{A}\in [-1,1]$ such that
\begin{equation*}
    |D\tilde{u}(\tilde{x},\tilde{t})-\tilde{A}e_n|\leq C(|\tilde x|+|\tilde t|^{1/2})\leq C(|\tilde x|^{\bar{\alpha}}+|\tilde t|^{\frac{\bar{\alpha}}{2-\bar{\alpha}\gamma}})
\quad\mbox{in } Q_{\tau}^{\tau^{\bar\alpha}+}\subset Q_{1/4}^+,
\end{equation*}
where $C$ is universal. Here, we have used
\begin{equation*}
\frac{\bar{\alpha}}{2-\bar{\alpha}\gamma}\leq \frac{1}{2},
\end{equation*}
which holds by the choice of $\bar{\alpha}$.

By rescaling back to $u$,
\begin{equation*}
    |Du(x,t)-A|\leq C(|x|^{\bar{\alpha}}+|t|^{\frac{\bar{\alpha}}{2-\bar{\alpha}\gamma}})
\quad\mbox{in } Q^{\tau^{(k_0+1)\bar{\alpha}}+}_{\tau^{k_0+1}},
\end{equation*}
where
\begin{equation*}
    A=\tau^{k_0\bar{\alpha}}\tilde{A}.
\end{equation*}
By combining with \eqref{e4.9}, we have
\begin{equation}\label{e4.13}
    |Du(x,t)-A|\leq C(|x|^{\bar{\alpha}}+|t|^{\frac{\bar{\alpha}}{2-\bar{\alpha}\gamma}})
\quad\mbox{in } Q_1^+.
\end{equation}

\noindent \textbf{Case 2:} \eqref{e4.6-0} fails for $k_0$ with
\begin{equation*}
k_0-1\leq \frac{\ln (2|a|)}{\ln (1-\delta)}=\min\left\{\frac{\ln (2\varepsilon)}{\ln (1-\delta)}, \frac{\ln (2|a|)}{\ln (1-\delta)}\right\}<k_0.
\end{equation*}
Take the same transformation as \eqref{e4.10}. Then $\tilde{u}$ satisfies \eqref{e4.15} with
\begin{equation*}
0<\tilde{\varepsilon}< 1, \quad \tilde{a}_n=0 \quad\text{and}\quad  1/2 <|\tilde{a}|< 1.
\end{equation*}

From \Cref{le3.4}, $\tilde{u}\in C^{1,1/2}(0,0)$. By transforming back to $u$ and combining with \eqref{e4.9} (note $\bar{\alpha}<1/2$), we have $u\in C_{\gamma}^{1,\bar{\alpha}}(0,0)$. That is, there exists $A\in [-1,1]$ such that
\begin{equation}\label{e5.3}
	|u(x,t)-Ax_n|\leq  C (|x|^{1+\bar\alpha}+|t|^{\frac{1+\bar\alpha}{2-\bar\alpha\gamma}}) \quad \text{for all } (x,t)\in Q_{1/2}^+.
\end{equation}

\noindent \textbf{Case 3:} \eqref{e4.6-1} or \eqref{e4.6} fails for $k_0$. Without loss of generality, we assume that \eqref{e4.6-1} fails (the other case can be treated similarly). That is,
\begin{equation}\label{e4.17}
    \left|\{(x,t)\in Q^{\tau^{k_0\bar{\alpha}}+}_{\tau^{k_0}}:u_n(x,t)\leq l\tau^{k_0\bar{\alpha}}\}\right|\leq \mu \left|Q^{\tau^{k_0\bar{\alpha}}+}_{\tau^{k_0}}\right|.
\end{equation}
Under the transformation \eqref{e4.10} again, $\tilde{u}$ satisfies \eqref{e4.15} with
\begin{equation*}
 0<\tilde{\varepsilon}< 1/2 \quad\text{and}\quad   |\tilde{a}|< 1/2.
\end{equation*}
Moreover, \eqref{e4.17} is equivalent to
\begin{equation*}
\left|\{(\tilde x,\tilde t)\in Q^{+}_{1}:\tilde u_n(\tilde x,\tilde t)\leq l\}\right|\leq \mu \left|Q_1^+\right|,
\end{equation*}
which implies that \eqref{e4.8} holds for $\tilde{u}_n$. From \Cref{le3.6},
\begin{equation*}
	|\tilde u(\tilde x,\tilde t)-\tilde x_n| \leq \eta \quad \text{for all } (\tilde x,\tilde t) \in Q^+_{1/2}.
\end{equation*}
Consider $\bar{u}(\tilde x,\tilde t)=\tilde u(\tilde x,\tilde t)-\tilde x_n$. Then $\bar{u}$ is a solution of
\begin{equation*}
\left\{\begin{aligned}
    P^{\tilde \varepsilon}_{\tilde a+e_n} \bar u&= 0&& \mbox{in }Q_{1}^+\\
\bar u&=0&& \mbox{on }S_1.
\end{aligned}\right.
\end{equation*}
Since $\tilde{a}_n=0$ and $|\tilde{a}|<1/2$, we have
\begin{equation*}
1\leq |\tilde a+e_n|\leq 2.
\end{equation*}
Then by \Cref{L.K2-2} with
\begin{equation*}
\nu=1, \quad \beta=1/2, \quad  f\equiv 0, \quad g\equiv 0, \quad \Omega=Q_1^+,
\end{equation*}
we have $\bar{u}\in C^{1,1/2}(0,0)$. By transforming back to $u$ and combining with \eqref{e4.9}, we have
$u\in C_{\gamma}^{1,\bar{\alpha}}(0,0)$ and for some $A\in [-1,1]$,
\begin{equation}\label{e5.4}
	|u(x,t)-Ax_n|\leq  C (|x|^{1+\bar\alpha}+|t|^{\frac{1+\bar\alpha}{2-\bar\alpha\gamma}}) \quad \text{for all } (x,t)\in Q_{1/2}^+.
\end{equation}

Finally, by combining \eqref{e4.13}, \eqref{e5.3} and \eqref{e5.4}, we conclude that $u\in C_{\gamma}^{1,\bar\alpha}(0,0)$ with \eqref{e5.5}.
\end{proof}

Up to a normalization, we have
\begin{lemma}\label{le5.1}
Let $u$ be a smooth solution of \eqref{e2.0} with $0<\varepsilon\leq 1$ and $a_n=0$. Then $u\in C^{1,\bar\alpha}_{\gamma}(0,0)$, i.e., there exists $A\in \mathbb{R}$ such that
\begin{equation}\label{e5.6}
	|u(x,t)-Ax_n|\leq  C(|x|^{1+\bar\alpha}+|t|^{\frac{1+\bar\alpha}{2-\bar\alpha\gamma}}) \quad \text{for all } (x,t)\in Q_{1/2}^+,
\end{equation}
and
\begin{equation}\label{e5.7}
|A|\leq C,
\end{equation}
where $\bar\alpha$ is as in \Cref{le3.7} and $C>0$ depends only on $n$, $p$, $\gamma$, $|a|$ and $\|u\|_{L^{\infty}(Q_1^+)}$.
\end{lemma}
\begin{proof}
By the global $C^{0,1}$ estimate \Cref{co3.1},
\begin{equation*}
\|Du\|_{L^{\infty}(Q_{1/2}^+)}\leq K,
\end{equation*}
where $K$ depends only on $n$, $p$, $\gamma$, $|a|$ and $\|u\|_{L^{\infty}(Q_1^+)}$. Consider the following transformation:
\begin{equation*}
  \begin{aligned}
&r=1/2, \quad\rho=K+|a|+1, \quad
\tilde{x}=\frac{x}{r}, \quad \tilde{t}=\frac{t}{\rho^{-\gamma}r^2} \quad\text{and}\quad
 \tilde{u}(\tilde{x},\tilde{t})=\frac{u(x,t)}{\rho r}.\\
  \end{aligned}
\end{equation*}
Then $\tilde{u}$ is a solution of
\begin{equation*}
	\left\{\begin{aligned}
		P^{\tilde \varepsilon}_{\tilde a}\tilde u& =0&& \mbox{in } Q_1^+ \\
		\tilde u&=0&& \mbox{on } S_1,
	\end{aligned}\right.
\end{equation*}
where
\begin{equation*}
\|D\tilde{u}\|_{L^{\infty}(Q_1^+)}=\rho^{-1}\|Du\|_{L^{\infty}(Q_{1/2}^+)} \leq 1, \quad
|\tilde{a}|=\rho^{-1}|a|\leq 1 \quad\text{and}\quad 0<\tilde{\varepsilon}=\rho^{-1}\varepsilon\leq 1.
\end{equation*}
Thus, $\tilde{u}$ satisfies the assumptions of \Cref{le3.7}. Then $\tilde{u}\in C^{1,\bar{\alpha}}_{\gamma}(0,0)$. By transforming back to $u$, we obtain $u\in C^{1,\bar{\alpha}}_{\gamma}(0,0)$ and \eqref{e5.6}, \eqref{e5.7} hold.
\end{proof}


By an approximation (see \cite[Section 5]{MR3918407}), we have
\begin{theorem}[\textbf{$C^{1,\alpha}$ regularity}]\label{th3.1}
Let $u$ be a viscosity solution of
\begin{equation*}
\left\{\begin{aligned}
    P_{a}u&=0&& \mbox{in } Q_1^+ \\
u&=0&&  \mbox{on } S_1,
\end{aligned}\right.
\end{equation*}
where $a_n=0$. Then $u\in C_{\gamma}^{1,\bar\alpha}(0,0)$, i.e., there exists $A\in \mathbb{R}$ such that
\begin{equation*}
	|u(x,t)-Ax_n|\leq C (|x|^{1+\bar\alpha}+|t|^{\frac{1+\bar\alpha}{2-\bar\alpha\gamma}}) \quad \text{for all }(x,t)\in Q_{1/2}^+,
\end{equation*}
and
\begin{equation*}
|A|\leq C,
\end{equation*}
where $C>0$ depends only on $n$, $p$, $\gamma$, $|a|$ and $\|u\|_{L^{\infty}(Q_1^+)}$.
\end{theorem}

%
%
\section{Boundary and global \texorpdfstring{$C^{1,\alpha}$}{C1,a} regularity for general problems}\label{S7}
In this section, we prove the boundary $C^{1,\alpha}$ regularity on a general domain by the perturbation technique. Throughout this section, we assume that $(0,0) \in \partial_p \Omega$ and prove the boundary pointwise $C^{1,\alpha}$ regularity at $(0,0)$. If we use $\partial \Omega\in C^{1,\alpha}_{\gamma}(0,0)$ (or $\underset{Q_1}{\mathrm{osc}} \, \partial_p\Omega$), we always assume that \eqref{e-re} and \eqref{e-re2} hold (or \eqref{e-re3} and \eqref{e-re4} hold).

We first prove the following lemma, which provides the ``equicontinuity'' up to the boundary of solutions. This will be used to show the continuity up to the boundary of the limit solution (see the proof of \Cref{L.K1}).
\begin{lemma}\label{L21}
Let $0<\theta<1/4$. Suppose that $u$ is a viscosity solution of
\begin{equation*}
\left\{\begin{aligned}
	P_{a}u &=f&& \mbox{in } \Omega_1 \\
	u&=g&& \mbox{on } (\partial_p \Omega)_1
\end{aligned}\right.
\end{equation*}
with
\begin{align*}
    &a_n=0,\quad |a|\leq 1, \quad
	\|u\|_{L^{\infty}(\Omega_1)}\leq 2, \quad
	\|f\|_{L^{\infty}(\Omega_1)}\leq\theta, \\
	&\|g\|_{L^{\infty}((\partial_p \Omega)_1)}\leq \theta \quad\text{and}\quad
	\underset{Q_1}{\mathrm{osc}} \, \partial_p\Omega \leq \theta.
\end{align*}

Then
\begin{equation}\label{e1.6}
 |u(x,t)|\leq C(x_n+\theta) \quad \text{for all } (x,t) \in \Omega_{1/4},
\end{equation}
where $C>0$ is universal.
\end{lemma}

\begin{proof}
Let
\begin{equation*}
v(x,t)=C\left(1-|x+(1+\theta )e_n|^{-\beta}\right)-t+\theta.
\end{equation*}
As in the proof of \Cref{le2.1}, by taking $\beta$ and $C$ large enough, $v$ satisfies
\begin{equation*}
\left\{\begin{aligned}
	P_av &> \|f\|_{L^{\infty}(\Omega_1)}&& \mbox{in } \Omega\cap Q^{+}_{1}(-\theta e_n,0) \\
	v&\geq  \|g\|_{L^{\infty}((\partial_p \Omega)_1)} && \mbox{on } \partial_p \Omega \cap Q^{+}_{1}(-\theta e_n,0)\\
	v&\geq \|u\|_{L^{\infty}(\Omega_1)} && \mbox{on } \Omega \cap \partial_p Q^{+}_{1}(-\theta e_n,0).
\end{aligned}\right.
\end{equation*}
By the comparison principle and a direct calculation, we have (noting $v(-\theta e_n,0)=\theta$)
\begin{equation*}
	-C(x_n+\theta)\leq -v\leq u\leq v\leq C(x_n+\theta) \quad \mbox{in } \{(x',x_n,t)\in \Omega \cap Q^{+}_{1/2}(-\theta e_n,0)\mid x'=0,\,t=0\}.
\end{equation*}
By considering $v(x'-x'_0,x_n,t-t_0)$  for $(x'_0,0,t_0)\in S_{1/2}$ and similar arguments, we have
\begin{equation*}
	-C(x_n+\theta)\leq u\leq C(x_n+\theta) \quad \mbox{in } \Omega\cap Q^{+}_{1/2}(-\theta e_n,0).
\end{equation*}
Note that $\Omega_{1/4}\subset Q^{+}_{1/2}(-\theta e_n,0)$ and we have \eqref{e1.6}.
\end{proof}

Now, we can prove the key step towards the boundary $C^{1,\alpha}$ regularity.
\begin{lemma}\label{L.K1}
For any $0<\alpha<\bar{\alpha}$ and $0<\eta<1$, there exists $\theta \in (0, \eta]$ depending only on $n$, $p$, $\gamma$, $\alpha$ and $\eta$ such that if $u$ is a viscosity solution of
\begin{equation*}
\left\{\begin{aligned}
	P_au &=f&& \mbox{in } \Omega_1 \\
	u&=g&&  \mbox{on } (\partial_p \Omega)_1
\end{aligned}\right.
\end{equation*}
with
\begin{align*}
    &a_n=0,\quad |a|\leq 1, \quad
\|u\|_{L^{\infty}(\Omega_1)}\leq 2, \quad \|f\|_{L^{\infty}(\Omega_1)}\leq\theta, \\
    &\|g\|_{L^{\infty}((\partial_p \Omega)_1)}\leq \theta \quad\text{and}\quad\underset{Q_1}{\mathrm{osc}} \,\partial_p\Omega \leq \theta,
\end{align*}
then there exists $A\in \mathbb{R}$ such that
\begin{equation*}\label{e.lK1.1}
	\|u-Ax_n\|_{L^{\infty}(\Omega_{\tau}^{\tau^{\alpha}})}\leq \eta\tau^{1+\alpha}
\end{equation*}
and
\begin{equation*}\label{e.lK1.2}
	|A|\leq \bar C,
\end{equation*}
where $\bar C>0$ is universal and $\tau \in (0,1/8)$ is a constant depending only on $n$, $p$, $\gamma$, $\alpha$ and $\eta$.
\end{lemma}

\begin{proof}
We prove the lemma by contradiction. Suppose that the lemma is false, that is, there exist $0<\alpha<\bar{\alpha}$,
$0<\eta<1$ and sequences of $a_m$, $u_m$, $f_m$, $g_m$ and $\Omega_m$ such that
\begin{equation*}
\left\{\begin{aligned}
	P_{a_m}u_m &=f_m&&  \mbox{in } \Omega_m\cap Q_1 \\
	u_m&=g_m&&  \mbox{on } \partial_p \Omega_m\cap Q_1
\end{aligned}\right.
\end{equation*}
with
\begin{equation*}
\begin{aligned}
    & (a_m)_n=0,\quad |a_m|\leq 1, \quad
    \|u_m\|_{L^{\infty}(\Omega_m\cap Q_1)}\leq 2, \\
	&\|f_m\|_{L^{\infty}(\Omega_m\cap Q_1)}\leq \frac{1}{m}, \quad
	\|g_m\|_{L^{\infty}(\partial_p \Omega_m\cap Q_1)}\leq \frac{1}{m} \quad\text{and}\quad
	\underset{Q_1}{\mathrm{osc}}\,\partial_p\Omega_m \leq \frac{1}{m}.\\
  \end{aligned}
\end{equation*}
Furthermore,
\begin{equation}\label{e.key1}
	\|u_m-Ax_n\|_{L^{\infty}(\Omega_m\cap Q_{\tau}^{\tau^{\alpha}})}> \eta\tau^{1+\alpha} \quad \text{for all }|A|\leq \bar{C},
\end{equation}
where $\tau \in (0,1/8)$ and $\bar C>0$ will be specified later.

Up to a subsequence, $a_m\to \bar{a}$ and $|\bar{a}|\leq 1$ for some $\bar{a}$. Clearly, $u_m$ are uniformly bounded. In addition, by the interior $C^{1,\alpha}$ regularity \Cref{th2.2}, $u_m$ are equicontinuous in any compact subset of $Q_1^+$. Hence, there exist a subsequence (denoted by $u_m$ again) and $\bar u$ such that
\begin{equation*}
     u_m \rightarrow \bar u ~\mbox{ locally uniformly in } Q_1^+.
\end{equation*}
Note that $\|f_m\|_{L^{\infty}(\Omega_m\cap Q_1)}\rightarrow 0$. Hence, by the stability of viscosity solutions (see \cite[Theorem 6.1]{MR1443043} or \cite[Proposition 3]{MR2804550}), $\bar u\in C(Q_1^+)$ is a viscosity solution of
\begin{equation*}
\begin{aligned}
	P_{\bar{a}}\bar u=0 \quad\mbox{in } Q_{1}^+.
\end{aligned}
\end{equation*}

Next, by \Cref{L21},
\begin{equation*}
	|u_m(x,t)|\leq C(x_n+1/m) \quad \text{for all }(x,t) \in \Omega_m\cap Q_{1/4}.
\end{equation*}
For any $(x,t)\in Q^+_{1/4}$, by taking $m\rightarrow \infty$, we have
\begin{equation*}
	|\bar u(x,t)|\leq Cx_n,
\end{equation*}
which implies that $\bar u$ is continuous up to $S_{1/4}$ and $\bar u\equiv 0$ on $S_{1/4}$.

By \Cref{th3.1}, there exist constants $\bar{A}$ and $\bar{C}$ such that
\begin{equation*}
	|\bar u(x,t)-\bar{A}x_n|\leq \bar{C} (|x|^{1+\bar\alpha}+|t|^{\frac{1+\bar\alpha}{2-\bar\alpha\gamma}})  \quad \text{for all } (x,t)\in Q_{1/8}^+
\end{equation*}
and
\begin{equation*}
  |\bar{A}|\leq \bar{C},
\end{equation*}
where $\bar{C}$ (fixed from now) is universal. Thus,
\begin{equation*}
 \|\bar u-\bar Ax_n\|_{L^{\infty}(Q_{\tau}^{\tau^{\alpha}+})}
 \leq \bar{C}\left(\tau^{1+\bar\alpha}
 +\tau^{(2-\alpha\gamma)\cdot\frac{1+\bar\alpha}{2-\bar\alpha\gamma}}\right)
 = \bar{C}\tau^{1+\alpha}\left(\tau^{\bar{\alpha}-\alpha}
 +\tau^{\frac{(2-\alpha\gamma)(1+\bar\alpha)}{2-\bar\alpha\gamma}-(1+\alpha)}\right).
\end{equation*}
By taking $\tau$ small enough such that
\begin{equation*}
\bar{C}\left(\tau^{\bar{\alpha}-\alpha}
 +\tau^{\frac{(2-\alpha\gamma)(1+\bar\alpha)}{2-\bar\alpha\gamma}-(1+\alpha)}\right)<\frac{\eta}{2},
\end{equation*}
we have
\begin{equation}\label{e.key3}
  \|\bar u-\bar{A}x_n\|_{L^{\infty}(Q_{\tau}^{\tau^{\alpha}+})}\leq \frac{\eta}{2}\tau^{1+\alpha}.
\end{equation}
By setting $A=\bar{A}$ and letting $m\rightarrow \infty$ in \eqref{e.key1}, we have
\begin{equation*}
    \|\bar u-\bar{A}x_n\|_{L^{\infty}(Q_{\tau}^{\tau^{\alpha}+})}\geq \eta \tau^{1+\alpha},
\end{equation*}
which contradicts \eqref{e.key3}.
\end{proof}

The scaling version of the above lemma reads:
\begin{lemma}\label{L.K1-2}
Let $0<\alpha<\bar{\alpha}$, $0<\eta<1$ and $\theta>0$ be as in \Cref{L.K1}. Let $u$ be a viscosity solution of
\begin{equation*}
\left\{\begin{aligned}
	P_{a}u&=f&&  \mbox{in } \Omega_r^{r^{\alpha}} \\
	u&=g&&\mbox{on } (\partial_p \Omega)_r^{r^{\alpha}}
\end{aligned}\right.
\end{equation*}
for some $0<r\leq 1$. Suppose that
\begin{equation*}
\begin{aligned}
    &a_n=0,\quad |a|\leq r^{\alpha}, \quad \|u\|_{L^{\infty}(\Omega_r^{r^\alpha})}\leq 2r^{1+\alpha}, \quad \|f\|_{L^{\infty}(\Omega_r^{r^\alpha})}\leq\theta,\\
    &\|g\|_{L^{\infty}((\partial_p \Omega)_{r}^{r^{\alpha}})}\leq \theta r^{1+\alpha}
    \quad\text{and}\quad
\underset{Q_r^{r^\alpha}}{\mathrm{osc}} \, \partial_p\Omega \leq \theta r^{1+\alpha}.
\end{aligned}
\end{equation*}
Then there exists a constant $A$ such that
\begin{equation*}
  \|u-Ax_n\|_{L^{\infty}(\Omega_{\tau r}^{(\tau r)^{\alpha}})}\leq \eta(\tau r)^{1+\alpha}
\end{equation*}
and
\begin{equation*}
|A|\leq \bar Cr^{\alpha},
\end{equation*}
where $\tau \in (0,1/8)$ and $\bar C>0$ are as in \Cref{L.K1}.
\end{lemma}
\begin{proof}
Let
\begin{equation*}
	\tilde{x}=\frac{x}{r}, \quad
	\tilde t=\frac{t}{r^{2-\alpha\gamma}} \quad \text{and} \quad
	\tilde u(\tilde x, \tilde t)=\frac{u(x,t)}{r^{1+\alpha}}.
\end{equation*}
Then $\tilde{u}$ is a solution of
\begin{equation*}
\left\{\begin{aligned}
	P_{\tilde{a}}\tilde u &=\tilde f&& \mbox{in } \tilde\Omega_1 \\
	\tilde u &=\tilde g&& \mbox{on } (\partial_p \tilde\Omega)_1,
\end{aligned}\right.
\end{equation*}
where
\begin{equation*}
\tilde{a}=\frac{a}{r^{\alpha}}, \quad
	\tilde{f}(\tilde x, \tilde t)= \frac{f(x,t)}{r^{(1+\gamma)\alpha-1}}, \quad
	\tilde{g}(\tilde x, \tilde t)=\frac{g(x,t)}{r^{1+\alpha}} \quad\text{and}\quad
    \tilde{\Omega}=\left\{(\tilde x, \tilde t): (r \tilde{x},r^{2-\alpha\gamma}\tilde{t})\in \Omega_r^{r^\alpha}\right\}.\\
\end{equation*}
From the assumptions and noting $\alpha\leq 1/(1+\gamma)$, we have
\begin{align*}
    &\tilde{a}_n=0,\quad|\tilde{a}|\leq 1, \quad
	\|\tilde u\|_{L^{\infty}(\tilde \Omega_1)}\leq 2, \quad
	\|\tilde f\|_{L^{\infty}(\tilde \Omega_1)}\leq\theta, \\
	&\|\tilde g\|_{L^{\infty}((\partial_p \tilde \Omega)_1)}\leq \theta \quad\text{and}\quad
	\underset{Q_1}{\mathrm{osc}} \, \partial_p\tilde \Omega \leq \theta.
\end{align*}
Then by \Cref{L.K1}, there exists $A\in \mathbb{R}$ such that
\begin{equation*}
  \|\tilde u-A\tilde{x}_n\|_{L^{\infty}(\tilde \Omega_{\tau}^{\tau^{\alpha}})}\leq \eta\tau^{1+\alpha}
\end{equation*}
and
\begin{equation*}
|A|\leq \bar C.
\end{equation*}
By taking the rescaling back to $u$, we arrive at the conclusion.
\end{proof}

%

Now, we can prove our main result, i.e., the boundary pointwise $C^{1,\alpha}$ regularity on a general boundary. We present the regularity for $\gamma\geq 0$ and $\gamma<0$,  respectively.
\begin{lemma}[$C^{1,\alpha}$ regularity for $\gamma\geq 0$]\label{th3.2}
For given $\alpha\in(0,\bar{\alpha})$, let $\eta\in (0,1)$ be the constant as in \Cref{L.K2-2} (where we choose $\beta=\alpha$ and $\|(\partial_p \Omega)_1\|_{C^{1,\beta}(0, 0)}\leq 1$) and let $\theta \in (0,\eta]$ be the constant as in \Cref{L.K1} (with respect to this $\eta$). Let $u$ be a viscosity solution of
\begin{equation*}
\left\{\begin{aligned}
	P_au &=f&&  \mbox{in } \Omega_1 \\
	u&=g&& \mbox{on } (\partial_p\Omega)_1
\end{aligned}\right.
\end{equation*}
with $\gamma\geq 0$. Suppose that
\begin{equation*}
\begin{aligned}
	&a_n=0,\quad |a|\leq 1, \quad \|u\|_{L^{\infty}(\Omega_1)}\leq 1, \quad
	\|f\|_{L^{\infty}(\Omega_1)}\leq\theta
\end{aligned}
\end{equation*}
and
\begin{equation*}
    \|g\|_{L^{\infty}((\partial_p \Omega)_{r}^{r^{\alpha}})}\leq \frac{1}{2}\theta r^{1+\alpha}, \quad
\underset{Q_r^{r^\alpha}}{\mathrm{osc}} \, \partial_p\Omega \leq \frac{1}{2\bar{C}}\theta r^{1+\alpha}\quad \text{for all }r \in (0, 1],
\end{equation*}
where $\bar C>0$ is as in \Cref{L.K1}.

Then $u\in C^{1,\alpha}(0,0)$, i.e., there exists a constant $A$ such that
\begin{equation*}
	|u(x,t)-Ax_n|\leq C(|x|^{1+\alpha}+|t|^{\frac{1+\alpha}{2}}) \quad \text{for all } (x,t)\in \Omega_1
\end{equation*}
and
\begin{equation*}
|A|\leq  C
\end{equation*}
where $C>0$ is a constant depending only on $n$, $p$, $\gamma$ and $\alpha$.
\end{lemma}

\begin{proof}

Let us consider the following iteration: there exists a sequence of constants $A_k$ such that for any $k\geq 0$,
\begin{equation}\label{e3.8}
  \|u-A_kx_n\|_{L^{\infty}(\Omega_{\tau^k}^{\tau ^{k\alpha}})}\leq \tau^{k(1+\alpha)}, \quad |A_k|\leq  \tau^{k\alpha} \quad \text{and} \quad |a|\leq \tau^{k\alpha},
\end{equation}
where  $\tau \in (0,1/8)$ is the constant as in \Cref{L.K1}.

Clearly \eqref{e3.8} holds for $k=0$ by taking $A_0=0$. If \eqref{e3.8} holds for any $k\geq 0$, we obtain immediately that $u\in C^{1,\alpha}_{\gamma}(0,0)$. Indeed, for any $(x,t)\in \Omega_1$, there exists $k\geq 1$ such that
\begin{equation*}
	(x,t)\in \Omega_{\tau^{k-1}}^{\tau ^{(k-1)\alpha}}\setminus \Omega_{\tau^k}^{\tau ^{k\alpha}}.
\end{equation*}
Then
\begin{equation*}
\begin{aligned}
	|u(x,t)| &\leq \|u-A_kx_n\|_{L^{\infty}(\Omega_{\tau^{k-1}}^{\tau ^{(k-1)\alpha}})}+|A_k|\tau^{k-1}\\
	&\leq 2 \tau^{(k-1)(1+\alpha)}
	\leq \frac{2}{\tau^{1+\alpha}} \left(|x|^{1+\alpha}+|t|^{\frac{1+\alpha}{2-\alpha\gamma}}\right).
  \end{aligned}
\end{equation*}
Since $\gamma\geq 0$, we have $C^{1,\alpha}_{\gamma}(0,0)\subset C^{1,\alpha}(0,0)$. Hence, $u\in C^{1,\alpha}(0,0)$.

We next assume that \eqref{e3.8} holds for any $k\leq k_0-1$ but not for $k=k_0$. Then we have
$|a| \leq \tau^{(k_0-1)\alpha}$ and
\begin{equation*}
  \|u\|_{L^{\infty}(\Omega_{\tau^{k_0-1}}^{\tau^{(k_0-1)\alpha}})}
  \leq\tau^{(k_0-1)(1+\alpha)}+|A_{k_0-1}|\tau^{k_0-1}
  \leq 2\tau^{(k_0-1)(1+\alpha)}.
\end{equation*}
By applying \Cref{L.K1-2} (with $r=\tau^{k_0-1}$), there exists $A_{k_0}$ such that
\begin{equation}\label{e10.10}
	\|u-A_{k_0}x_n\|_{L^{\infty}(\Omega_{\tau^{k_0}}^{\tau^{k_0\alpha}})}\leq \eta\tau^{k_0(1+\alpha)}
\end{equation}
and
\begin{equation}\label{e10.11}
|A_{k_0}|\leq \bar{C} \tau^{(k_0-1)\alpha}.
\end{equation}
Since \eqref{e3.8} does not hold for $k_0$, we have either
\begin{equation}\label{e4.1}
|A_{k_0}|>\tau^{k_0\alpha}\quad\mbox{ or }\quad |a|>\tau^{k_0\alpha}.
\end{equation}
Consider the following transformation ($\rho \in (0,\bar C+1)$ to be specified soon):
\begin{equation*}
    r=\tau^{k_0}, \quad \tilde{x}=\frac{x}{r}, \quad \tilde t=\frac{t}{\rho^{-\gamma}r^{2}} \quad \text{and}\quad
    \tilde u(\tilde x,\tilde t)=\frac{u(x,t)-A_{k_0}x_n}{r^{1+\alpha}}.
\end{equation*}
Then $\tilde{u}$ is a solution of
\begin{equation*}
\left\{\begin{aligned}
	P_{\tilde a,\tilde{\nu}}\tilde u &=\tilde f&& \mbox{in } \tilde\Omega_1 \\
	\tilde u&=\tilde g&& \mbox{on } (\partial_p \tilde\Omega)_1,
\end{aligned}\right.
\end{equation*}
where
\begin{equation*}
  \tilde a=\rho^{-1}\left(A_{k_0}e_n+a \right),  \quad \tilde{\nu} =\rho^{-1}r^{\alpha},
\end{equation*}
\begin{equation*}
	\tilde{f}(\tilde x,\tilde t)=\frac{f(x,t)}{\rho^{\gamma}r^{\alpha-1}}, \quad
	\tilde{g}(\tilde x,\tilde t)=\frac{g(x,t){-A_{k_0}x_n}}{r^{1+\alpha}} \quad\text{and}\quad
	\tilde{\Omega}=\left\{(\tilde x,\tilde t): (r\tilde x,\rho^{-\gamma}r^2\tilde t)\in \Omega\right\}.
\end{equation*}
In fact, we choose $\rho=|A_{k_0}e_n+a|$ so that $|\tilde{a}|=1$. Since $a_n=0$, $|a| \leq \tau^{(k_0-1)\alpha}$ and \eqref{e10.11}, \eqref{e4.1} hold, we know
\begin{equation*}
r^{\alpha} \leq \rho=|A_{k_0}e_n+a|\leq (\bar C+1)r^{\alpha},
\end{equation*}
which implies $\tilde{\nu}\leq 1$. Moreover, from this choice of $\rho$, the assumptions, \eqref{e10.10} and \eqref{e10.11}, it is easy to check that (note that $\alpha(1+\gamma)<1$)
\begin{equation*}
\|\tilde u\|_{L^{\infty}(\tilde \Omega_1)}\leq \eta, \quad
	\|\tilde f\|_{L^{\infty}(\tilde \Omega_1)}\leq\theta\leq 1,
\end{equation*}
\begin{equation*}
|\tilde{g}(\tilde x,\tilde t)|\leq \theta \left(|\tilde{x}|^{1+\alpha}+|\tilde{t}|^{\frac{1+\alpha}{2}}\right) \quad \text{for all }(\tilde{x},\tilde{t})\in (\partial_p \tilde{\Omega})_1
\end{equation*}
and
\begin{equation*}
\|(\partial_p \tilde\Omega)_1\|_{C^{1,\alpha}(0, 0)} \leq \theta\leq \eta\leq 1.
\end{equation*}
Hence,
\begin{equation*}
\tilde g(0,0)= |D\tilde{g}(0,0)|=0, \quad	\|\tilde g\|_{C^{1,\alpha}(0,0)}\leq \theta.
\end{equation*}
Then by \Cref{L.K2-2} (with $\varepsilon=0$ and $\beta=\alpha$), $\tilde{u}\in C^{1,\alpha}(0,0)$. By taking the rescaling back to $u$, we arrive at the conclusion.
\end{proof}
\begin{remark}\label{re4.1}
The strategy considering the iteration \eqref{e3.8} is motivated by \cite[Lemma 4.3]{MR4122677} and \cite[Corollary 3.3]{MR4128334}.
\end{remark}

Similarly, we have the following $C^{1,\alpha}$ regularity for $\gamma\leq 0$. Since the proof is similar to that of the above lemma, we omit it.
\begin{lemma}[$C^{1,\alpha}$ regularity for $\gamma\leq 0$]\label{th3.3}
For given $\alpha\in(0,\bar{\alpha})$, let $\eta\in (0,1)$ be the constant as in \Cref{L.K2-2} (where we choose $\beta=\alpha$ and $\|(\partial_p \Omega)_1\|_{C^{1,\beta}(0, 0)}\leq 1$) and let $\theta \in (0,\eta]$ be the constant as in \Cref{L.K1} (with respect to this $\eta$). Let $u$ be a viscosity solution of
\begin{equation*}
\left\{\begin{aligned}
	P_au &=f&&  \mbox{in } \Omega_1 \\
	u&=g&& \mbox{on } (\partial_p\Omega)_1
\end{aligned}\right.
\end{equation*}
with $\gamma\leq 0$. Suppose that
\begin{equation*}
\begin{aligned}
	&a_n=0,\quad |a|\leq 1, \quad \|u\|_{L^{\infty}(\Omega_1)}\leq 1, \quad
	\|f\|_{L^{\infty}(\Omega_1)}\leq\theta
\end{aligned}
\end{equation*}
and
\begin{equation*}
    \|g\|_{L^{\infty}((\partial_p \Omega)_r)}\leq \frac{1}{2}\theta r^{1+\alpha}, \quad
\underset{Q_r}{\mathrm{osc}} \, \partial_p\Omega \leq \frac{1}{2\bar{C}}\theta r^{1+\alpha}\quad \text{for all }r \in (0, 1],
\end{equation*}
where $\bar C>0$ is as in \Cref{L.K1}.

Then $u\in C^{1,\alpha}_{\gamma}(0,0)$, i.e., there exists a constant $A$ such that
\begin{equation*}
	|u(x,t)-Ax_n|\leq C\left(|x|^{1+\alpha}+|t|^{\frac{1+\alpha}{2-\alpha\gamma}}\right) \quad \text{for all } (x,t)\in \Omega_1
\end{equation*}
and
\begin{equation*}
|A|\leq  C
\end{equation*}
where $C>0$ is a constant depending only on $n$, $p$, $\gamma$ and $\alpha$.
\end{lemma}

We are now ready to prove \Cref{th1.1}.
\begin{proof}[Proof of \Cref{th1.1}]
We only consider the case $\gamma\geq 0$ and the proof for $\gamma\leq 0$ is similar and we omit it. In the following proof, we just need to use the two-parameter family of scaling to make some normalization such that the assumptions of \Cref{th3.2} are satisfied. Indeed, since $g\in C^{1,\alpha}_{\gamma}(0,0)$, there exists a linear polynomial $L(x)\coloneqq A+\sum_{i=1}^{n}B_ix_i$ such that
\begin{equation}\label{e6.1}
|g(x,t)-L(x)|\leq \|g\|_{C^{1,\alpha}_{\gamma}(0,0)}
\left(|x|^{1+\alpha}+|t|^{\frac{1+\alpha}{2-\alpha\gamma}}\right)\quad
\mbox{on}~~(\partial_p \Omega)_1.
\end{equation}
Since $\partial_p \Omega\in C^{1,\alpha}_{\gamma}(0,0)$,
\begin{equation}\label{e6.2}
    |x_n|\leq
\|(\partial_p \Omega)_1\|_{C^{1,\alpha}_{\gamma}(0,0)}
\left(|x|^{1+\alpha}+|t|^{\frac{1+\alpha}{2-\alpha\gamma}}\right)
\quad \mbox{on } (\partial_p \Omega)_1.
\end{equation}
We observe that $\bar{u}=u-\bar{L}$ is a viscosity solution of
\begin{equation*}
\left\{\begin{aligned}
    P_a\bar u&=f&& \mbox{in } \Omega_1\\
\bar u&=\bar g&& \mbox{on } (\partial_p \Omega)_1,
\end{aligned}\right.
\end{equation*}
where
\begin{equation*}
a=(B_1,\cdots,B_{n-1},0),\quad \bar{L}(x)=L(x',0) \quad \text{and} \quad \bar{g}=g-\bar{L}.
\end{equation*}
By \eqref{e6.1} and \eqref{e6.2},
\begin{equation*}
\begin{aligned}
    |\bar g(x,t)| &\leq |g(x,t)-\bar L(x)| + |B_n| x_n \\
    &\leq \left(1 + \|(\partial_p \Omega)_1\|_{C^{1,\alpha}_{\gamma}(0,0)}\right) \|g\|_{C^{1,\alpha}_{\gamma}(0,0)}\left(|x|^{1+\alpha}+|t|^{\frac{1+\alpha}{2-\alpha\gamma}}\right)
\quad \mbox{on } (\partial_p \Omega)_1.
\end{aligned}
\end{equation*}
Next, let
\begin{equation}\label{e10.12}
\begin{aligned}
    &\rho_0=\left(\| u\|_{L^{\infty}( \Omega_1)}+\|f\|_{L^{\infty}( \Omega_1)}
+\left(1 + \|(\partial_p \Omega)_1\|_{C^{1,\alpha}_{\gamma}(0,0)}\right)\|g\|_{C^{1,\alpha}_{\gamma}(0,0)}+1\right)^2\theta_0^{-2}, \\
    &r_0=\rho_0^{-1/2},\quad \tilde{x}=\frac{x}{r_0}, \quad \tilde{t}=\frac{t}{\rho_0^{-\gamma} r_0^2}
    \quad\text{and}\quad
    \tilde{u}(\tilde{x},\tilde{t})=\frac{\bar u(x,t)}{\rho_0 r_0},
  \end{aligned}
\end{equation}
where $\theta_0>0$ will be determined soon. Then $\tilde{u}$ is a viscosity solution of
\begin{equation*}
\left\{\begin{aligned}
    P_{\tilde a}\tilde u&=\tilde f&& \mbox{in } \tilde \Omega_1\\
\tilde u&=\tilde g&& \mbox{on } (\partial_p \tilde \Omega)_1,
\end{aligned}\right.
\end{equation*}
where
\begin{equation*}
\tilde{a}=\frac{a}{\rho_0}, \quad
	\tilde{f}(\tilde x,\tilde t)=\frac{f(x,t)}{\rho_0^{3/2+\gamma}}, \quad
	\tilde{g}(\tilde x,\tilde t)=\frac{\bar g(x,t)}{\rho_0^{1/2}}
    \quad\text{and}\quad
	\tilde{\Omega}=\left\{(\tilde x,\tilde t): (\rho_0^{-1/2}\tilde x,\rho_0^{-1-\gamma}\tilde t)\in \Omega\right\}.\\
\end{equation*}

Therefore, by choosing $\theta_0\leq \theta$ ($\theta$ is as in \Cref{th3.2}) small enough, the assumptions of \Cref{th3.2} can be satisfied for $\tilde{u}$. By \Cref{th3.2}, $\tilde{u}\in C^{1,\alpha}(0,0)$ and so $u\in C^{1,\alpha}(0,0)$.
\end{proof}

Finally, by combining the interior $C^{1,\alpha}$ regularity (see \Cref{th2.2}) and the boundary $C^{1,\alpha}$ regularity (see \Cref{th1.1}), we can prove the local $C^{1,\alpha}$ regularity on a general domain.
\begin{proof}[Proof of \Cref{th1.2}] As above, we only consider the case $\gamma\geq 0$. The strategy is to combine the interior regularity with the boundary regularity as usual. Up to a two-parameter transformation (cf. \eqref{e10.12} in the proof of \Cref{th1.1} above), we can assume
\begin{equation*}
    \|u\|_{L^{\infty}(\Omega_1)}\leq1,\quad
    \|f\|_{L^{\infty}(\Omega_1)}\leq1,\quad
    \|g\|_{C^{1,\alpha}_{\gamma}((\partial_p\Omega)_1)} \leq1
    \quad\text{and}\quad
    \|(\partial_p\Omega)_1\|_{C^{1,\alpha}_{\gamma}}\leq 1.
\end{equation*}

By the boundary $C^{1,\alpha}$ regularity \Cref{th1.1}, $u\in C^{1,\alpha}(x,t)$ for any $(x,t)\in \partial_p \Omega\cap Q_{1/2}$ and
\begin{equation*}
\|u\|_{C^{1,\alpha}(x,t)} \leq C,
\end{equation*}
where $C>0$ depends only on $n$, $p$, $\gamma$ and $\alpha$. Then we can consider a transformation like \eqref{e10.12} again such that after the transformation, the new solution $\tilde u$ satisfies
\begin{equation*}
\|\tilde u\|_{C^{1,\alpha}(x,t)} \leq 1 \quad \mbox{for all } (x,t)\in \partial_p \Omega\cap Q_{1/2}.
\end{equation*}
Hence, without loss of generality, we also assume
\begin{equation*}
\|u\|_{C^{1,\alpha}(x,t)} \leq 1 \quad \mbox{for all } (x,t)\in \partial_p \Omega\cap Q_{1/2}
\end{equation*}
throughout this proof.

To prove the conclusion, we only need to show that for any $(x_0,t_0)\in \overline{\Omega_{1/2}}$, there exists a linear polynomial $L$ such that
\begin{equation}\label{e7.1}
|u(x,t)-L(x)|\leq C \left(|x-x_0|^{1+\alpha}+|t-t_0|^{\frac{1+\alpha}{2}}\right) \quad \text{for all }(x,t)\in \Omega\cap Q_{1/4}(x_0,t_0).
\end{equation}

If $(x_0,t_0)\in \partial_p \Omega$, then the estimate \eqref{e7.1} follows from \Cref{th1.1}. In the following, we assume $(x_0,t_0)\in \Omega_{1/2}$. Let
\begin{equation*}
    r_0=\sup\left\{r\in(0,1) \mid B_r(x_0)\times\{t_0\}\subset \Omega\right\}
\end{equation*}
and choose $(x_1,t_0)$ such that
\begin{equation*}
    (x_1,t_0)\in \partial_p \Omega\cap (\partial B_{r_0}(x_0)\times\{t_0\}).
\end{equation*}
Without loss of generality, we assume $r_0<1/8$. By the boundary $C^{1,\alpha}$ regularity at $(x_1,t_0)$, there exists a linear polynomial $L_1$ such that
\begin{equation}\label{e7.2}
|u(x,t)-L_1(x)|\leq |x-x_1|^{1+\alpha}+|t-t_0|^{\frac{1+\alpha}{2}}\quad \text{for all }(x,t)\in \Omega\cap Q_{1/4}(x_1,t_0).
\end{equation}
and
\begin{equation*}
    |DL_1|\leq 1.
\end{equation*}

Next, we prove the conclusion according to two cases:

\noindent\textbf{Case 1:} $|DL_1|\leq r_0^{\alpha}$. Consider the following transformation:
\begin{equation*}
    \tilde{x}=\frac{x-x_0}{r_0}, \quad
    \tilde t=\frac{t-t_0}{r_0^{2-\alpha\gamma}}, \quad
    \tilde a=\frac{DL_1 }{r_0^{\alpha}}
    \quad \text{and}\quad
    \tilde u(\tilde x,\tilde t)=\frac{u(x,t)-L_1(x)}{r_0^{1+\alpha}}.
\end{equation*}
Then $\tilde{u}$ is a solution of
\begin{equation*}
	P_{\tilde a}\tilde u =\tilde f\quad \mbox{in } Q_1,
\end{equation*}
where
\begin{equation*}
	\tilde{f}(\tilde x,\tilde t)=r_0^{1-(1+\gamma)\alpha} f(x,t).
\end{equation*}
By \eqref{e7.2} and noting $(1+\gamma)\alpha<1$,
\begin{equation*}
    |\tilde{a}|\leq 1, \quad 	
    \|\tilde u\|_{L^{\infty}(Q_1)}\leq 1 \quad\text{and}\quad 	
    \|\tilde f\|_{L^{\infty}(Q_1)}\leq1.
\end{equation*}
From the interior $C^{1,\bar\alpha}_{\gamma}$ regularity (see \Cref{th2.2}), there exists a linear polynomial $\tilde{L}$ such that (note that $\gamma\geq 0$ and $\alpha\leq\bar\alpha/2$)
\begin{equation*}
    |\tilde u(\tilde x,\tilde t)-\tilde L(\tilde x)|\leq C \left(|\tilde x|^{1+\bar\alpha}
+|\tilde t|^{\frac{1+\bar\alpha}{2-\bar\alpha\gamma}}\right)
\leq C \left(|\tilde x|^{1+\alpha}
+|\tilde t|^{\frac{1+\alpha}{2}}\right)
 \quad \text{for all }(\tilde x,\tilde t)\in Q_{1/2}.
\end{equation*}
and
\begin{equation*}
|D\tilde L|\leq C,
\end{equation*}
where $C$ is universal. By rescaling back to $u$ and recalling that $Q_{r_0/2}(x_0,t_0) \subset Q_{r_0/2}^{r_0^{\alpha}}(x_0,t_0)$ when $\gamma \geq 0$, we have
 \begin{equation}\label{e7.3}
|u(x,t)-L(x)|\leq C \left(|x-x_0|^{1+\alpha}
+|t-t_0|^{\frac{1+\alpha}{2}}\right) \quad \text{for all }(x,t)\in Q_{r_0/2}(x_0,t_0),
\end{equation}
where
\begin{equation*}
L(x)=L_1(x)+r_0^{1+\alpha}\tilde{L}(\tilde{x}).
\end{equation*}

For any $(x,t)\in \Omega\cap Q_{1/4}(x_0,t_0)$, if $(x,t)\in Q_{r_0/2}(x_0,t_0)$, we obtain \eqref{e7.1} immediately from \eqref{e7.3}. If $(x,t)\in (\Omega\cap Q_{1/4}(x_0,t_0))\setminus Q_{r_0/2}(x_0,t_0)$, we have either
\begin{equation*}
    |x-x_0|\geq \frac{r_0}{2} \quad\mbox{or}\quad |t-t_0|\geq \left(\frac{r_0}{2}\right)^{2}.
\end{equation*}
Then
\begin{equation*}
  \begin{aligned}
    |u(x,t)-L(x)|&\leq |u(x,t)-L_1(x)|+|r_0^{1+\alpha}\tilde{L}(\tilde{x})|\\
    &\leq |x-x_1|^{1+\alpha}+|t-t_0|^{\frac{1+\alpha}{2}}
    +C\left(r_0^{1+\alpha}+r_0^{\alpha}|x-x_0|\right)\\
    &\leq C \left(|x-x_0|^{1+\alpha}+|t-t_0|^{\frac{1+\alpha}{2}}\right).
  \end{aligned}
\end{equation*}

\noindent\textbf{Case 2:} $|DL_1|\geq r_0^{\alpha}$. In this case, there exists $ \beta\in[0, \alpha]$ such that
\begin{equation*}
|DL_1|r_0^{-\beta}=1.
\end{equation*}
Now, we use the following transformation:
\begin{equation*}
    \tilde{x}=\frac{x-x_0}{r_0}, \quad
    \tilde t=\frac{t-t_0}{r_0^{2-\beta\gamma}} \quad\text{and}\quad
    \tilde u(\tilde x,\tilde t)=\frac{u(x,t)-L_1(x)}{r_0^{1+\alpha}}.
\end{equation*}
Then $\tilde{u}$ is a solution of
\begin{equation*}
	P_{\tilde a,\tilde \nu}\tilde u =\tilde f\quad \mbox{in } Q_1,
\end{equation*}
where
\begin{equation*}
\tilde a=r_0^{-\beta} DL_1, \quad
    \tilde \nu=r_0^{\alpha-\beta}\quad\text{and}\quad
    \tilde{f}(\tilde x,\tilde t)=r_0^{1-\alpha-\beta\gamma} f(x,t).
\end{equation*}
Thus,
\begin{equation*}
    |\tilde{a}|=1, \quad \tilde{\nu}\leq 1, \quad 	
    \|\tilde u\|_{L^{\infty}(Q_1)}\leq 1\quad\text{and}\quad 	
    \|\tilde f\|_{L^{\infty}(Q_1)}\leq1.
\end{equation*}
From the interior $C^{1,\bar\alpha}_{\gamma}$ regularity (see \Cref{th2.2}), there exists a linear polynomial $\tilde{L}$ such that
 \begin{equation*}
|\tilde u(\tilde x,\tilde t)-\tilde L(\tilde x)|\leq C \left(|\tilde x|^{1+\alpha}
+|\tilde t|^{\frac{1+\alpha}{2}}\right) \quad \text{for all }(\tilde x,\tilde t)\in Q_{1/2}.
\end{equation*}
and
\begin{equation*}
|D\tilde L|\leq C.
\end{equation*}
By rescaling back to $u$, we have
 \begin{equation*}
|u(x,t)-L(x)|\leq C \left(|x-x_0|^{1+\alpha}
+|t-t_0|^{\frac{1+\alpha}{2}}\right) \quad \text{for all }(x,t)\in Q_{r_0/2}(x_0,t_0),
\end{equation*}
where
\begin{equation*}
L(x)=L_1(x)+r_0^{1+\alpha}\tilde{L}(\tilde{x}).
\end{equation*}

Then as in \textbf{Case 1}, for any $(x,t)\in \Omega\cap Q_{1/4}(x_0,t_0)$, we have
\begin{equation*}
  \begin{aligned}
|u(x,t)-L(x)|\leq C \left(|x-x_0|^{1+\alpha}+|t-t_0|^{\frac{1+\alpha}{2}}\right).
  \end{aligned}
\end{equation*}
\end{proof}
%
%
\begin{appendix}
    \section{Small perturbation regularity}\label{S9}
In the proof of the interior $C^{1,\alpha}$ regularity for viscosity solutions (see \Cref{th2.1}), we need the following small perturbation regularity. See \cite{Sav07} and \cite{Wan13} for similar results in the elliptic and parabolic settings, respectively.
\begin{lemma}\label{L.K2-1}
Let $u$ be a viscosity solution of \eqref{eq-approx} with
\begin{equation*}
1/4\leq |a| \leq 2 \quad\text{and}\quad  0\leq \varepsilon\leq 1.
\end{equation*}
Given $\beta \in (0,1)$, suppose that $0\leq \nu\leq \eta$, where $\eta \in (0,1)$ depends only on $n$, $p$, $\gamma$, $\beta$, $\|u\|_{L^{\infty}(Q_1)}$ and $\|f\|_{L^{\infty}(Q_1)}$.

Then $u\in C^{1,\beta}(0,0)$, i.e., there exists a linear polynomial $L$ such that
\begin{equation*}
	 |u(x,t)-L(x)|\leq C(|x|^{1+\beta}+|t|^{\frac{1+\beta}{2}})\quad \text{for all }(x,t)\in Q_1
\end{equation*}
and
\begin{equation*}
|DL|\leq C,
\end{equation*}
where $C>0$ is a constant depending only on $n$, $p$, $\gamma$, $\beta$, $\|u\|_{L^{\infty}(Q_1)}$ and $\|f\|_{L^{\infty}(Q_1)}$.
\end{lemma}

On the other hand, in the proofs of \Cref{le3.7} and \Cref{th3.2}, we need the following \emph{boundary} version of small perturbation regularity.
\begin{lemma}\label{L.K2-2}
Let $u$ be a viscosity solution of
\begin{equation}\label{e5.1}
\left\{\begin{aligned}
    P^{\varepsilon}_{a, \nu}u&=f&& \mbox{in } \Omega_1 \\
	u&=g&& \mbox{on } (\partial_p \Omega)_1,
\end{aligned}\right.
\end{equation}
where
\begin{equation}\label{cond_lem_a2}
    0\leq\nu\leq 1, \quad
    1/2\leq |a| \leq 2 \quad\text{and}\quad
    0\leq \varepsilon\leq 1.
\end{equation}
Given $\beta\in(0,1)$, suppose that $(\partial_p \Omega)_1\in C^{1,\beta}(0,0)$ and
\begin{equation*}
    \|u\|_{L^{\infty}(\Omega_1)}\leq \eta, \quad
    \|f\|_{L^{\infty}(\Omega_1)}\leq 1, \quad g(0,0)=|Dg(0,0)|=0 \quad\text{and}\quad
    \|g\|_{C^{1,\beta}(0,0)}\leq \eta,
\end{equation*}
where $\eta \in (0,1)$ is a constant depending only on $n$, $p$, $\gamma$, $\beta$ and $\|(\partial_p \Omega)_1\|_{C^{1,\beta}(0, 0)}$.

Then $u\in C^{1,\beta}(0,0)$, i.e., there exists a linear polynomial $L$ such that
\begin{equation*}
	 |u(x,t)-L(x)|\leq C(|x|^{1+\beta}+|t|^{\frac{1+\beta}{2}})\quad \text{for all }(x,t)\in \Omega_1
\end{equation*}
and
\begin{equation*}
|DL|\leq C\eta,
\end{equation*}
where $C>0$ is a constant depending only on $n$, $p$, $\gamma$ and $\beta$.
\end{lemma}

\begin{remark}\label{re10.2}
The small perturbation regularity is to regard the equation as a perturbation of the heat equation (see \eqref{e5.2}). Hence, this regularity falls into the framework of uniformly parabolic equation. Therefore, if the prescribed data are smooth, then we can obtain higher regularity (see \cite{lian2020pointwise,LZ_Parabolic} and \cite{LZ_locally} for related techniques).
\end{remark}

Since the proofs of \Cref{L.K2-1} and \Cref{L.K2-2} are similar and the proof of \Cref{L.K2-1} is simpler, we only give the proof of \Cref{L.K2-2}. First, we prove the key step by the compactness method.
\begin{lemma}\label{le3.1}
Let $0<\beta<1$, $(\partial_p \Omega)_1\in C^{1,\beta}(0,0)$ and $u\in C(\overline{\Omega_1})$ be a viscosity solution of \eqref{e5.1} with \eqref{cond_lem_a2}. Let $r\leq \delta_0 $ and assume that for some $A_0 \in \mathbb{R}$,
\begin{equation*}
	\|u-A_0x_n\|_{L^{\infty}(\Omega_{r})}\leq r^{1+\beta},\quad
	|A_0|\leq \tilde{C}\delta_0^{\beta}, \quad
    \|f\|_{L^{\infty}(\Omega_{1})}\leq 1 \quad\text{and}\quad
	\|g\|_{L^{\infty}((\partial_p\Omega)_{r})}\leq \delta_0 r^{1+\beta},
\end{equation*}
where $\tilde{C}>0$ and $\delta_0 \in (0,1)$ are constants depending only on $n$, $p$, $\gamma$, $\beta$ and $\|(\partial_p \Omega)_1\|_{C^{1,\beta}(0,0)}$.

Then there exists $A \in \mathbb{R}$ such that
\begin{equation*}
  \|u-Ax_n\|_{L^{\infty}(\Omega_{\tau r})}\leq  (\tau r) ^{1+\beta}\quad \text{and} \quad |A-A_0|\leq \tilde{C}(\tau r)^{\beta},
\end{equation*}
where $\tau \in (0,1/2)$ is a constant depending only on $n$, $p$, $\gamma$ and $\beta$.
\end{lemma}
\begin{proof}
We prove the lemma by contradiction. Suppose that there exist a constant $K>0$ and sequences
of $u_m$, $\varepsilon_m$, $a_m$, $\nu_m$, $f_m$, $g_m$, $\Omega_m$, $A_m$, $r_m$ such that $0<r_m\leq 1/m$,
\begin{equation*}
0\leq \varepsilon_m\leq 1, \quad
    1/2\leq |a_m| \leq 2, \quad
    0\leq\nu_m\leq 1,
\end{equation*}

\begin{equation*}
\left\{\begin{aligned}
    P^{\varepsilon_m}_{a_m, \nu_m}u_{m} &=f_m&& \mbox{in } (\Omega_{m})_1 \\
	u_m&=g_m&&  \mbox{on } (\partial_p \Omega_m)_{1},
\end{aligned}\right.
\end{equation*}
\begin{equation}\label{e10.3}
	\|u_m-A_mx_n\|_{L^{\infty}(\Omega_m\cap Q_{r_m})}\leq r_m^{1+\beta},\quad |A_m|\leq \frac{\tilde{C}}{m^
\beta}, \quad \|f_m\|_{L^{\infty}(\Omega_m\cap Q_1)}\leq 1,
\end{equation}
and
\begin{equation*}
    \|g_m\|_{L^{\infty}(\partial_p\Omega_m\cap Q_{r_m})}\leq \frac{r_m^{1+\beta}}{m}, \quad
	\|(\partial_p \Omega_m)_1\|_{C^{1,\beta}(0,0)}\leq K.
\end{equation*}
Moreover, for any $A\in \mathbb{R}$ with
\begin{equation*}
	|A-A_m|\leq \tilde{C}(\tau r_m)^{\beta},
\end{equation*}
we have
\begin{equation}\label{e2.2-2}
	\|u_m-Ax_n\|_{L^{\infty}(\Omega_m\cap Q_{\tau r_m})}> (\tau r_m)^{1+\beta},
\end{equation}
where $\tilde{C}$ and $0<\tau<1/2$ are to be specified later.

Let
\begin{equation*}
	\tilde{x}=\frac{x}{r_m},\quad
	\tilde{t}=\frac{t}{r_m^2} \quad \text{and} \quad
	\tilde{u}_m(\tilde{x},\tilde{t})=\frac{u_m(x,t)-A_mx_n}{r_m^{1+\beta}}.
\end{equation*}
Then $\tilde{u}_m$ are viscosity solutions of
\begin{equation*}
\left\{\begin{aligned}
    P^{\varepsilon_m}_{\tilde {a}_m,\tilde {\nu}_m} \tilde u_m &=\tilde f_m&&  \mbox{in } (\tilde \Omega_m)_1 \\
	\tilde u_m &=\tilde g_m&&  \mbox{on } (\partial_p \tilde \Omega_m)_{1},
\end{aligned}\right.
\end{equation*}
where
\begin{equation*}
\tilde{a}_m=a_m+A_me_n, \quad \tilde \nu_m =\nu_mr_m^{\beta}
\end{equation*}
and
\begin{equation*}
    \tilde{f}_m(\tilde{x},\tilde{t})=\frac{f_m(x,t)}{r_m^{\beta-1}}, \quad
    \tilde{g}_m(\tilde{x},\tilde{t})=\frac{g_m(x,t)-A_m x_n}{r_m^{1+\beta} } \quad \text{and} \quad\tilde{\Omega}_m= \{(\tilde{x},\tilde{t}): (r_m \tilde{x}, r_m^2 \tilde{t}) \in \Omega_m\}.
\end{equation*}

First, by \eqref{e10.3} and the definition of $\tilde{u}$, we have
\begin{equation*}
	\|\tilde{u}_m\|_{L^{\infty}((\tilde{\Omega}_m)_1)}\leq 1.
\end{equation*}
Next, it is easy to verify that there exist $\tilde{\varepsilon}\in [0,1]$ and $\tilde{a}\in \mathbb{R}^n$ with $1/2\leq |\tilde{a}|\leq 2$ such that (up to a subsequence and similarly hereinafter)
\begin{equation*}
\varepsilon_m\to \tilde{\varepsilon}, \quad \tilde a_m\to \tilde{a} \quad \text{and} \quad \tilde \nu_m\to 0.
\end{equation*}
Finally, due to $\tilde \nu_m\to 0$, we have $\tilde \nu_m  \leq \nu_2$ for $m$ large enough where $\nu_2$ is the small constant chosen in \Cref{lem-holder-time}. Thus, we can apply the interior H\"{o}lder estimate (\Cref{lem-lipschitz-space} and \Cref{lem-holder-time}) to verify that $\tilde{u}_m$ are equicontinuous. By Arzel\`{a}--Ascoli theorem, there exists $\tilde{u}\in C(Q_1^+)$ such that
$\tilde{u}_m\to \tilde{u}$ in $L^{\infty}(\Omega')$ for any compact subset $\Omega'$ of $Q_1^+$.

Now, we show that $\tilde{u}$ is a viscosity solution of
\begin{equation}\label{e2.5-2}
    P^{\tilde{\varepsilon}}_{\tilde{a},0}\tilde u=0 \quad\mbox{in } Q^+_{1}.
\end{equation}
Note that \eqref{e2.5-2} is a linear uniformly parabolic equation with constant coefficients. Thus, the definition of ``viscosity solution'' is in the classical sense (see \cite[Definition 3.4]{MR1135923}).

Given $(\tilde x_0, \tilde{t}_0)\in Q^+_{1}$ and $\varphi\in C^2$ touching $\tilde{u}$ strictly by above at $(\tilde x_0, \tilde{t}_0)$. Then there exist a sequence of $(\tilde x_m, \tilde{t}_m)\to (\tilde x_0, \tilde{t}_0)$ such that $\varphi+C_m$ touch $\tilde{u}_m$ by above at $(\tilde x_m, \tilde{t}_m)$ and $C_m\to 0$. In addition, for $m$ large enough,
\begin{equation*}
	\frac{1}{4}\leq \frac{|a_m|}{2}\leq |\nu_m r_m^{\beta}D\varphi+A_me_n+a_m|\leq  2|a_m|\leq  4.
\end{equation*}
Then by the definition of viscosity solution (in the sense of \Cref{de2.1}), for $m$ large enough, we have
\begin{equation*}
    P^{\varepsilon_m}_{\tilde a_m,\tilde \nu_m} \varphi(\tilde x_m, \tilde{t}_m)\leq \tilde f_m(\tilde x_m, \tilde{t}_m).
\end{equation*}
By letting $m\to \infty$, we have
\begin{equation*}
P^{\tilde\varepsilon}_{\tilde a,0} \varphi(\tilde x_0, \tilde{t}_0)\leq 0.
\end{equation*}
Hence, $\tilde{u}$ is a subsolution of \eqref{e2.5-2}. Similarly, we can prove that it is a viscosity supersolution
as well. That is, $\tilde{u}$ is a viscosity solution of \eqref{e2.5-2}.

Next, note that
\begin{equation*}
\begin{aligned}
\|\tilde f_m\|_{L^{\infty}((\tilde \Omega_m)_1)} &\leq r_m^{1-\beta}\|f_m\|_{L^{\infty}((\Omega_m)_{r_m})}\leq r_m^{1-\beta}, \quad\\
\|\tilde g_m\|_{L^{\infty}((\partial_p \tilde \Omega_m)_1)}&\leq r_m^{-(1+\beta)}
\left(\|g_m\|_{L^{\infty}((\partial_p \Omega_m)_{r_m})}+|A_m|\|x_n\|_{L^{\infty}((\partial_p \Omega_m)_{r_m})}\right)\leq
 \frac{1}{m}+\frac{\tilde{C}K }{m^{\beta}}, \quad\\
 \underset{Q_1}{\mathrm{osc}}\,\partial_p\tilde\Omega_m &\leq 2Kr_m^{\beta}.
\end{aligned}
\end{equation*}
Then by using a barrier similar to \eqref{e10.7} and an argument similar to the proof of \Cref{L21}, we have
\begin{equation}\label{e10.8}
	 |\tilde u_m(\tilde x,\tilde t)|\leq C(\tilde x_n+\tilde{C}_m) \quad \text{for all  }(\tilde x,\tilde t) \in (\tilde\Omega_m)_{1/4},
\end{equation}
where $C>0$ is universal and $\tilde{C}_m\to 0$ as $m\to \infty$. Letting $m\rightarrow \infty$ in \eqref{e10.8}, we have
\begin{equation*}
  |\tilde u(\tilde x,\tilde t)|\leq C\tilde x_n \quad \text{for all  }(\tilde x,\tilde t) \in Q^+_{1/4},
\end{equation*}
Hence, $\tilde u$ is continuous up to $S_{1/4}$ and $\tilde  u\equiv 0$ on $S_{1/4}$. Therefore, $\tilde{u}$ is a viscosity solution of
\begin{equation}\label{e5.2}
\left\{\begin{aligned}
    P^{\tilde{\varepsilon}}_{\tilde{a},0}\tilde u&=0&& \mbox{in } Q^+_{1/4} \\
	\tilde u&=0&&  \mbox{on } S_{1/4}.
\end{aligned}\right.
\end{equation}

Since \eqref{e5.2} is a linear uniformly parabolic equation with constant coefficients, $\tilde{u}\in C^{\infty}(\overline{Q^+_{1/8}})$. Then
there exists $\tilde{A}\in \mathbb{R}$ such that for any $\tau \in (0,1/8)$,
\begin{equation*}
  \|\tilde{u}-\tilde A\tilde x_n\|_{L^{\infty}(Q^+_{\tau})}\leq
  \hat{C}\tau^{2}\|\tilde{u}\|_{L^{\infty}(Q^+_{1/4})}\leq
  \hat{C}\tau^{2},
\end{equation*}
and
\begin{equation*}
  \|\tilde A\|\leq \hat{C}\|\tilde{u}\|_{L^{\infty}(Q^+_{1/4})}\leq \hat{C},
\end{equation*}
where $\hat{C}>0$ is universal. If we take $\tau$ small and $\tilde{C}$ large so that
\begin{equation}\label{e10.9}
\tau^{\beta}\leq \frac{1}{2}, \quad  \hat{C}\tau^{1-\beta}\leq \frac{1}{2}\quad \text{and} \quad \hat{C}\leq \tilde{C}\tau^{\beta},
\end{equation}
then we have
\begin{equation}\label{e2.6-2}
    \|\tilde{u}-\tilde A\tilde x_n\|_{L^{\infty}(Q^+_{\tau})}\leq \frac{1}{2}\tau^{1+\beta}\quad\text{and}\quad
    \|\tilde A\|\leq \tilde{C}\tau^{\beta}.
\end{equation}

Furthermore, if we let
\begin{equation*}
    B_m=A_m+r_m^{\beta}\tilde A,
\end{equation*}
then
\begin{equation*}
	|B_m-A_m|\leq \tilde{C}(\tau r_m)^{\beta}.
\end{equation*}
Hence, \eqref{e2.2-2} holds for $B_m$:
\begin{equation*}
	\|u_m-B_mx_n\|_{L^{\infty}(\Omega_m\cap Q_{\tau r_m})}> (\tau r_m)^{1+\beta},
\end{equation*}
or equivalently,
\begin{equation*}
	\|\tilde u_m-\tilde A\tilde{x}_n\|_{L^{\infty}(\tilde\Omega_m\cap Q_{\tau})}> \tau^{1+\beta}.
\end{equation*}
By letting $m\to \infty$, we conclude that
\begin{equation*}
\|\tilde u-\tilde A\tilde{x}_n\|_{L^{\infty}(Q^+_{\tau})}\geq \tau^{1+\beta},
\end{equation*}
which contradicts to \eqref{e2.6-2}.
\end{proof}

Now, we prove \Cref{L.K2-2}.
\begin{proof}[Proof of \Cref{L.K2-2}] It is enough to show that there exists a
sequence of $A_k$ ($k\geq -1$) such that for all $k\geq 0$,
\begin{equation}\label{e3.2-1}
    \|u-A_kx_n\|_{L^{\infty }(\Omega_{\tau^k \delta_0})}\leq (\tau^k \delta_0 )^{1+\beta},
\end{equation}
and
\begin{equation}\label{e3.3-1}
    |A_k-A_{k-1}|\leq \tilde{C}(\tau^k \delta_0 )^{\beta},
\end{equation}
where $\tau$, $\delta_0 $ and $\tilde{C}$ are as in \Cref{le3.1}.

We prove the above by induction. For $k=0$, by setting $A_0\equiv A_{-1}\equiv 0$ and choosing $\eta$ small enough (such that $\eta (1+\|(\partial_p \Omega)_1\|_{C^{1,\beta}(0, 0)})\leq \delta_0$), \eqref{e3.2-1} and \eqref{e3.3-1} hold clearly since we have assumed $\|u\|_{L^{\infty}(\Omega_1)}\leq \eta$. Suppose that the conclusion holds for $k\leq k_0$. By \eqref{e3.3-1} and the first inequality in \eqref{e10.9},
\begin{equation*}
    |A_{k_0}|\leq \sum_{i=1}^{k_0}|A_i-A_{i-1}|\leq
\tilde{C}\delta_0 ^{\beta}\frac{\tau^{\beta}}{1-\tau^{\beta}}\leq \tilde{C}\delta_0 ^{\beta}.
\end{equation*}
In addition, since $g(0,0)=|Dg(0,0)|=0$,
\begin{equation*}
\|g\|_{L^{\infty}((\partial_p\Omega)_r)} \leq  \|g\|_{C^{1,\beta}(0,0)}r^{1+\beta}\leq \eta r^{1+\beta}\leq \delta_0r^{1+\beta}.
\end{equation*}

By \Cref{le3.1} (with $r=\tau^{k_0}\delta_0$ and $A_0=A_{k_0}$), there exists $A_{k_0+1}\in \mathbb{R}$ such that \eqref{e3.2-1} and \eqref{e3.3-1} hold for $k=k_0+1$. By induction, the proof is completed.
\end{proof}

\end{appendix}

\bibliographystyle{abbrv}
\bibliography{reference}

\end{document}